\documentclass{article}[12pt]

\usepackage[russian]{babel}
\usepackage[cp1251]{inputenc}
\usepackage{amsmath}
\usepackage{amssymb}
\usepackage{cite,graphicx}
\usepackage{amscd}

\newtheorem{theorem}{Теорема}
\newtheorem{lemma}{Лемма}
\newtheorem{remark}{Замечание}

\def\Re{\mathop{\rm Re}\nolimits}
\def\Im{\mathop{\rm Im}\nolimits}
\begin{document}
\begin{center}

\textbf{S.~R.~Nasyrov \\
GEOMETRIC AND ANALYTIC PROBLEMS \\ OF THE THEORY OF RAMIFIED
COVERINGS OF THE SPHERE }

\end{center}

This is a course of lectures given for students of the Regional
Mathematical Center of the Novosibirsk State University from
October 20 to November 3, 2017.

The course is devoted to some geometric problems of ramified
coverings of the Riemann sphere. A special attention is payed to
compact surfaces of genus one (complex tori). In the first section
we give a short introduction to the theory of elliptic functions.

Section 2 is devoted to one-parametric families of holomorphic and
meromorphic functions. We recall the role of such families on
Loewner's equation in solving some problems of the theory of
univalent functions. Further we deduce a system of ODEs expressing
dependence of critical points of a family of rational functions
from their critical values. This gives an approximate method to
find a conformal mapping of the Riemann sphere onto a given
simply-connected compact Riemann surface over the sphere.
Thereafter a similar problem is solved for elliptic functions
uniformizing complex tori over the Riemann sphere. The
corresponding system of ODEs also contains a differential equation
for the modules of complex tori. In Section 2 we follow to our
papers \cite{nas}, \cite{nas3}, \cite{nas4}.

In Section 3 we apply complex tori in some problems connected to
Pade-Hermit approximations. We study the partition of a 3-sheeted
Riemann surface into sheets induced by some abelian integrals
(so-called Nattall's partition). In the symmetric case we describe
the structure of trajectories of quadratic differentials connected
with the abelian integrals. \medskip

\begin{center}
\textbf{С.~Р.~Насыров \\ ГЕОМЕТРИЧЕСКИЕ И АНАЛИТИЧЕСКИЕ ПРОБЛЕМЫ\\
ТЕОРИИ РАЗВЕТВЛЕННЫХ НАКРЫТИЙ СФЕРЫ}\footnote{Работа выполнена при
поддержке Регионального математического центра, г.~Новосибирск.}
\end{center}




Содержание этой статьи составляет небольшой курс лекций,
прочитанных для слушателей Регионального математического центра
Новосибирского государственного университета с 20 октября по 3
ноября 2017~г.

Курс посвящен некоторым геометрическим аспектам теории разветвленных
накрытий сферы Римана. Особое внимание уделяется компактным
поверхностям рода один (комплексным торам). В связи с этим, в первом
разделе сначала дается  введение в теорию эллиптических
(двоякопериодических мероморфных) функций. 

Второй раздел посвящен однопараметрическим семействам голоморфных и
мероморфных функций. Напоминается роль таких семейств и
дифференциального уравнения Левнера для решения экстремальных задач,
связанных с однолистными функциями. Далее выводится система
дифференциальных уравнений, выражающая зависимость критических точек
семейства рациональных функций от их критических значений. Это дает
возможность нахождения приближенного конформного отображения сферы
Римана на заданную компактную односвязную риманову поверхность,
разветвленно накрывающую сферу. Затем аналогичная задача решается
для эллиптических функций, т.~е. для комплексных торов над сферой
Римана. При этом, в соответствующую систему входит и
дифференциальное уравнение для определения модуля соответствующего
комплексного тора. Материал этого раздела изложен в статьях автора
\cite{nas}, \cite{nas3}, \cite{nas4}.

Третий раздел посвящен применению комплексных торов в задачах,
связанных с аппроксимациями Паде-Эрмита. Изучается разбиение
трехлистной поверхности рода один на листы, индуцированное
некоторыми абелевыми интегралами (так называемое разложение
Наттолла). В симметричном случае исследуется поведение траекторий
квадратичных дифференциалов, связанных с данными абелевыми
интегралами.

\section{Эллиптические функции}

Теория эллиптических функций --- это теория двоякопериодических
мероморфных функций на комплексной плоскости $\mathbb{C}$. В этом
разделе мы излагаем некоторые факты теории эллиптических функций,
которые являются, в основном, классическими и будут необходимы нам
во втором и третьем разделах. Более подробную информацию можно
найти в известных монографиях и учебниках, напр.,  в
\cite{akhiezer,gur_kur}.

\subsection{Определение эллиптической  функции. Комплексные торы.}

Функция $f$, заданная на всей комплексной плоскости называется
двоякопериодической, если существуют два комплексных числа
$\omega_1$ и $\omega_2$, линейно независимых над полем
$\mathbb{R}$, таких что
$$
f(z+\omega_k)=f(z)\quad \forall z\in \mathbb{C}, \ k=1,2.
$$

Отметим, что примеры мероморфных функций с одним периодом можно
построить достаточно просто: это функции $w=\sin z$, $w=\cos z$,
$w=e^z$ и т.~д. Существование непостоянных двоякопериодических
функций не так очевидно. Отметим также следующее утверждение.

\begin{theorem} [Якоби]
На комплексной плоскости у непостоянной мероморфной функции не может
быть трех линейно независимых над $\mathbb{Q}$ периодов. Если таких
периодов два, то они линейно независимы над $\mathbb{R}$
\end{theorem}

Таким образом, случай двух периодов является наиболее интересным.

Очевидно, что если функция является двоякопериодической с периодами
$\omega_1$ и $\omega_2$, то периодами этой функции будут также все
ненулевые элементы решетки
$$\mbox{\boldmath$\omega$}:=\{m_1\omega_1+m_2\omega_2,\, m_1,m_2\in \mathbb{Z}\}.$$
Множество {\boldmath$\omega$} является свободной абелевой группой с
двумя образующими $\omega_1$ и $\omega_2$. Ясно, что в качестве
образующих решетки $\mbox{\boldmath$\omega$}$ можно взять любые
числа вида
$$
\omega'_1=m_{11}\omega_1+m_{12}\omega_2,\quad
\omega'_2=m_{21}\omega_1+m_{22}\omega_2,
$$
где числа $m_{jk}\in \mathbb{Z}$ и матрица, составленная из них,
унимодулярна, т.е. $m_{11}m_{22}-m_{12}m_{21}=\pm 1$.

Заметим также, что у двоякопериодической функции $f$, заданной на
комплексной плоскости, значение $f(z)$ зависит от класса
эквивалентности точки $z$ по модулю решетки
$\mbox{\boldmath$\omega$}$, поэтому можно определить
фактор-отображение
$\widetilde{f}:\mathbb{C}/\mbox{\boldmath$\omega$}\to \mathbb{C}$ по
правилу $\widetilde{f}(z+\mbox{\boldmath$\omega$})=f(z)$. На
фактор-множестве $T:=\mathbb{C}/\mbox{\boldmath$\omega$}$
определяется топология, превращающая $T$ в поверхность, точнее, в
ориентированную компактную поверхность рода один (тор). Кроме того,
с помощью фактор-отображения $\mathbb{C}\to T$ на $T$ можно ввести
комплексную структуру, превращающую его в риманову поверхность
(комплексный тор). При этом, это фактор-отображение является
универсальным накрытием тора. Отметим, что для разных решеток торы,
как правило, неэквивалентны, т.~е. в общем случае их нельзя
отобразить конформным гомеоморфизмом друг на друга.

\begin{remark}
Достаточно просто найти необходимое и достаточное условие
эквивалентности двух комплексных торов $T'$ и $T''$, выражающееся в
терминах их решеток $\mbox{\boldmath$\omega$}'$ и
$\mbox{\boldmath$\omega$}''$.
\end{remark}

Если двоякопериодическая функция $f$ мероморфна, то
фактор-отображение $\widetilde{f}$ является мероморфной функцией
на торе $T$. Таким образом, по-существу, теория эллиптических
функций является теорией мероморфных функций на торе.

Пусть $\omega_1$ и $\omega_2$ --- любые два периода, порождающие
решетку $\mbox{\boldmath$\omega$}$. Рассмотрим любой
параллелограмм $\Pi$, построенный на  векторах $\omega_1$ и
$\omega_2$, отложенных из произвольной точки $\alpha$. Он
называется \textit{параллелограммом периодов} или
\textit{фундаментальным параллелограммом.} Для удобства в
дальнейшем будем считать, что параллелограмм $\Pi$ содержит не все
граничные точки, а ровно по одной из каждого класса
эквивалентности по модулю решетки. Например, можно добавить к
внутренности параллелограмма две смежные стороны параллелограмма,
представляющие собой отрезки, соединяющие точку $\alpha$ с точками
$\alpha+\omega_1$ и $\alpha+\omega_2$ без этих концевых точек
$\alpha+\omega_1$ и $\alpha+\omega_2$.

Поскольку при замене $\omega_2$ на $(-\omega_2)$ решетка не
меняется, в дальнейшем всегда будем считать, что выполнено условие
$\Im (\omega_2/\omega_1)>0$.

\subsection{Простейшие свойства эллиптических функций. Теоремы Лиувилля}

Отметим следующее свойство эллиптических функций.

\begin{theorem}\label{res}
Сумма вычетов любой эллиптической функции, взятая по всем полюсам,
лежащем в параллелограмме периодов $\Pi$, равна нулю.
\end{theorem}

Доказательство основано на том, что если выбрать параллелограмм
периодов ${\Pi}$ так, чтобы его стороны не содержали полюсов
эллиптической функции, то интеграл по границе параллелограмма от
этой функции равен нулю, поскольку интегралы по противоположным
сторонам параллелограмма сокращаются.  С другой стороны, интеграл
равен сумме вычетов по всем полюсам, лежащим в $\Pi$. Наконец, в
каждом параллелограмме содержится ровно по одному полюсу из каждого
класса эквивалентности, а вычет зависит только от этого класса
эквивалентности.

Напомним, что в окрестности своего нуля $z_0$ любая непостоянная
мероморфная функция $f$ имеет разложение
$$
f(z)=\sum_{k=m}^\infty \alpha_k(z-z_0)^k, \quad
$$
где $m\ge 1$ и $\alpha_m\neq 0$. Число $m$ называется кратностью
этого нуля. В общем случае, если $f(z_0)=a$, точка $z_0$ называется
$a$-точкой функции $f$, и кратность этой точки, по определению, это
кратность нуля функции $f(z)-a$.

\begin{theorem}[Лиувилль]\label{val}
$1)$ Любая непостоянная эллиптическая функция $f$ в каждом
параллелограмме периодов имеет одинаковое число нулей и полюсов (с
учетом их кратности). Более того, для любого $a\in \mathbb{C}$
число $a$-точек (с учетом их кратности) совпадает с числом полюсов
функции $f$ и является характеристикой функции $f$.

$2)$ Если эллиптическая функция $f$ не имеет полюсов в
параллелограмме периодов, то она постоянна.

$3)$ Непостоянная эллиптическая функция имеет не менее двух
полюсов (с учетом их кратности) в каждом параллелограмме периодов.

\end{theorem}

\textit{Доказательство.} 1)  Утверждение следует из теоремы
~\ref{res}, если применить ее к функции $f'/f$, т.~е. является
своеобразной версией принципа аргумента. В случае $a$-точек нужно
вместо $f$ рассмотреть функцию $f-a$.

Утверждение 2) сразу следует из 1). Для доказательства 3) заметим,
что если у эллиптической функции ровно один полюс в
параллелограмме периодов, то она принимает любое значение в этом
параллелограмме ровно один раз. Значит, функция должна быть
однолистна в этом параллелограмме. Но тогда полюс должен быть
простым, и вычет функции $f$ в этом полюсе отличен от нуля. Это
противоречит утверждению теоремы~\ref{res}.\hfill
$\square$\bigskip

Пусть $n$ --- число полюсов непостоянной функции $f$ (с учетом их
кратности). Это число называется порядком эллиптической функции
$f$. С точки зрения теории римановых поверхностей
теорема~\ref{val} утверждает, что функция $f$ осуществляет
$n$-кратное разветвленное накрытие сферы Римана. При этом, $n\ge
2$.

Фиксируем точку $a\in \mathbb{C}$.  Пусть $z_1$, $z_2,\ldots,z_n$
--- $a$-точки  функции $f$ в параллелограмме периодов, причем если точка имеет кратность $m$, то она записывается $m$ раз.
Обозначим через $p_1$, $p_2,\ldots,p_n$ полюсы функции $f$ в
параллелограмме периодов, также выписанные с учетом их кратности.

\begin{theorem}[Лиувилль]\label{val1}
Сумма $a$-точек функции $f$ в параллелограмме периодов сравнима с
суммой полюсов по модулю решетки, т.~е.
$$
\sum_{k=1}^n z_k\equiv\sum_{k=1}^n p_k\quad (\mbox{\rm mod}\
\mbox{\boldmath$\omega$}).
$$
\end{theorem}

Доказательство следует из рассмотрения интеграла
$$
\frac{1}{2\pi i}\int_{\partial \Pi}\frac{zf'(z)dz}{f(z)-a}\,,
$$
взятого по границе параллелограмме периодов.

\subsection{${\mathfrak P}$-функция Вейерштрасса}

Наиболее естественный путь для построения эллиптических функций
состоит в использовании двойных рядов. На этом пути Вейерштрасс
построил важнейшую эллиптическую функцию, которая называется
${\mathfrak P}$-функцией Вейерштрасса и традиционно обозначается
готической буквой ${\mathfrak P}$ (<<пэ>>):
\begin{equation}\label{p}
\mathfrak{P}(z)=\frac{1}{z^2}+\sum\nolimits'\left[\frac{1}{(z-\omega)^2}-\frac{1}{\omega^2}\right];
\end{equation}
здесь и далее штрих при знаке суммирования означает, что оно ведется
по всем ненулевым периодам $\omega=m_1\omega_1+m_2\omega_2$ решетки
$\mbox{\boldmath$\omega$}$. Фактически ряд в (\ref{p}) является
двойным, суммирование ведется по всем целым $m_1$, $m_2$ таким, что
$m_1^2+m_2^2\neq 0$.

Сходимость ряда обеспечивается следующим утверждением.

\begin{lemma}\label{series2}
Ряд
\begin{equation}\label{serper}
\sum\nolimits'\frac{1}{|\omega|^p}
\end{equation}
сходится тогда и только тогда, когда $p>2$.
\end{lemma}

\textit{Доказательство.} Сначала покажем, что можно свести дело к
случаю $\omega_1=1$,  $\omega_2=i$,  Функция
$f(x,y)=x\omega_1+y\omega_2$ непрерывна на единичной окружности
$S^1=\{x^2+y^2=1\}$ и по теореме Вейерштрасса ее модуль принимает
там свои наибольшее и наименьшие значения. При этом, наименьшее
значение отлично от нуля, т.~к. $\omega_1$ и $\omega_2$ линейно
независимы над $\mathbb{R}$. Следовательно, существуют константы
$0<c_1<c_2<+\infty$ такие, что $ c_1\le|f(x,y)|\le c_2$, $(x,y)\in
S^1$. В силу однородности $f$ отсюда получаем:
$$
c_1\sqrt{x^2+y^2}\le|f(x,y)|\le c_2\sqrt{x^2+y^2},\quad (x,y)\in
\mathbb{R}^2.
$$
Если применить это неравенство к $\omega=m_1\omega_1+m_2\omega_2$,
то получим
$$
c_1\sqrt{m_1^2+m_2^2}\le |\omega|\le c_2\sqrt{m_1^2+m_2^2}.
$$
Следовательно, ряд (\ref{serper}) сходится тогда и только тогда.
когда сходится ряд
\begin{equation}\label{serper1}
\sum_{m_1^2+m_2^2\neq 0}\frac{1}{(m_1^2+m_2^2)^{p/2}}.
\end{equation}
Наконец, сходимость ряда (\ref{serper1}) эквивалентна сходимости
двойного интеграла
$$
\int\!\!\!\int_{x^2+y^2>1}
\frac{dxdy}{(x^2+y^2)^{p/2}}=\int_0^{2\pi}d\varphi
\int_1^\infty\frac{r\,dr}{r^p}.
$$
Последний интеграл сходится тогда и только тогда, когда
$p>2$.\hfill $\square$\bigskip

Лемма~\ref{series2} обеспечивает абсолютную сходимость ряда в
правой части  (\ref{p}), так как
$$
\frac{1}{(z-\omega)^2}-\frac{1}{\omega^2}=\frac{z(2\omega-z)}{\omega^2(z-\omega)^2}\sim
\frac{2z}{\omega^3}, \quad \omega\to\infty.
$$
Кроме того, ряд сходится равномерно на любом компакте, не
содержащем точек решетки $\mbox{\boldmath$\omega$}$,
следовательно, определяет некоторую голоморфную функцию в
$\mathbb{C}\setminus \mbox{\boldmath$\omega$}$. Очевидно, что в
каждой точек решетки ровно одно слагаемое ряда имеет особенность
--- полюс второго порядка; после удаления этого слагаемого ряд
сходится абсолютно и равномерно и в окрестности этой точки.
Следовательно, (\ref{p}) определяет мероморфную функцию во всей
плоскости.

Покажем, что ${\mathfrak P}$-функция Вейерштрасса действительно
является двоякопериодической с периодами $\omega_1$ и $\omega_2$.
Имеем с учетом абсолютной сходимости ряда
$$
{\mathfrak P}(z)=\sum_{m_1\in
\mathbb{Z}}\frac{1}{(z-m_1\omega_1)^2}-\sum_{m_1\neq
0}\frac{1}{m_1^2\omega_1^2}+\sum_{m_2\neq 0}\left[\sum_{m_1\in
\mathbb{Z}}\frac{1}{(z-m_1\omega_1-m_2\omega_2)^2}-\sum_{m_1\in
\mathbb{Z}}\frac{1}{(m_1\omega_1+m_2\omega_2)^2}\right].
$$
Очевидно, что сумма
\begin{equation}\label{sin}
\sum_{m_1\in
\mathbb{Z}}\frac{1}{(z-m_1\omega_1)^2}=\frac{\pi^2}{\omega_1^2\sin^2\frac{\pi
z}{\omega_1}}
\end{equation}
является периодической с периодом $\omega_1$, поэтому ${\mathfrak
P}$-функция также имеет период $\omega_1$. Аналогично
показывается, что функция ${\mathfrak P}$ имеет период $\omega_2$.

С учетом (\ref{sin}) можно записать функцию ${\mathfrak P} $ в
виде
\begin{equation*}
    {\mathfrak P}(z)=\frac{\pi^2}{\omega_1^2}\,\left[
-\frac{1}{3}+\sin^{-2}\frac{\pi z}{\omega_1}+\sum_{n\neq
0}\left(\sin^{-2}\frac{\pi
(z-n\omega_2)}{\omega_1}-\sin^{-2}\frac{\pi
n\omega_2}{\omega_1}\right)
    \right].
\end{equation*}

Отметим некоторые очевидные и важные свойства функции
Вейерштрасса.

1) Функция ${\mathfrak P}$ --- четная. Это сразу видно из ее
определения.

2) Функция ${\mathfrak P}$ в некотором смысле однородна. Для этого
обозначим ее через ${\mathfrak P}(z; \omega_1,\omega_2)$, явно
указывая ее зависимость от периодов. Легко видеть, что
\begin{equation}\label{homogen}
{\mathfrak P}(\alpha z; \alpha
\omega_1,\alpha\omega_2)=\alpha^{-2}{\mathfrak P}(z;
\omega_1,\omega_2)
\end{equation} для любого $\alpha\neq 0$, поэтому по
сути дела поведение этой функции зависит от отношения периодов
$\omega_2/\omega_1$.

3) Функция ${\mathfrak P}$ имеет единственный полюс второго порядка
в параллелограмме периодов. Следовательно, в силу теоремы Лиувилля
любое значение она принимает два раза (с учетом кратности). Другими
словами, ${\mathfrak P}$-функция Вейерштрасса осуществляет
двулистное разветвленное накрытие тором расширенной комплексной
плоскости.

Изучим, в каких точках производная функции ${\mathfrak P}$
обращается в нуль. Для этого выведем дифференциальное уравнения
для этой функции.

Имеем в окрестности начала координат
\begin{equation*}
\frac{1}{(z-\omega)^2}-\frac{1}{\omega^2}=\frac{1}{\omega^2}\left[\frac{1}{(1-
z/\omega)^2}-1\right]=\sum_{k=1}^\infty\frac{(k+1)z^k}{\omega^{k+2}},
\end{equation*}
поэтому с учетом четности ${\mathfrak P}$-функция получаем
\begin{equation*}
{\mathfrak
P}(z)=\frac{1}{z^2}+\sum_{n=1}^\infty(2n+1)z^{2n}\,\sum\nolimits'\frac{1}{\omega^{2n+2}}\,.
\end{equation*}
В теории эллиптических функций используются обозначения
$$
g_2=60\sum\nolimits'\frac{1}{\omega^{4}},\quad
g_3=140\sum\nolimits'\frac{1}{\omega^{6}}.
$$
Числа $g_2$ и $g_3$ называются инвариантами Вейерштрасса; они
зависят от периодов $\omega_1$ и $\omega_2$ и являются однородными
функциями от них порядка $(-4)$ и $(-6)$ соответственно. С
использованием этих обозначений запишем
$$
{\mathfrak
P}(z)=\frac{1}{z^2}+\frac{g_2}{20}\,z^2+\frac{g_3}{28}\,z^4+\ldots,
$$
$$
{\mathfrak
P}'(z)=-\frac{2}{z^3}+\frac{g_2}{10}\,z+\frac{g_3}{7}\,z^3+\ldots
$$
С учетом этих разложений получим при $z\to 0$:
$$
({\mathfrak
P}'(z))^2=\frac{4}{z^6}-\frac{2g_2}{5z^2}-\frac{4g_3}{7}\,+o(1),
$$
$$
({\mathfrak P}(z))^3=
\frac{1}{z^6}+\frac{3g_2}{20z^2}+\frac{3g_3}{28}+o(1),
$$
поэтому
$$
({\mathfrak P}'(z))^2-4({\mathfrak P}(z))^3+g_2{\mathfrak
P}(z)=-g_3+o(1).
$$
Теперь заметим, что $({\mathfrak P}'(z))^2-4({\mathfrak
P}(z))^3+g_2{\mathfrak P}(z)$ является эллиптической функцией,
которая имеет только устранимые особенности в точках решетки. В
силу теоремы Лиувилля эта функция --- тождественная константа,
равная $(-g_3)$. Итак, доказана

\begin{theorem}\label{difeq}
Функция Вейерштрасса удовлетворяет дифференциальному уравнению
\begin{equation}\label{difeq1}
{\mathfrak P}'(z))^2=4({\mathfrak P}(z))^3-g_2{\mathfrak
P}(z)-g_3.
\end{equation}
\end{theorem}

Теперь найдем нули функции ${\mathfrak P}'(z)$. Так как ${\mathfrak
P}(z)$ --- четная функция, ее производная ${\mathfrak P}'(z)$
нечетна. Тогда, используя этот факт и периодичность, получим
$$-{\mathfrak P}'(\omega_1/2)={\mathfrak
P}'(-\omega_1/2)={\mathfrak P}'(-\omega_1/2+\omega_1)={\mathfrak
P}'(\omega_1/2).$$ Таким образом, ${\mathfrak P}'(\omega_1/2)=0$.
Аналогично доказывается, что ${\mathfrak P}'(\omega_2/2)=0$,
${\mathfrak P}'((\omega_1+\omega_2)/2)=0$. Итак, эллиптическая
функция ${\mathfrak P}'(z)$ имеет в параллелограмме периодов нули
в точках $\omega_1/2$,   $\omega_2/2$ и  $(\omega_1+\omega_2)/2$.
Так как ее порядок равен $3$, все эти нули --- простые.

В силу (\ref{difeq1}) $$({\mathfrak P}'(z))^2=4({\mathfrak
P(z)}-e_1)({\mathfrak P(z)}-e_2)({\mathfrak P(z)}-e_3),
$$
где $e_1$, $e_2$ и $e_3$ --- корни многочлена $4w^3-g_2w-g_3$.
Значит, можно считать, что $$ {\mathfrak P}(\omega_1/2)=e_1,\
{\mathfrak P}(\omega_2/2)=e_2,\ {\mathfrak
P}((\omega_1+\omega_2)/2)=e_3.
$$
Отметим, что значения $e_1$, $e_2$ и $e_3$ попарно различны,
поскольку эллиптическая функция ${\mathfrak P}$ имеет порядок два,
а кратность этих значений также равна двум.

\subsection{$\zeta$-функция Вейерштрасса}

$\zeta$-функция Вейерштрасса вводится равенством
\begin{equation}\label{zeta}
\zeta(z)=\frac{1}{z}+\sum\nolimits'\left[\frac{1}{z-\omega}+\frac{1}{\omega}+\frac{z}{\omega^2}\right].
\end{equation}

Поскольку
$$
\frac{1}{z-\omega}+\frac{1}{\omega}+\frac{z}{\omega^2}=\frac{z^2}{\omega^2(z-\omega)}\,\sim\,-\,
\frac{z^2}{\omega^3}\,,\quad \omega\to\infty,
$$
ряд (\ref{zeta}) сходится абсолютно; он также сходится равномерно на
компактах вне точек решетки и его сумма представляет собой
мероморфную функцию во всей плоскости. В каждом параллелограмме
периодов эта функция имеет единственный простой полюс с вычетом $1$.

Нетрудно показать, что $\zeta$-функция обладает свойством:
$\zeta'(z)=-\mathfrak{P}(z)$, таким образом, производная
$\zeta$-функции является четной эллиптической функцией.
Следовательно, $\zeta$-функция  нечетна. К сожалению, эллиптической
она не является. Это можно вывести их того, что не существует
эллиптической функции с одним простым полюсом в параллелограмме
периодов.

Однако при изменении аргумента $z$ на период значение функции
изменяется на константу:
\begin{equation}\label{periodzeta}
\zeta(z+\omega_k)=\zeta(z)+\eta_k,\quad k=1,2,
\end{equation}
где $\eta_k=2\zeta(\omega_k/2)$. Действительно, $$
(\zeta(z+\omega_k)-\zeta(z))'=-{\mathfrak
P}(z+\omega_k)+{\mathfrak P}(z)=0,
$$
поэтому $\zeta(z+\omega_k)-\zeta(z)\equiv \mbox{\rm const}$, а
значение константы легко определяется подставлением значения
$z=-\omega_k/2$.

\begin{remark}
Можно доказать, что
$$
\zeta((\omega_1+\omega_2)/2)=\zeta(\omega_1/2)+\zeta(\omega_2/2).
$$
\end{remark}

На основании равенства (\ref{periodzeta}) получаем следующее
утверждение.

\begin{theorem}\label{poles2}
Для любых точек $a$, $b$, не сравнимых по модулю решетки
$\mbox{\boldmath$\omega$}$, мероморфная функция
$\zeta(z-a)-\zeta(z-b)$ является эллиптической функцией, имеющий
два простых полюса в каждом параллелограмме периодов.
\end{theorem}

Теперь установим, что
\begin{equation}\label{etaomega}
\eta_1\omega_2-\eta_2\omega_1=2\pi i.
\end{equation}
Для этого выберем параллелограмм $\Pi$, на границе которого нет
точек решетки $\mbox{\boldmath$\omega$}$, с вершинами в точках
$z_0$, $z_1=z_0+\omega_1$, $z_2=z_0+\omega_1+\omega_2$,
$z_3=z_0+\omega_2$. (Здесь мы используем, что
$\Im(\omega_2/\omega_1)>0$, т.~е. при положительном обходе границы
параллелограмма точки встречаются в порядке: $z_0$, $z_1$, $z_2$,
$z_3$.) Внутри него содержится ровно один простой полюс $\omega$
функции $\zeta$ с вычетом, равным $1$, поэтому с помощью теоремы о
вычетах получаем
$$2\pi i= \mbox{\rm res}_{z=\omega}\zeta(z)=\int_{\partial
\Pi}\zeta(t)dt=\int_{z_0}^{z_1}\zeta(t)dt+\int_{z_1}^{z_2}\zeta(t)dt+\int_{z_2}^{z_3}\zeta(t)dt+\int_{z_3}^{z_0}\zeta(t)dt=$$$$=
\int_{z_0}^{z_1}\zeta(t)dt-\int_{z_3}^{z_2}\zeta(t)dt+\int_{z_1}^{z_2}\zeta(t)dt-\int_{z_0}^{z_3}\zeta(t)dt.
$$
Здесь криволинейные интегралы берутся по отрезкам, соединяющим
вершины параллелограмма. Учитывая (\ref{periodzeta}) и делая
замену переменных, получаем
 $$
\int_{z_0}^{z_1}\zeta(t)dt-\int_{z_3}^{z_2}\zeta(t)dt=\int_{z_0}^{z_1}\zeta(t)dt-\int_{z_0+\omega2}^{z_1+\omega_2}\zeta(t)dt=\int_{z_0}^{z_1}(\zeta(t)-\zeta(t+\omega_2))dt=-\eta_2\omega_1.
 $$
Аналогично
$$
\int_{z_1}^{z_2}\zeta(t)dt-\int_{z_0}^{z_3}\zeta(t)dt=\eta_1\omega_2,
$$
что и доказывает соотношение~(\ref{etaomega}).

\subsection{$\sigma$-функция Вейерштрасса}

Эта функция Вейерштрасса определяется равенством
\begin{equation}\label{sigma}
\sigma(z)=z\prod\nolimits'\left(1-\frac{z}{\omega}\,\right)e{\raisebox{3mm}{$\,\frac{z}{\omega}+\frac{z^2}{2\,\omega^2}$}}
\,.
\end{equation}
Сходимость произведения следует из того, что
\begin{equation}\label{log}
\ln\Bigl\{\prod\nolimits'\left(1-\frac{z}{\omega}\,\right)e{\raisebox{3mm}{$\,\frac{z}{\omega}+\frac{z^2}{2\omega^2}$}}\Bigr\}=
\sum\nolimits'\left[\ln\left(1-\frac{z}{\omega}\,\right)+\,\frac{z}{\omega}+\frac{z^2}{2\,\omega^2}\right].
\end{equation}
При достаточно больших $\omega$ точки $1-{z}/{\omega}$ лежат в
сколь угодно малой окрестности единицы, и мы выбираем ветви
логарифмов в (\ref{log}) так, чтобы $\ln 1=0$. При
$\omega\to\infty$
$$
\ln\left(1-\frac{z}{\omega}\,\right)+\,\frac{z}{\omega}+\frac{z^2}{2\,\omega^2}=\frac{z^3}{3\,\omega^3}+\frac{z^4}{4\,\omega^4}+\ldots=\frac{z^3}{3\,\omega^3}+o(\omega^{-3})
$$
равномерно на компактах на плоскости, поэтому ряд (\ref{log})
сходится абсолютно на компактах, не содержащих точек решетки
$\mbox{\boldmath$\omega$}$. Следовательно, на таких компактах
сходится и бесконечное произведение в правой части (\ref{sigma}).
Но поскольку каждый сомножитель произведения является голоморфной
функцией, по принципу максимума произведение сходится на любом
компакте, даже если он и включает точки решетки
$\mbox{\boldmath$\omega$}$.

Из  (\ref{sigma}) сразу следует, что
\begin{equation}\label{lnsigma}
\frac{\sigma'(z)}{\sigma(z)}=\zeta(z).
\end{equation}
Отсюда заключаем, что $\sigma(z)$ -- нечетная функция. Она имеет
простые нули в точках решетки $\mbox{\boldmath$\omega$}$.

Функция $\sigma(z)$, как и $\zeta(z)$, не является эллиптической, и
при добавлении к аргументу периодов $\omega_k$ она изменяется по
определенному закону. Для нахождения этого закона заметим, что в
силу (\ref{periodzeta}) и (\ref{lnsigma})
$$
\frac{\sigma'(z+\omega_1)}{\sigma(z+\omega_1)}=\frac{\sigma'(z)}{\sigma(z)}+\eta_1.
$$
Интегрируя, получаем
$$ \ln\sigma(z+\omega_1)=\ln\sigma(z) +\eta_1
z+c,
$$
значит,
$$ \sigma(z+\omega_1)=\sigma(z) e^{\eta_1
z+c},
$$
Подставляя вместо $z$ точку  $(-\omega_1/2)$, получаем
$$ \sigma(\omega_1/2)=\sigma(-\omega_1/2) e^{-\eta_1
\omega_1/2+c}=-\sigma(\omega_1/2) e^{-\eta_1 \omega_1/2+c}.
$$
Учитывая, что $\sigma(\omega_1/2)\neq 0$,  находим
$$e^c=-e^{\eta_1 \omega_1/2}.
$$
Значит,
\begin{equation}\label{persi1}
    \sigma(z+\omega_1)=-\sigma(z) e^{\eta_1
(z+\omega_1/2)}.
\end{equation}
Аналогично доказывается, что
\begin{equation}\label{persi2}
\sigma(z+\omega_2)=-\sigma(z) e^{\eta_2 (z+\omega_2/2)}.
\end{equation}

\subsection{Представление эллиптических функций в виде
частного двух произведений $\sigma$-функций}

Пусть даны две системы точек $a_1$, $a_2,\ldots,a_n$ и $b_1$,
$b_2,\ldots,b_n$. Будем предполагать, что никакие две точки из
разных систем не эквивалентны по модулю решетки периодов
$\mbox{\boldmath$\omega$}$. Рассмотрим функцию
\begin{equation}\label{prodsig}
f(z)=\frac{\prod_{k=1}^n\sigma(z-a_k)}{\prod_{k=1}^n\sigma(z-b_k)}
\end{equation}
и поставим вопрос: когда эта функция будет эллиптической?

Необходимое условие для этого дает теорема Лиувилля: необходимо,
чтобы \begin{equation}\label{equiv} a_1+a_2+\ldots+a_n\equiv
b_1+b_2+\ldots+b_n\quad (\mbox{\rm mod}\,
\mbox{\boldmath$\omega$}).
\end{equation} Найдем условие, при котором это
сравнение является и достаточным условием. Используя
(\ref{persi1}), получаем
$$
f(z+\omega_1)=\frac{\prod_{k=1}^n\sigma(z+\omega_1-a_k)}{\prod_{k=1}^n\sigma(z+\omega_1-b_k)}=
\frac{(-1)^n\prod_{k=1}^n\sigma(z-a_k)e^{\eta_1(z-a_k+\omega_1/2)}}{(-1)^n\prod_{k=1}^n\sigma(z-b_k)e^{\eta_1(z-b_k+\omega_1/2)}}=$$$$=
f(z)\prod_{k=1}^ne^{\eta_1(b_k-a_k)}=f(z)e^{\eta_1\sum_{k=1}^n(b_k-a_k)}.
$$
Мы видим что если $\sum_{k=1}^n(b_k-a_k)=0$, то
$f(z+\omega_1)=f(z)$. Аналогично показывается, что при выполнении
этого равенства $f(z+\omega_2)=f(z)$.

Таким образом, если заменить сравнение (\ref{equiv}) на равенство
\begin{equation}\label{equal}
a_1+a_2+\ldots+a_n= b_1+b_2+\ldots+b_n,
\end{equation}
то функция (\ref{prodsig}) становится двоякопериодической.

Теперь установим теорему о восстановлении эллиптической функции по
ее нулям и полюсам.

\begin{theorem}\label{zeroes}
Пусть непостоянная эллиптическая функция $g$ в некотором
параллелограмме периодов имеет $n$ нулей в точках $a_1$,
$a_2,\ldots,a_n$ и $n$ полюсов в точках $b_1$, $b_2,\ldots,b_n$
(нули и полюсы выписываются столько раз, какова их кратность).
Обозначим $$\omega=\sum_{k=1}^n(b_k-a_k).$$ Тогда $\omega$
является элементом решетки $\mbox{\boldmath$\omega$}$ и существует
ненулевая константа $C$ такая, что
$$
g(z)=C\,\frac{\prod_{k=1}^n\sigma(z-a_k)}{\prod_{k=1}^n\sigma(z-b_k^*)},
$$
где $b_k^*=b_k$, $1\le k\le n-1$, и $b_n^*=b_n-\omega$.
\end{theorem}

\textit{Доказательство.} То, что $\omega\in
\mbox{\boldmath$\omega$}$, следует из теоремы Лиувилля. Имеем
$$a_1+a_2+\ldots+a_n= b_1^*+b_2^*+\ldots+b_n^*,$$
поэтому функция
$$
\widetilde{g}(z)=\,\frac{\prod_{k=1}^n\sigma(z-a_k)}{\prod_{k=1}^n\sigma(z-b_k^*)},
$$
является эллиптической и имеет нули и полюсы в тех же точках что и
$g$. Частное $\widetilde{g}/g$ является эллиптической функцией без
нулей и полюсов, поэтому это --- константа $C\neq0$.\hfill
$\square$

Дадим некоторые следствия теоремы~\ref{zeroes}.

Рассмотрим функцию $\mathfrak{P}(z)-\mathfrak{P}(a)$,  где $a\not\in
\mbox{\boldmath$\omega$}$. Эта эллиптическая функций имеет нули в
точках $a$ и $(-a)$, а также полюс второго порядка в начале
координат. Согласно теореме~\ref{zeroes} она представима в виде
$$\mathfrak{P}(z)-\mathfrak{P}(a)=C\,\frac{\sigma(z-a)\,\sigma(z+a)}{\sigma^2(z)}\,.$$
Разлагая левую и правую часть в ряд Лорана в окрестности нуля,
получаем: $C=-1/\sigma^2(a)$.

Итак,
$$\mathfrak{P}(z)-\mathfrak{P}(a)=-\,\frac{\sigma(z-a)\,\sigma(z+a)}{\sigma^2(a)\,\sigma^2(z)}\,.$$
Беря логарифмическую производную, получаем
$$
\frac{\mathfrak{P}'(z)}{\mathfrak{P}(z)-\mathfrak{P}(a)}\,=\zeta(z-a)+\zeta(z+a)-2\zeta(z).
$$
Меняя местами $z$ и $a$, находим

$$
\zeta(z-a)-\zeta(z+a)=-2\zeta(a)+\,\frac{\mathfrak{P}'(a)}{\mathfrak{P}(z)-\mathfrak{P}(a)}\,.
$$
Это равенство выражает эллиптическую функцию, имеющую простые
полюсы в точках $\pm a$, через $\mathfrak{P}$-функцию
Вейерштрасса.

\subsection{Зависимость функций Вейерштрасса от периодов}

В дальнейшем нам понадобятся выражения для частных производных
функции $$\ln \sigma(z)=\ln \sigma(z;\omega_1,\omega_2)$$ по
периодам $\omega_1$ и $\omega_2$. Заодно найдем частные
производные по периодам от функций $\zeta(z;\omega_1,\omega_2)$ и
$\mathfrak{P}(z;\omega_1,\omega_2)$. Разумеется, эти производные с
помощью (\ref{p}), (\ref{zeta}) и (\ref{sigma})  можно представить
в явном виде через ряды, однако мы выразим эти производные через
$\mathfrak{P}$- и $\zeta$-функции Вейерштрасса. Излагаемые здесь
результаты получены в \cite{nas4}.

\begin{theorem}\label{defzeta}
Частные производные функции $\zeta(z)=\zeta(z;\omega_1,\omega_2)$
по периодам $\omega_1$ и $\omega_2$ равны
\begin{equation}\label{zetaom1}
\frac{\partial\zeta(z)}{\partial \omega_1}=\frac{1}{2\pi i}\left[
\frac{1}{2}\,\omega_2\mathfrak{P}'(z)+(\omega_2\zeta(z)-\eta_2z)\mathfrak{P}(z)+
\eta_2 \zeta(z)-(\omega_2 g_2/12)z\right],
\end{equation}
\begin{equation} \label{zetaom2}\frac{\partial\zeta(z)}{\partial
\omega_2}=-\frac{1}{2\pi i}\left[
\frac{1}{2}\,\omega_1\mathfrak{P}'(z)+(\omega_1\zeta(z)-\eta_1z)\mathfrak{P}(z)+
\eta_1 \zeta(z)-(\omega_1 g_2/12)z\right]. \end{equation}
\end{theorem}

Доказательство. Запишем (\ref{periodzeta}) в виде
$$\zeta(z+\omega_k;\omega_1,\omega_2)-\zeta(z;\omega_1,\omega_2)=\eta_k:=2\zeta\left({\omega_k}/{2};\omega_1,\omega_2\right), \quad k=1,2.
$$
Продифференцируем эти равенства по $\omega_2$. Имеем
$$\frac{\partial\zeta(z+\omega_1;\omega_1,\omega_2)}{\partial \omega_2}-\frac{\partial\zeta(z;\omega_1,\omega_2)}{\partial \omega_2}=
2\frac{\partial\zeta\left({\omega_1}/{2};\omega_1,\omega_2\right)}{\partial
\omega_2},$$
\begin{multline*}
$$-\mathfrak{P}(z+\omega_2;\omega_1,\omega_2)+\displaystyle
\frac{\partial\zeta(z+\omega_2;\omega_1,\omega_2)}{\partial
\omega_2}-\frac{\partial\zeta(z;\omega_1,\omega_2)}{\partial
\omega_2}=\\=-\mathfrak{P}\left({\omega_2}/{2};\omega_1,\omega_2\right)+
2\frac{\partial\zeta\left({\omega_2}/{2};\omega_1,\omega_2\right)}{\partial
\omega_2}. \end{multline*}
 Если обозначить
\begin{equation}\label{phi}
\Phi(z):=\frac{\partial\zeta(z;\omega_1,\omega_2)}{\partial
\omega_2}\,,
\end{equation}
то получим
$$
\Phi(z+\omega_1)-\Phi(z)= 2\Phi\left({\omega_1}/{2}\right),\quad
\Phi(z+\omega_2)-\Phi(z)= \mathfrak{P}(z)-e_2+
2\Phi\left({\omega_2}/{2}\right),
$$
где $e_2=\mathfrak{P}\left({\omega_2}/{2}\right)$.

Из соотношений
$$
\zeta(z+\omega_1)\mathfrak{P}(z+\omega_1)-\zeta(z)\mathfrak{P}(z)=\eta_1
\mathfrak{P}(z),\quad
(z+\omega_1)\mathfrak{P}(z+\omega_1)-z\mathfrak{P}(z)=\omega_1
\mathfrak{P}(z),
$$
$$
\zeta(z+\omega_2)\mathfrak{P}(z+\omega_2)-\zeta(z)\mathfrak{P}(z)=\eta_2
\mathfrak{P}(z),\quad
(z+\omega_2)\mathfrak{P}(z+\omega_1)-z\mathfrak{P}(z)=\omega_2
\mathfrak{P}(z),
$$
с учетом (\ref{etaomega}) получаем для функции
\begin{equation}\label{psi}
\Psi(z):=\frac{\omega_1\zeta(z)-\eta_1z}{\omega_1\eta_2-\omega_2\eta_1}\,\mathfrak{P}(z)=-\frac{1}{2\pi
i}\,(\omega_1\zeta(z)-\eta_1z)\,\mathfrak{P}(z)
\end{equation}
соотношения
$$ \Psi(z+\omega_1)-\Psi(z)= 0,\quad
\Psi(z+\omega_2)-\Psi(z)= \mathfrak{P}(z),
$$
поэтому функция
\begin{equation}\label{x}
X(z):=\Phi(z)-\Psi(z)
\end{equation}
удовлетворяет условиям
$$
X(z+\omega_1)-X(z)= 2\Phi\left({\omega_1}/{2}\right),\quad
X(z+\omega_2)-X(z)= -e_2+ 2\Phi\left({\omega_2}/{2}\right).
$$
Подберем константы $\alpha$ и $\beta$ так, чтобы выполнялись
равенства
$$
\alpha \bigl(\zeta(z+\omega_1)-\zeta(z)\bigr) +\beta \omega_1=
2\Phi\left({\omega_1}/{2}\right),\quad \alpha
\bigl(\zeta(z+\omega_2)-\zeta(z)\bigr) +\beta \omega_2= -e_2+
2\Phi\left({\omega_2}/{2}\right).
$$
Тогда функция $X(z)-\alpha \zeta(z)-\beta z$ будет являться
двоякопериодической на плоскости с периодами $\omega_1$ и
$\omega_2$. Учитывая, что в параллелограмме периодов эта функция
имеет единственную особенность в точке $0$, получаем
\begin{equation*}
X(z)-\alpha\zeta(z)-\beta
z=\Phi(z)+\frac{\omega_1\zeta(z)-\eta_1z}{2\pi
i}\,\mathfrak{P}(z)-\alpha\zeta(z)-\beta z
= \frac{\omega_1}{2\pi i}\frac{1}{z^3}-\left(\frac{\eta_1}{2\pi
i}+\alpha\right)\frac{1}{z}+O(z),
\end{equation*}
$z\to 0$. Поэтому двоякопериодическая функция
\begin{equation}\label{funct}
 X(z)-\alpha\zeta(z)-\beta
z+\frac{\omega_1}{4\pi i}\mathfrak{P}'(z)
\end{equation}
в параллелограмме периодов имеет единственный простой полюс в
нуле. Поскольку сумма вычетов любой двоякопериодической функции в
точках, лежащих в параллелограмме периодов, равна нулю, заключаем,
что ее вычет в этой точке равен нулю. Таким образом, эта функция
не имеет особых точек, т.~е. является константой $C$, и равенство
нулю вычета в точке $z=0$ дает
\begin{equation}\label{alpha}
\alpha=-\frac{\eta_1 }{2\pi i}\,.
\end{equation}
Из нечетности функции (\ref{funct}) следует, что константа $C$
равна нулю. С учетом (\ref{phi}), (\ref{psi}) и (\ref{x}) получаем
\begin{multline}\label{funct1}
 X(z)-\alpha\zeta(z)-\beta
z+\frac{\omega_1}{4\pi i}\mathfrak{P}'(z)=
\Phi(z)-\Psi(z)-\alpha\zeta(z)-\beta z+\frac{\omega_1}{4\pi
i}\mathfrak{P}'(z)=\\
=\frac{\partial\zeta(z;\omega_1,\omega_2)}{\partial
\omega_2}+\frac{1}{2\pi
i}\,(\omega_1\zeta(z)-\eta_1z)\,\mathfrak{P}(z)-\alpha\zeta(z)-\beta
z+\frac{\omega_1}{4\pi i}\mathfrak{P}'(z)\equiv 0.
\end{multline}
Используя разложения в нуле 
$$
\mathfrak{P}(z)=\frac{1}{z^2}+\frac{g_2}{20}\,z^2+\frac{g_3}{28}\,z^4+\ldots,\quad
\zeta(z)=\frac{1}{z}-\frac{g_2}{60}\,z^3-\frac{g_3}{140}\,z^5+\ldots,
$$
из равенства нулю коэффициента при первой степени $z$ в
(\ref{funct1}) найдем
\begin{equation}\label{beta}
\beta=\frac{\omega_1g_2 }{24\pi i}.
\end{equation}
Из (\ref{funct1}) с учетом (\ref{alpha}) и (\ref{beta}) получаем
(\ref{zetaom2}). Справедливость (\ref{zetaom1}) устанавливается
аналогично. \hfill $\square$ \smallskip

Интегрируя по $z$ полученные в теореме~\ref{defzeta} соотношения,
с учетом равенств
$$\frac{\partial\ln\sigma(0;\omega_1,\omega_2)}{\partial
\omega_2}=\frac{\partial\ln\sigma(0;\omega_1,\omega_2)}{\partial
\omega_1}=0$$ приходим к следующему результату.

\begin{theorem}
Частные производные функции
$\ln\sigma(z)=\ln\sigma(z;\omega_1,\omega_2)$ по периодам равны
\begin{equation}\label{lnsom1}
 \frac{\partial\ln\sigma(z)}{\partial
\omega_1}\,=\phantom{-}\frac{1}{2\pi
i}\left[\frac{1}{2}\,\omega_2\left(\mathfrak{P}(z)-(\zeta(z))^2\right)+\eta_2(z\zeta(z)-1)+\omega_2\frac{g_2}{24}\,z^2\right],
\end{equation}\begin{equation}\label{lnsom2}
 \frac{\partial\ln\sigma(z)}{\partial
\omega_2}\,=-\frac{1}{2\pi
i}\left[\frac{1}{2}\,\omega_1\left(\mathfrak{P}(z)-(\zeta(z))^2\right)+\eta_1(z\zeta(z)-1)-\omega_1\frac{g_2}{24}\,z^2\right].
\end{equation}
\end{theorem}

Дифференцированием по $z$ соотношений (\ref{zetaom1})  и
(\ref{zetaom2}), получаем

\begin{theorem}
Частные производные функции
$\mathfrak{P}(z)=\mathfrak{P}(z;\omega_1,\omega_2)$ по периодам
равны
\begin{equation}\label{p1}
 \frac{\partial\mathfrak{P}(z)}{\partial
\omega_1}\,=-\frac{1}{2\pi
i}\left[2\omega_2\mathfrak{P}^2(z)-2\eta_2\mathfrak{P}(z)+(\omega_2\zeta(z)-\eta_2z)\mathfrak{P}'(z)-\omega_2\frac{g_2}{3}\right],
\end{equation}\begin{equation}\label{p2}
\frac{\partial\mathfrak{P}(z)}{\partial
\omega_2}\,=\phantom{-}\frac{1}{2\pi
i}\left[2\omega_1\mathfrak{P}^2(z)-2\eta_1\mathfrak{P}(z)+(\omega_1\zeta(z)-\eta_1z)\mathfrak{P}'(z)-\omega_1\frac{g_2}{3}\right].
\end{equation}
\end{theorem}

\subsection{Геометрия функций Вейерштрасса в симметричных случаях}

Поставим вопрос: когда функции Вейерштрасса отображают границы
параллелограмма периодов с вершиной в начале координат в
прямолинейные участки? Начнем с $\mathfrak{P}$-функции
Вейерштрасса. Мы уже отмечали однородность (\ref{homogen}) этой
функции. В силу нее достаточно рассмотреть случай, когда период
$\omega_1$ --- положительное вещественное число. Ясно, что если
решетка $\mbox{\boldmath$\omega$}$ является симметричным
относительно вещественной оси множеством, то
$\mathfrak{P}$-функция Вейерштрасса принимает на вещественной оси
вещественные значения, так как если $x\in \mathbb{R}$ и $\omega$
--- вещественный элемент решетки, то
$$\frac{1}{(x-\omega)^2}-\frac{1}{\omega^2}$$
вещественное число, а если $\omega$ не вещественно, то
$\overline{\omega}$ --- также элемент решетки и
$$\left[\frac{1}{(x-\omega)^2}-\frac{1}{\omega^2}\right]+\left[\frac{1}{(x-\overline{\omega})^2}-
\frac{1}{\overline{\omega}^2}\right]=2\Re\left[\frac{1}{(x-\omega)^2}-\frac{1}{\omega^2}\right]$$
вещественное число.

Найдем условия, при которых решетка $\mbox{\boldmath$\omega$}$
симметрична относительно вещественной оси. Ясно, что это будет
тогда и только тогда, когда $\overline{\omega}_2$ принадлежит
решетке, т.~е. существуют $m_1$, $m_2\in \mathbb{Z}$ такие, что
$\overline{\omega}_2=m_1\omega_1+m_2\omega_2$.  Из сравнения
мнимых частей делаем вывод, что $m_2=-1$ и тогда $\Re
\omega_2=(m_1/2)\omega_1$. Поскольку  мы можем менять базис
решетки, не изменяя функцию Вейерштрасса, заменим период
$\omega_2$ на $\omega_2-n \omega_1$, где $n$ --- целая часть числа
$m_1/2$, и тогда дело сводится к рассмотрению двух случаев: $\Re
\omega_2=0$ и $\Re \omega_2=(1/2)\omega_1$.

1) Если $\Re \omega_2=0$, то решетка симметрична не только
относительно действительной, но и относительно мнимой оси. Пусть
$\omega_1=a>0$, $\omega_2=bi$, $b>0$ ,  Следовательно, если мы
рассмотрим фундаментальный прямоугольник $\Pi=\{x+iy\mid 0<x<a,
0<y<b\}$, то функция Вейерштрасса принимает на границе $\Pi$
вещественные значения.

Решетка симметрична также относительно прямых $\{x=a/2\}$ и
$\{y=b/2\}$, поэтому точки этих прямых также отображаются в точки
вещественной оси. Рассмотрим четверть прямоугольника $\Pi$:
$\widetilde{\Pi}=\{x+iy\mid 0<x<a/2, 0<y<b /2\}$. Его граница
отображается в вещественную ось. Заметим, что функция
$\mathfrak{P}$ четная и имеет периоды $\omega_1$ и $\omega_2$,
поэтому она принимает одинаковые значения в точках, симметричных
относительно центра прямоугольника $\Pi$.  Поскольку  функция
$\mathfrak{P}$ принимает каждое значение в фундаментальном
прямоугольнике $\Pi$ два раза (с учетом кратности), делаем вывод,
что она однолистна в $\widetilde{\Pi}$. Следовательно,
$\mathfrak{P}$ отображает $\widetilde{\Pi}$ на полуплоскость.
Поскольку производная
$$\mathfrak{P}'(x)=-\,\frac{2}{x^3}-2\sum\nolimits'\frac{1}{(x-\omega)^3}$$ при достаточно малых
положительных $x$ отрицательна, делаем вывод, что $\mathfrak{P}$
отображает $\widetilde{\Pi}$ на нижнюю полуплоскость (рис.~1).

Точка $0$ переходит в бесконечность, а точки $e_1$, $e_2$ и $e_3$,
соответствующие трем остальным вершинам $\widetilde{\Pi}$,  лежат
на вещественной оси, причем $e_2<e_3<e_1$.

\begin{remark}\label{rect}
С использованием принципа  симметрии Римана-Шварца нетрудно
определить, куда отображает $\mathfrak{P}$-функция Вейерштрасса
прямоугольники $\Pi_1=\{x+iy\mid 0<x<a/2, 0<y<b\}$,
$\Pi_2=\{x+iy\mid 0<x<a, 0<y<b/2\}$ и $\Pi$, а также тор
$\mathbb{C}/\mbox{\boldmath$\omega$}$.
\end{remark}

\begin{remark}\label{rect2}
Используя геометрический смысл производной и принцип симметрии,
можем также установить, что  $\zeta$-функция Вейерштрасса
отображает прямоугольники $\widetilde{\Pi}$, $$\Pi_1=\{x+iy\mid
0<x<a/2,\, 0<y<b\}, \quad \Pi_2=\{x+iy\mid 0<x<a,\, 0<y<b/2\}$$ и
$\Pi$ на многоугольные области, границы которых лежат на прямых,
параллельных координатным осям.
\end{remark}

\hskip 1 cm \unitlength 1 mm
\linethickness{0.4pt}
\ifx\plotpoint\undefined\newsavebox{\plotpoint}\fi 
\begin{picture}(141.5,47.75)(0,0)
\put(4.75,10.25){\vector(1,0){65.5}}
\put(14.25,4){\vector(0,1){43.75}}
\put(14.25,32.25){\line(1,0){37}}
\put(51.25,32.25){\line(0,-1){21.75}}
\put(14.18,21.18){\line(1,0){.9932}}
\put(16.166,21.18){\line(1,0){.9932}}
\put(18.153,21.18){\line(1,0){.9932}}
\put(20.139,21.18){\line(1,0){.9932}}
\put(22.126,21.18){\line(1,0){.9932}}
\put(24.112,21.18){\line(1,0){.9932}}
\put(26.099,21.18){\line(1,0){.9932}}
\put(28.085,21.18){\line(1,0){.9932}}
\put(30.072,21.18){\line(1,0){.9932}}
\put(32.058,21.18){\line(1,0){.9932}}
\put(34.045,21.18){\line(1,0){.9932}}
\put(36.031,21.18){\line(1,0){.9932}}
\put(38.018,21.18){\line(1,0){.9932}}
\put(40.004,21.18){\line(1,0){.9932}}
\put(41.991,21.18){\line(1,0){.9932}}
\put(43.977,21.18){\line(1,0){.9932}}
\put(45.963,21.18){\line(1,0){.9932}}
\put(47.95,21.18){\line(1,0){.9932}}
\put(49.936,21.18){\line(1,0){.9932}}
\put(32.68,10.5){\line(0,1){1.}} \put(32.68,12.5){\line(0,1){1.}}
\put(32.68,14.5){\line(0,1){1.}} \put(32.68,16.5){\line(0,1){1.}}
\put(32.68,18.5){\line(0,1){1.}} \put(32.68,20.5){\line(0,1){1.}}
\put(32.68,22.5){\line(0,1){1.}} \put(32.68,24.5){\line(0,1){1.}}
\put(32.68,26.5){\line(0,1){1.}} \put(32.68,28.5){\line(0,1){1.}}
\put(32.68,30.5){\line(0,1){1.}}
\put(23.75,15.5){\makebox(0,0)[cc]{$\widetilde\Pi$}}
\put(69.75,8.25){\makebox(0,0)[cc]{$x$}}
\put(12.25,46.25){\makebox(0,0)[cc]{$y$}}
\put(11.75,7.5){\makebox(0,0)[cc]{$O$}}
\put(51,7.75){\makebox(0,0)[cc]{$\omega_1$}}
\put(11,32.25){\makebox(0,0)[cc]{$\omega_2$}}
\put(88,20.5){\vector(1,0){53.5}} \put(104.25,20.5){\circle*{1.}}
\put(115.5,20.5){\circle*{1}} \put(126,20.5){\circle*{1.}}
\put(32.75,10.15){\circle*{1.}} \put(14.25,21.15){\circle*{1}}
\put(32.75,21.15){\circle*{1.}}
\put(104.25,25){\makebox(0,0)[cc]{$e_2$}}
\put(115.5,25.){\makebox(0,0)[cc]{$e_3$}}
\put(126.25,25.){\makebox(0,0)[cc]{$e_1$}}
\put(81.75,19.5){\vector(3,-1){.07}}\qbezier(60.5,20.25)(73.875,22.625)(81.75,19.5)
\put(70.25,26.5){\makebox(0,0)[cc]{$\mathfrak{P}$}}
\put(114.25,13.5){\makebox(0,0)[cc]{$\mathfrak{P}(\widetilde\Pi)$}}
\put(69,0.5){\makebox(0,0)[cc]{Рис.~1}}
\end{picture}
\vskip 1 cm

2) Если $\Re \omega_2=(1/2)\omega_1$, то все углы в треугольнике
$\Delta$, построенном на векторах $\omega_1$ и $\omega_2$, ---
острые. Меняя в предыдущих рассуждениях местами $\omega_1$ и
$\omega_2$, получаем, что граница параллелограмма переходит в
прямолинейные участки, если этот треугольник --- правильный.

Опишем геометрию функций $\mathfrak{P}$ и $\zeta$ в этом случае.
Фундаментальный параллелограмм состоит из двух правильных
треугольников, причем в каждом треугольнике функция $\mathfrak{P}$
принимает каждое значение ровно один раз. Значит, $\mathfrak{P}$
однолистна в треугольнике $ABD$.

Найдем, куда отображается этот треугольник функцией
$\mathfrak{P}$. На стороне $AD$ функция $\mathfrak{P}$ принимает
вещественные значения. Из вида функции $\mathfrak{P}$ следует, что
$\mathfrak{P}(x)\to+\infty$, $x\to 0+$. Функция $\mathfrak{P}'$
вещественна и не обращается в нуль на интервале $(0,\omega/2)$,
соответствующем стороне $AD$ треугольника, поэтому
$\mathfrak{P}(x)$ монотонно убывает от $=\infty$ до значения
$e_1=\mathfrak{P}(\omega_1/2)$ на этом интервале. Решетка
$\mbox{\boldmath$\omega$}$ симметрична относительно прямой, на
которой лежит сторона $AB$ треугольника, поэтому  в силу
однородности $\mathfrak{P}$-функции Вейерштрасса получаем, что
образ этой стороны лежит на прямой $\alpha:\Im w=-\sqrt{3}\Re w$.
При этом, поскольку $\mathfrak{P}'(z)$ не обращается  в нуль в
неконцевых точках отрезка $AB$, с использованием геометрического
смысла производной заключаем, что при движении $z$ по отрезку от
$A$ до $D$ точка $w=\mathfrak{P}(z)$ движется по прямой $\alpha$
так, что $\Im w$ возрастает, а $\Re w$ убывает.

\vskip 1 cm

\hskip -1.2 cm
\unitlength 1mm 
\linethickness{0.4pt}
\ifx\plotpoint\undefined\newsavebox{\plotpoint}\fi 
\begin{picture}(176.25,85.75)(0,0)
\put(10.5,28.75){\vector(1,0){88}}
\put(11,27.75){\vector(0,1){56.5}}
\multiput(11.25,29.25)(.03371934605,.05892370572){734}{\line(0,1){.05892370572}}
\put(36,72.5){\line(1,0){50}}
\multiput(61,28.5)(.03373819163,.05937921727){741}{\line(0,1){.05937921727}}
\multiput(36,72.5)(.03371934605,-.05926430518){734}{\line(0,-1){.05926430518}}
\put(90.5,72.75){\makebox(0,0)[cc]{$G$}}
\put(31.5,74){\makebox(0,0)[cc]{$F$}}
\put(61.75,25.75){\makebox(0,0)[cc]{$E$}}
\put(13.5,26.25){\makebox(0,0)[cc]{$A$}}
\put(33.75,48){\makebox(0,0)[cc]{$B$}}
\put(35.75,26.25){\makebox(0,0)[cc]{$D$}}
\put(19.5,52.25){\makebox(0,0)[cc]{$K$}}
\put(51,53.5){\makebox(0,0)[cc]{$H$}}
\put(100,49.5){\vector(1,0){73.25}}
\put(132.5,22.75){\vector(0,1){62.75}}
\put(147.5,49.75){\line(0,1){1.25}}
\put(147.5,51){\line(1,0){23.75}}
\multiput(124,63.75)(.0416667,.0333333){30}{\line(1,0){.0416667}}
\multiput(125.25,64.75)(-.0337243402,.0593841642){341}{\line(0,1){.0593841642}}
\multiput(124.25,34.75)(-.0333333,.0333333){30}{\line(0,1){.0333333}}
\multiput(123.25,35.75)(-.0336391437,-.0596330275){327}{\line(0,-1){.0596330275}}
\multiput(113.75,16.5)(.0336538462,.0592948718){312}{\line(0,1){.0592948718}}
\multiput(124.25,64.25)(-.0337243402,.0571847507){341}{\line(0,1){.0571847507}}
\put(132.5,23.25){\line(0,-1){6.75}}
\put(90.5,26){\makebox(0,0)[cc]{$x$}}
\put(8.25,83.5){\makebox(0,0)[cc]{$y$}}
\put(165.5,55){\makebox(0,0)[cc]{$E(\infty)$}}
\put(121.5,83){\makebox(0,0)[cc]{$E(\infty)$}}
\put(109.5,78){\makebox(0,0)[cc]{$F(\infty)$}}
\put(109.5,24){\makebox(0,0)[cc]{$F(\infty)$}}
\put(121.5,19){\makebox(0,0)[cc]{$A(\infty)$}}
\put(145,53.75){\makebox(0,0)[cc]{$D$}}
\put(128.25,52.25){\makebox(0,0)[cc]{$B$}}
\put(165.25,46.25){\makebox(0,0)[cc]{$A(\infty)$}}
\put(128.5,66.25){\makebox(0,0)[cc]{$H$}}
\put(122.25,37.75){\makebox(0,0)[cc]{$K$}}
\multiput(123.93,35.18)(.0322917,.0526042){16}{\line(0,1){.0526042}}
\multiput(124.963,36.863)(.0322917,.0526042){16}{\line(0,1){.0526042}}
\multiput(125.996,38.546)(.0322917,.0526042){16}{\line(0,1){.0526042}}
\multiput(127.03,40.23)(.0322917,.0526042){16}{\line(0,1){.0526042}}
\multiput(128.063,41.913)(.0322917,.0526042){16}{\line(0,1){.0526042}}
\multiput(129.096,43.596)(.0322917,.0526042){16}{\line(0,1){.0526042}}
\multiput(130.13,45.28)(.0322917,.0526042){16}{\line(0,1){.0526042}}
\multiput(131.163,46.963)(.0322917,.0526042){16}{\line(0,1){.0526042}}
\multiput(132.196,48.646)(.0322917,.0526042){16}{\line(0,1){.0526042}}
\multiput(133.23,50.33)(.0322917,.0526042){16}{\line(0,1){.0526042}}
\multiput(134.263,52.013)(.0322917,.0526042){16}{\line(0,1){.0526042}}
\multiput(135.296,53.696)(.0322917,.0526042){16}{\line(0,1){.0526042}}
\multiput(136.33,55.38)(.0322917,.0526042){16}{\line(0,1){.0526042}}
\multiput(137.363,57.063)(.0322917,.0526042){16}{\line(0,1){.0526042}}
\multiput(138.396,58.746)(.0322917,.0526042){16}{\line(0,1){.0526042}}
\multiput(139.43,60.43)(.0322917,.0526042){16}{\line(0,1){.0526042}}
\multiput(140.463,62.113)(.0322917,.0526042){16}{\line(0,1){.0526042}}
\multiput(141.496,63.796)(.0322917,.0526042){16}{\line(0,1){.0526042}}
\multiput(142.53,65.48)(.0322917,.0526042){16}{\line(0,1){.0526042}}
\multiput(143.563,67.163)(.0322917,.0526042){16}{\line(0,1){.0526042}}
\multiput(144.596,68.846)(.0322917,.0526042){16}{\line(0,1){.0526042}}
\multiput(145.63,70.53)(.0322917,.0526042){16}{\line(0,1){.0526042}}
\multiput(146.663,72.213)(.0322917,.0526042){16}{\line(0,1){.0526042}}
\multiput(147.696,73.896)(.0322917,.0526042){16}{\line(0,1){.0526042}}
\multiput(148.73,75.58)(.0322917,.0526042){16}{\line(0,1){.0526042}}
\multiput(149.763,77.263)(.0322917,.0526042){16}{\line(0,1){.0526042}}
\multiput(150.796,78.946)(.0322917,.0526042){16}{\line(0,1){.0526042}}
\multiput(151.83,80.63)(.0322917,.0526042){16}{\line(0,1){.0526042}}
\multiput(152.863,82.313)(.0322917,.0526042){16}{\line(0,1){.0526042}}
\multiput(153.896,83.996)(.0322917,.0526042){16}{\line(0,1){.0526042}}
\multiput(124.18,63.93)(.0318182,-.053125){16}{\line(0,-1){.053125}}
\multiput(125.198,62.23)(.0318182,-.053125){16}{\line(0,-1){.053125}}
\multiput(126.216,60.53)(.0318182,-.053125){16}{\line(0,-1){.053125}}
\multiput(127.234,58.83)(.0318182,-.053125){16}{\line(0,-1){.053125}}
\multiput(128.252,57.13)(.0318182,-.053125){16}{\line(0,-1){.053125}}
\multiput(129.271,55.43)(.0318182,-.053125){16}{\line(0,-1){.053125}}
\multiput(130.289,53.73)(.0318182,-.053125){16}{\line(0,-1){.053125}}
\multiput(131.307,52.03)(.0318182,-.053125){16}{\line(0,-1){.053125}}
\multiput(132.325,50.33)(.0318182,-.053125){16}{\line(0,-1){.053125}}
\multiput(133.343,48.63)(.0318182,-.053125){16}{\line(0,-1){.053125}}
\multiput(134.362,46.93)(.0318182,-.053125){16}{\line(0,-1){.053125}}
\multiput(135.38,45.23)(.0318182,-.053125){16}{\line(0,-1){.053125}}
\multiput(136.398,43.53)(.0318182,-.053125){16}{\line(0,-1){.053125}}
\multiput(137.416,41.83)(.0318182,-.053125){16}{\line(0,-1){.053125}}
\multiput(138.434,40.13)(.0318182,-.053125){16}{\line(0,-1){.053125}}
\multiput(139.452,38.43)(.0318182,-.053125){16}{\line(0,-1){.053125}}
\multiput(140.471,36.73)(.0318182,-.053125){16}{\line(0,-1){.053125}}
\multiput(141.489,35.03)(.0318182,-.053125){16}{\line(0,-1){.053125}}
\multiput(142.507,33.33)(.0318182,-.053125){16}{\line(0,-1){.053125}}
\multiput(143.525,31.63)(.0318182,-.053125){16}{\line(0,-1){.053125}}
\multiput(144.543,29.93)(.0318182,-.053125){16}{\line(0,-1){.053125}}
\multiput(145.562,28.23)(.0318182,-.053125){16}{\line(0,-1){.053125}}
\multiput(146.58,26.53)(.0318182,-.053125){16}{\line(0,-1){.053125}}
\multiput(147.598,24.83)(.0318182,-.053125){16}{\line(0,-1){.053125}}
\multiput(148.616,23.13)(.0318182,-.053125){16}{\line(0,-1){.053125}}
\multiput(149.634,21.43)(.0318182,-.053125){16}{\line(0,-1){.053125}}
\multiput(150.652,19.73)(.0318182,-.053125){16}{\line(0,-1){.053125}}
\multiput(151.671,18.03)(.0318182,-.053125){16}{\line(0,-1){.053125}}
\multiput(10.93,28.43)(.0555556,.032963){15}{\line(1,0){.0555556}}
\multiput(12.596,29.419)(.0555556,.032963){15}{\line(1,0){.0555556}}
\multiput(14.263,30.407)(.0555556,.032963){15}{\line(1,0){.0555556}}
\multiput(15.93,31.396)(.0555556,.032963){15}{\line(1,0){.0555556}}
\multiput(17.596,32.385)(.0555556,.032963){15}{\line(1,0){.0555556}}
\multiput(19.263,33.374)(.0555556,.032963){15}{\line(1,0){.0555556}}
\multiput(20.93,34.363)(.0555556,.032963){15}{\line(1,0){.0555556}}
\multiput(22.596,35.352)(.0555556,.032963){15}{\line(1,0){.0555556}}
\multiput(24.263,36.341)(.0555556,.032963){15}{\line(1,0){.0555556}}
\multiput(25.93,37.33)(.0555556,.032963){15}{\line(1,0){.0555556}}
\multiput(27.596,38.319)(.0555556,.032963){15}{\line(1,0){.0555556}}
\multiput(29.263,39.307)(.0555556,.032963){15}{\line(1,0){.0555556}}
\multiput(30.93,40.296)(.0555556,.032963){15}{\line(1,0){.0555556}}
\multiput(32.596,41.285)(.0555556,.032963){15}{\line(1,0){.0555556}}
\multiput(34.263,42.274)(.0555556,.032963){15}{\line(1,0){.0555556}}
\multiput(35.93,43.263)(.0555556,.032963){15}{\line(1,0){.0555556}}
\multiput(37.596,44.252)(.0555556,.032963){15}{\line(1,0){.0555556}}
\multiput(39.263,45.241)(.0555556,.032963){15}{\line(1,0){.0555556}}
\multiput(40.93,46.23)(.0555556,.032963){15}{\line(1,0){.0555556}}
\multiput(42.596,47.219)(.0555556,.032963){15}{\line(1,0){.0555556}}
\multiput(44.263,48.207)(.0555556,.032963){15}{\line(1,0){.0555556}}
\multiput(45.93,49.196)(.0555556,.032963){15}{\line(1,0){.0555556}}
\multiput(47.596,50.185)(.0555556,.032963){15}{\line(1,0){.0555556}}
\put(35.93,72.43){\line(0,-1){.9886}}
\put(35.93,70.452){\line(0,-1){.9886}}
\put(35.93,68.475){\line(0,-1){.9886}}
\put(35.93,66.498){\line(0,-1){.9886}}
\put(35.93,64.521){\line(0,-1){.9886}}
\put(35.93,62.543){\line(0,-1){.9886}}
\put(35.93,60.566){\line(0,-1){.9886}}
\put(35.93,58.589){\line(0,-1){.9886}}
\put(35.93,56.612){\line(0,-1){.9886}}
\put(35.93,54.634){\line(0,-1){.9886}}
\put(35.93,52.657){\line(0,-1){.9886}}
\put(35.93,50.68){\line(0,-1){.9886}}
\put(35.93,48.702){\line(0,-1){.9886}}
\put(35.93,46.725){\line(0,-1){.9886}}
\put(35.93,44.748){\line(0,-1){.9886}}
\put(35.93,42.771){\line(0,-1){.9886}}
\put(35.93,40.793){\line(0,-1){.9886}}
\put(35.93,38.816){\line(0,-1){.9886}}
\put(35.93,36.839){\line(0,-1){.9886}}
\put(35.93,34.862){\line(0,-1){.9886}}
\put(35.93,32.884){\line(0,-1){.9886}}
\put(35.93,30.907){\line(0,-1){.9886}}
\multiput(60.93,28.93)(-.057197,.0325758){15}{\line(-1,0){.057197}}
\multiput(59.214,29.907)(-.057197,.0325758){15}{\line(-1,0){.057197}}
\multiput(57.498,30.884)(-.057197,.0325758){15}{\line(-1,0){.057197}}
\multiput(55.782,31.862)(-.057197,.0325758){15}{\line(-1,0){.057197}}
\multiput(54.066,32.839)(-.057197,.0325758){15}{\line(-1,0){.057197}}
\multiput(52.35,33.816)(-.057197,.0325758){15}{\line(-1,0){.057197}}
\multiput(50.634,34.793)(-.057197,.0325758){15}{\line(-1,0){.057197}}
\multiput(48.918,35.771)(-.057197,.0325758){15}{\line(-1,0){.057197}}
\multiput(47.202,36.748)(-.057197,.0325758){15}{\line(-1,0){.057197}}
\multiput(45.487,37.725)(-.057197,.0325758){15}{\line(-1,0){.057197}}
\multiput(43.771,38.702)(-.057197,.0325758){15}{\line(-1,0){.057197}}
\multiput(42.055,39.68)(-.057197,.0325758){15}{\line(-1,0){.057197}}
\multiput(40.339,40.657)(-.057197,.0325758){15}{\line(-1,0){.057197}}
\multiput(38.623,41.634)(-.057197,.0325758){15}{\line(-1,0){.057197}}
\multiput(36.907,42.612)(-.057197,.0325758){15}{\line(-1,0){.057197}}
\multiput(35.191,43.589)(-.057197,.0325758){15}{\line(-1,0){.057197}}
\multiput(33.475,44.566)(-.057197,.0325758){15}{\line(-1,0){.057197}}
\multiput(31.759,45.543)(-.057197,.0325758){15}{\line(-1,0){.057197}}
\multiput(30.043,46.521)(-.057197,.0325758){15}{\line(-1,0){.057197}}
\multiput(28.327,47.498)(-.057197,.0325758){15}{\line(-1,0){.057197}}
\multiput(26.612,48.475)(-.057197,.0325758){15}{\line(-1,0){.057197}}
\multiput(24.896,49.452)(-.057197,.0325758){15}{\line(-1,0){.057197}}
\put(90,8){\makebox(0,0)[cc]{Рис.~2}}
\end{picture}

В силу того, что функция $\mathfrak{P}(z)$ отображает точки,
симметричные относительно действительной оси, в точки, симметричные
относительно действительной оси, с использованием четности и
периодичности этой функции получаем, что
$\overline{\mathfrak{P}(z)}=\mathfrak{P}(\overline{z})=\mathfrak{P}(-\overline{z})=\mathfrak{P}(\omega_1-\overline{z})$.

Поскольку в точках прямой, на которой лежит сторона $BD$
треугольника $ABD$ выполняется равенство $z=\omega_1-\overline{z}$,
заключаем, что образ отрезка $BD$ лежит на вещественной оси. С
учетом того, что угол треугольника  $ABD$ в точке $D$ прямой и
производная $\mathfrak{P}'(z)$ имеет нуль первого порядка в точке
$\omega_1/2$, соответствующей точке $D$, убеждаемся, что при
конформном отображении треугольника $\mathfrak{P}$ угловой точке $D$
соответствует граничная  точка в образе, лежащая на вещественной
оси, угол в которой равен $\pi$. Используя еще раз то, что
производная $\mathfrak{P}'(z)$ не обращается в нуль в  точках
границы треугольника $ABD$, за исключением разве лишь вершин,
убеждаемся, что отрезку $AB$ соответствует при отображении
$\mathfrak{P}$ отрезок действительной оси от некоторой
точки $\xi_0$ до $e_1$. Поскольку точка $\xi_0$  лежит также на
прямой $\alpha$, заключаем, что $\xi_0=0$, а образ треугольника есть
бесконечный угол $\{-\pi/3<\arg w<0\}$.

Применяя несколько раз принцип симметрии, убеждаемся, что образ
треугольника $AEF$ есть плоскость, разрезанная вдоль трех лучей, с
вершинами в точках  $e_1>0$, $e_2=e^{-2\pi i/3}$ и $e_3=e^{2\pi
i/3}$, на продолжении которых содержится начало координат
(рис.~2).

\begin{remark}\label{rect3}
С использованием принципа симметрии можно показать, что
 $\mathfrak{P}$-функция Вейерштрасса отображает правильный
шестиугольник с центром в начале координат, одна из вершин
которого расположена в точке $B$,  на двулистную риманову
поверхность над сферой Римана. При этом, граница шестиугольника
переходит в объединение разрезов вдоль отрезков, соединяющих точки
$e_1$, $e_2$ и $e_3$ с началом координат. Дзета-функция
Вейерштрасса отображает этот правильный шестиугольник на внешность
правильного шестиугольника, гомотетичного исходному.
\end{remark}

\section{Однопараметрические семейства мероморфных функций}

В этом разделе мы рассматриваем однопараметрические семейства
римановых поверхностей, разветвленно накрывающих сферу Римана,
рода $0$ и $1$. Заметим, что семейства римановых поверхностей
постоянно привлекают внимание специалистов в области комплексного
анализа и его приложений, при этом рассматриваются не только
семейства разветвленных накрытий, но и семейства абстрактных
римановых поверхностей (см., напр., статью \cite{chirka},
посвященную голоморфным семействам таких поверхностей). Наше
внимание сосредоточено на описании гладких семейств поверхностей с
заданным типом ветвления и их зависимости от траекторий проекций
точек ветвления на комплексную плоскость. В пп.~\ref{21} и
\ref{22} кратко описывается роль параметрического метода в теории
однолистных функций и уравнения Левнера и Левнера-Куфарева,
которые послужили толчком к рассмотрению однопараметрических
семейств многолистных функций и соответствующих им римановых
поверхностей. Далее рассматриваются однопараметрические семейства
рациональных функций (поверхности рода нуль) и эллиптических
функций (поверхности рода один или комплексные торы).

\subsection{Экстремальные задачи геометрической теории функций. Гипотеза
Бибербаха}\label{21}

Одними из основных в геометрической теории функций являются задачи,
связанные с описанием значений однолистных функций, их производных
или других аналитических объектов, которые принято называть
\textit{функционалами.} Напомним, что мероморфная функция $f$
называется \textit{однолистной} в области $D$, лежащей в расширенной
комплексной плоскости $\overline{\mathbb{C}}$, если она инъективна.
Название связано с тем, что согласно идее Римана непостоянная в
области $D$ мероморфная функция осуществляет конформное отображение
этой области на некоторую риманову поверхность, расположенную над
$\overline{\mathbb{C}}$. Если число листов римановой поверхности
равно единице, то есть над каждой точкой из $\overline{\mathbb{C}}$
располагается не более одной точки поверхности, то это равносильно
инъективности, т.~е. однолистности функции.

В качестве области $D$ часто выступают единичный круг $E:=\{z\in
\mathbb{C}\mid |z|<1\}$, $E^-:=\{z\in \mathbb{C}\mid |z|>1\}$ или
верхняя полуплоскость $H:=\{z\in \mathbb{C}\mid \Im z>0\}$.

Опишем некоторые задачи, связанные с функциями, заданными в круге
$E$. Обычно в $E$ рассматривают голоморфные однолистные функции.
Поскольку одно из простейших преобразований  --- линейное  $w\mapsto
aw+b$, $a\neq0$, которое сохраняет однолистность функции $w=f(z)$ и
имеет очевидное геометрическое описание (суперпозиция гомотетии,
поворота и сдвига), всегда можно таким преобразованием добиться,
чтобы однолистная функция имела в нуле тейлоровское разложение вида
\begin{equation}\label{s}
f(z)=z+a_2z^2+a_3z^3+\ldots
\end{equation}
Множество однолистных голоморфных функций с разложением (\ref{s})
обозначается через $S$; его принято называть \textit{классом}~$S$
(эта терминология используется и при описании других множеств
однолистных функций). В качестве основных функционалов в классе $S$
выступают такие величины как $f(z)$, $f'(z)$, $f''(z)/f'(z)$,
производная Шварца $\{f,z\}:=(f''(z)/f'(z))'-(f''(z)/f'(z))^2$
($z$~--- произвольная фиксированная точка единичного круга).
Представляет интерес исследование функционалов в подклассах класса
$S$, имеющих геометрическое или аналитическое описание. В качестве
примеров таких подклассов можно привести классы функций с
вещественными коэффициентами $a_n$, классы ограниченных функций
$S_M:=\{f\in S\mid |f(z)|<M, z\in E\}$, $M>1$ --- фиксированное
число, классы однолистных функций, отображающих конформно единичный
круг на выпуклые, звездные относительно нуля области и многие
другие.

На классе $S$ достаточно естественной является топология, сходимость
в которой равносильна равномерной сходимости на компактах в~$E$.

Сложность исследования функционалов в классах однолистных функций
состоит в том, что они нелинейны, и простейшие алгебраические
операции нарушают однолистность. Однако если функционал на классе
$S$ (или каком-то его подклассе)  непрерывен относительно топологии
локально равномерной сходимости, а большинство исследуемых
функционалов являются таковыми, то достаточно изучить его на
некотором  плотном подмножестве.

В качестве наиболее интересных функционалов, кроме отмеченных выше,
выступают коэффициенты $a_n$ в разложении~(\ref{s}). По сути дела,
это (с точностью до множителя) производные функции $f$, подсчитанные
в точке~$0$:
$$
a_n=\,\frac{f^{(n)}(0)}{n!}\,.
$$

В 1916 г. Л.~Бибербах поставил следующую проблему. \medskip

\textbf{Гипотеза Бибербаха.}  \textit{Установить, верно ли, что дл
любой функции $f\in S$ ее коэффициенты удовлетворяют неравенству
$|a_n|\le n$ для любого $n\ge1$ $?$}
\medskip

Заметим, что если положить в (\ref{s}) $a_n=n$, то получается
функция
$$
K_0(z):=z+2z^2+3z^3+\ldots=z(1+z+z^2+z^3+\ldots)'=z\,\left(\frac{1}{1-z}\right)'=\frac{z}{(1-z)^2},
$$
которая называется функцией Кёбе. Функция $K_0$ однолистно
отображает конформно единичный круг $E$ на плоскость с разрезом
вдоль луча, идущего по вещественной оси влево от точки
$K_0(-1)=-1/4$ до $-\infty$. В этом легко убедиться, рассмотрев
функцию $$1/K_0(z)=(z-1)^2/z= z+1/z-2.$$ Эта функция выражается
через функцию Жуковского $w=(1/2)(z+1/z)$, отображающую круг на
внешность отрезка вещественной оси $[-1;1]$.

Рассмотрим  функции
$$
K_\theta(z):=\frac{z}{(1-e^{i \theta}z)^2},
$$
где $\theta$ --- вещественный параметр. Функция $K_\theta(z)$
отображает конформно единичный круг на внешность луча с вершиной в
точке $(-1/4)e^{-i\theta}$, который на своем продолжении содержит
начало координат. Функции
$K_\theta(z)=e^{-i\theta}K_0(e^{i\theta}z)$ называются вращениями
функции Кёбе $K_0(z)$.

Часто гипотеза Бибербаха дополнялась предположением, что равенство
$|a_n|=n$ хотя бы для одного $n\le2$ влечет, что данная однолистная
функция совпадает с одной из функций $K_\theta(z)$.

Гипотеза Бибербаха долгое время являлась путеводной звездой для
специалистов в области геометрической теории функций комплексного
переменного. Она была доказана в 1984~г. американским ученым Луи
де Бранжем \cite{brange} при серьезном содействии советских
математиков, представителей ленинградской школы. В доказательстве
де Бранжа основное место занимал так называемый параметрический
метод, который мы сейчас опишем. поле подобно о параметрическом
методе можно посмотреть в \cite{aleks,goluzin,gutl_ryaz}.

\subsection{Параметрический метод в теории однолистных функций. Уравнения Левнера и
Левнера-Куфарева}\label{22}

Рассмотрим некоторую область $D$, которая получается из расширенной
плоскости проведением разреза вдоль простой гладкой дуги,
соединяющей точки $z_0\in \mathbb{C}$ и $\infty$ и не проходящей
через начало координат. Если рассмотреть множество однолистных
функций отображающих конформно единичный круг  на такие области $D$,
нормированные условиями $f(0)=0$, $f'(0)>0$ (это означает, что
$f'(0)$ --- положительное вещественное число), то его подмножество,
состоящее из функций, для которых $f'(0)=1$, является плотным в $S$.
Обозначим это подмножество через~$\widetilde{S}$.

Функции $f$ класса $\widetilde{S}$ можно описать с помощью
дифференциального уравнения, предложенного Левнером. Пусть $f$
отображает $E$ конформно на $\mathbb{C}\setminus L$, где $L$ ---
дуга, удовлетворяющая свойствам, описанным выше, и $w=w(\tau)$,
$0\le \tau<+\infty$, --- ее уравнение. Тогда $w(0)=z_0$,
$\lim_{t\to+\infty} w(t)=\infty$. Для любого $t\ge0$ обозначим через
$L_t$ поддугу $L$ с уравнением $w=w(\tau)$, $t\le \tau<+\infty$.
Обозначим $D_t=\mathbb{C}\setminus L_t$, $f_t$ --- конформное
отображение $E$ на $D_t$, нормированное так, что $f_t(0)=0$,
$f'_t(0)>0$. В силу леммы Шварца $f'_t$ является строго возрастающей
функцией, причем $\lim_{t\to+\infty}f'(t)=\infty$.
Перепараметризацией можно добиться, чтобы $f'_t(0)=e^tf'(0)$, $t\ge
0$; будем в дальнейшем считать это условие выполненным.

Обозначим $\varphi_t=f_t^{-1}\circ f$; эта функция конформно
отображает единичный круг $E$ на его подобласть, которая получается
из круга $E$ проведением некоторого разреза. \vskip 1 cm

\unitlength 1mm 
\linethickness{0.4pt}
\ifx\plotpoint\undefined\newsavebox{\plotpoint}\fi 
\begin{picture}(156.508,51.75)(0,0)
\put(44.258,26.5){\line(0,1){.7575}}
\put(44.239,27.257){\line(0,1){.7556}}
\put(44.183,28.013){\line(0,1){.7519}}
\multiput(44.089,28.765)(-.03276,.18657){4}{\line(0,1){.18657}}
\multiput(43.958,29.511)(-.033585,.147769){5}{\line(0,1){.147769}}
\multiput(43.79,30.25)(-.0292,.104229){7}{\line(0,1){.104229}}
\multiput(43.586,30.98)(-.030045,.089819){8}{\line(0,1){.089819}}
\multiput(43.345,31.698)(-.030638,.078415){9}{\line(0,1){.078415}}
\multiput(43.07,32.404)(-.031043,.069118){10}{\line(0,1){.069118}}
\multiput(42.759,33.095)(-.031306,.061356){11}{\line(0,1){.061356}}
\multiput(42.415,33.77)(-.031453,.054749){12}{\line(0,1){.054749}}
\multiput(42.037,34.427)(-.031507,.0490338){13}{\line(0,1){.0490338}}
\multiput(41.628,35.064)(-.0314807,.0440229){14}{\line(0,1){.0440229}}
\multiput(41.187,35.681)(-.0336273,.0424058){14}{\line(0,1){.0424058}}
\multiput(40.716,36.274)(-.0333116,.0379719){15}{\line(0,1){.0379719}}
\multiput(40.217,36.844)(-.0329583,.0340045){16}{\line(0,1){.0340045}}
\multiput(39.689,37.388)(-.0346058,.0323264){16}{\line(-1,0){.0346058}}
\multiput(39.136,37.905)(-.0385791,.0326065){15}{\line(-1,0){.0385791}}
\multiput(38.557,38.394)(-.043018,.0328405){14}{\line(-1,0){.043018}}
\multiput(37.955,38.854)(-.0480257,.0330233){13}{\line(-1,0){.0480257}}
\multiput(37.33,39.283)(-.05374,.033148){12}{\line(-1,0){.05374}}
\multiput(36.686,39.681)(-.060348,.033207){11}{\line(-1,0){.060348}}
\multiput(36.022,40.047)(-.068114,.033187){10}{\line(-1,0){.068114}}
\multiput(35.341,40.378)(-.07742,.033072){9}{\line(-1,0){.07742}}
\multiput(34.644,40.676)(-.088837,.032837){8}{\line(-1,0){.088837}}
\multiput(33.933,40.939)(-.103266,.032441){7}{\line(-1,0){.103266}}
\multiput(33.21,41.166)(-.122207,.031821){6}{\line(-1,0){.122207}}
\multiput(32.477,41.357)(-.148363,.030858){5}{\line(-1,0){.148363}}
\put(31.735,41.511){\line(-1,0){.7486}}
\put(30.987,41.628){\line(-1,0){.7535}}
\put(30.233,41.708){\line(-1,0){.7565}}
\put(29.477,41.751){\line(-1,0){.7577}}
\put(28.719,41.756){\line(-1,0){.757}}
\put(27.962,41.723){\line(-1,0){.7544}}
\put(27.208,41.653){\line(-1,0){.75}}
\multiput(26.458,41.545)(-.148746,-.028953){5}{\line(-1,0){.148746}}
\multiput(25.714,41.4)(-.122605,-.030251){6}{\line(-1,0){.122605}}
\multiput(24.978,41.219)(-.103673,-.031114){7}{\line(-1,0){.103673}}
\multiput(24.252,41.001)(-.089251,-.031695){8}{\line(-1,0){.089251}}
\multiput(23.538,40.747)(-.077837,-.032077){9}{\line(-1,0){.077837}}
\multiput(22.838,40.459)(-.068534,-.032311){10}{\line(-1,0){.068534}}
\multiput(22.153,40.135)(-.060769,-.03243){11}{\line(-1,0){.060769}}
\multiput(21.484,39.779)(-.05416,-.032457){12}{\line(-1,0){.05416}}
\multiput(20.834,39.389)(-.0484452,-.0324048){13}{\line(-1,0){.0484452}}
\multiput(20.204,38.968)(-.0434356,-.0322862){14}{\line(-1,0){.0434356}}
\multiput(19.596,38.516)(-.038994,-.0321092){15}{\line(-1,0){.038994}}
\multiput(19.011,38.034)(-.0350174,-.03188){16}{\line(-1,0){.0350174}}
\multiput(18.451,37.524)(-.0333916,-.0335791){16}{\line(0,-1){.0335791}}
\multiput(17.917,36.987)(-.0316835,-.0351953){16}{\line(0,-1){.0351953}}
\multiput(17.41,36.424)(-.0318904,-.0391731){15}{\line(0,-1){.0391731}}
\multiput(16.932,35.836)(-.0320426,-.0436156){14}{\line(0,-1){.0436156}}
\multiput(16.483,35.226)(-.0321331,-.0486258){13}{\line(0,-1){.0486258}}
\multiput(16.065,34.594)(-.032153,-.054341){12}{\line(0,-1){.054341}}
\multiput(15.679,33.941)(-.03209,-.060949){11}{\line(0,-1){.060949}}
\multiput(15.326,33.271)(-.031927,-.068714){10}{\line(0,-1){.068714}}
\multiput(15.007,32.584)(-.031641,-.078016){9}{\line(0,-1){.078016}}
\multiput(14.722,31.882)(-.031195,-.089427){8}{\line(0,-1){.089427}}
\multiput(14.473,31.166)(-.030534,-.103846){7}{\line(0,-1){.103846}}
\multiput(14.259,30.439)(-.029564,-.122772){6}{\line(0,-1){.122772}}
\multiput(14.082,29.703)(-.02812,-.148906){5}{\line(0,-1){.148906}}
\put(13.941,28.958){\line(0,-1){.7506}}
\put(13.838,28.208){\line(0,-1){.7548}}
\put(13.772,27.453){\line(0,-1){1.5148}}
\put(13.752,25.938){\line(0,-1){.7562}}
\put(13.799,25.182){\line(0,-1){.753}}
\multiput(13.883,24.429)(.03037,-.18697){4}{\line(0,-1){.18697}}
\multiput(14.005,23.681)(.031688,-.148188){5}{\line(0,-1){.148188}}
\multiput(14.163,22.94)(.032504,-.122027){6}{\line(0,-1){.122027}}
\multiput(14.358,22.208)(.033019,-.103082){7}{\line(0,-1){.103082}}
\multiput(14.589,21.486)(.033333,-.088652){8}{\line(0,-1){.088652}}
\multiput(14.856,20.777)(.033505,-.077233){9}{\line(0,-1){.077233}}
\multiput(15.157,20.082)(.033568,-.067927){10}{\line(0,-1){.067927}}
\multiput(15.493,19.403)(.033544,-.060161){11}{\line(0,-1){.060161}}
\multiput(15.862,18.741)(.033449,-.053553){12}{\line(0,-1){.053553}}
\multiput(16.263,18.098)(.0332916,-.0478401){13}{\line(0,-1){.0478401}}
\multiput(16.696,17.476)(.0330808,-.0428335){14}{\line(0,-1){.0428335}}
\multiput(17.159,16.877)(.032822,-.0383959){15}{\line(0,-1){.0383959}}
\multiput(17.652,16.301)(.0325196,-.0344243){16}{\line(0,-1){.0344243}}
\multiput(18.172,15.75)(.0341884,-.0327675){16}{\line(1,0){.0341884}}
\multiput(18.719,15.226)(.0381578,-.0330985){15}{\line(1,0){.0381578}}
\multiput(19.291,14.729)(.0425934,-.0333894){14}{\line(1,0){.0425934}}
\multiput(19.888,14.262)(.0475983,-.0336364){13}{\line(1,0){.0475983}}
\multiput(20.506,13.824)(.0492094,-.031232){13}{\line(1,0){.0492094}}
\multiput(21.146,13.418)(.054924,-.031147){12}{\line(1,0){.054924}}
\multiput(21.805,13.045)(.06153,-.030962){11}{\line(1,0){.06153}}
\multiput(22.482,12.704)(.06929,-.030656){10}{\line(1,0){.06929}}
\multiput(23.175,12.397)(.078585,-.030198){9}{\line(1,0){.078585}}
\multiput(23.882,12.126)(.089986,-.029542){8}{\line(1,0){.089986}}
\multiput(24.602,11.889)(.121789,-.033385){6}{\line(1,0){.121789}}
\multiput(25.333,11.689)(.147955,-.032757){5}{\line(1,0){.147955}}
\multiput(26.073,11.525)(.18675,-.03172){4}{\line(1,0){.18675}}
\put(26.82,11.398){\line(1,0){.7524}}
\put(27.572,11.309){\line(1,0){.7559}}
\put(28.328,11.257){\line(1,0){1.5149}}
\put(29.843,11.265){\line(1,0){.7553}}
\put(30.598,11.326){\line(1,0){.7513}}
\multiput(31.349,11.424)(.149105,.027043){5}{\line(1,0){.149105}}
\multiput(32.095,11.559)(.122983,.028676){6}{\line(1,0){.122983}}
\multiput(32.833,11.731)(.104064,.029782){7}{\line(1,0){.104064}}
\multiput(33.561,11.94)(.08965,.030548){8}{\line(1,0){.08965}}
\multiput(34.278,12.184)(.078242,.031076){9}{\line(1,0){.078242}}
\multiput(34.983,12.464)(.068943,.03143){10}{\line(1,0){.068943}}
\multiput(35.672,12.778)(.06118,.031649){11}{\line(1,0){.06118}}
\multiput(36.345,13.126)(.054572,.031759){12}{\line(1,0){.054572}}
\multiput(37,13.507)(.0488567,.031781){13}{\line(1,0){.0488567}}
\multiput(37.635,13.92)(.043846,.0317266){14}{\line(1,0){.043846}}
\multiput(38.249,14.364)(.0394025,.0316065){15}{\line(1,0){.0394025}}
\multiput(38.84,14.839)(.0377849,.0335236){15}{\line(1,0){.0377849}}
\multiput(39.407,15.341)(.0338195,.0331482){16}{\line(1,0){.0338195}}
\multiput(39.948,15.872)(.0321322,.0347862){16}{\line(0,1){.0347862}}
\multiput(40.462,16.428)(.0323901,.038761){15}{\line(0,1){.038761}}
\multiput(40.948,17.01)(.0325992,.0432012){14}{\line(0,1){.0432012}}
\multiput(41.404,17.615)(.032754,.0482098){13}{\line(0,1){.0482098}}
\multiput(41.83,18.241)(.032847,.053924){12}{\line(0,1){.053924}}
\multiput(42.224,18.888)(.032869,.060533){11}{\line(0,1){.060533}}
\multiput(42.586,19.554)(.032805,.068299){10}{\line(0,1){.068299}}
\multiput(42.914,20.237)(.032638,.077604){9}{\line(0,1){.077604}}
\multiput(43.207,20.936)(.032339,.089019){8}{\line(0,1){.089019}}
\multiput(43.466,21.648)(.031863,.103446){7}{\line(0,1){.103446}}
\multiput(43.689,22.372)(.031136,.122383){6}{\line(0,1){.122383}}
\multiput(43.876,23.106)(.030027,.148533){5}{\line(0,1){.148533}}
\put(44.026,23.849){\line(0,1){.7492}}
\put(44.139,24.598){\line(0,1){.7539}}
\put(44.215,25.352){\line(0,1){1.148}}
\put(156.508,26.5){\line(0,1){.7575}}
\put(156.489,27.257){\line(0,1){.7556}}
\put(156.433,28.013){\line(0,1){.7519}}
\multiput(156.339,28.765)(-.03276,.18657){4}{\line(0,1){.18657}}
\multiput(156.208,29.511)(-.033585,.147769){5}{\line(0,1){.147769}}
\multiput(156.04,30.25)(-.0292,.104229){7}{\line(0,1){.104229}}
\multiput(155.836,30.98)(-.030045,.089819){8}{\line(0,1){.089819}}
\multiput(155.595,31.698)(-.030638,.078415){9}{\line(0,1){.078415}}
\multiput(155.32,32.404)(-.031043,.069118){10}{\line(0,1){.069118}}
\multiput(155.009,33.095)(-.031306,.061356){11}{\line(0,1){.061356}}
\multiput(154.665,33.77)(-.031453,.054749){12}{\line(0,1){.054749}}
\multiput(154.287,34.427)(-.031507,.0490338){13}{\line(0,1){.0490338}}
\multiput(153.878,35.064)(-.0314807,.0440229){14}{\line(0,1){.0440229}}
\multiput(153.437,35.681)(-.0336273,.0424058){14}{\line(0,1){.0424058}}
\multiput(152.966,36.274)(-.0333116,.0379719){15}{\line(0,1){.0379719}}
\multiput(152.467,36.844)(-.0329583,.0340045){16}{\line(0,1){.0340045}}
\multiput(151.939,37.388)(-.0346058,.0323264){16}{\line(-1,0){.0346058}}
\multiput(151.386,37.905)(-.0385791,.0326065){15}{\line(-1,0){.0385791}}
\multiput(150.807,38.394)(-.043018,.0328405){14}{\line(-1,0){.043018}}
\multiput(150.205,38.854)(-.0480257,.0330233){13}{\line(-1,0){.0480257}}
\multiput(149.58,39.283)(-.05374,.033148){12}{\line(-1,0){.05374}}
\multiput(148.936,39.681)(-.060348,.033207){11}{\line(-1,0){.060348}}
\multiput(148.272,40.047)(-.068114,.033187){10}{\line(-1,0){.068114}}
\multiput(147.591,40.378)(-.07742,.033072){9}{\line(-1,0){.07742}}
\multiput(146.894,40.676)(-.088837,.032837){8}{\line(-1,0){.088837}}
\multiput(146.183,40.939)(-.103266,.032441){7}{\line(-1,0){.103266}}
\multiput(145.46,41.166)(-.122207,.031821){6}{\line(-1,0){.122207}}
\multiput(144.727,41.357)(-.148363,.030858){5}{\line(-1,0){.148363}}
\put(143.985,41.511){\line(-1,0){.7486}}
\put(143.237,41.628){\line(-1,0){.7535}}
\put(142.483,41.708){\line(-1,0){.7565}}
\put(141.727,41.751){\line(-1,0){.7577}}
\put(140.969,41.756){\line(-1,0){.757}}
\put(140.212,41.723){\line(-1,0){.7544}}
\put(139.458,41.653){\line(-1,0){.75}}
\multiput(138.708,41.545)(-.148746,-.028953){5}{\line(-1,0){.148746}}
\multiput(137.964,41.4)(-.122605,-.030251){6}{\line(-1,0){.122605}}
\multiput(137.228,41.219)(-.103673,-.031114){7}{\line(-1,0){.103673}}
\multiput(136.502,41.001)(-.089251,-.031695){8}{\line(-1,0){.089251}}
\multiput(135.788,40.747)(-.077837,-.032077){9}{\line(-1,0){.077837}}
\multiput(135.088,40.459)(-.068534,-.032311){10}{\line(-1,0){.068534}}
\multiput(134.403,40.135)(-.060769,-.03243){11}{\line(-1,0){.060769}}
\multiput(133.734,39.779)(-.05416,-.032457){12}{\line(-1,0){.05416}}
\multiput(133.084,39.389)(-.0484452,-.0324048){13}{\line(-1,0){.0484452}}
\multiput(132.454,38.968)(-.0434356,-.0322862){14}{\line(-1,0){.0434356}}
\multiput(131.846,38.516)(-.038994,-.0321092){15}{\line(-1,0){.038994}}
\multiput(131.261,38.034)(-.0350174,-.03188){16}{\line(-1,0){.0350174}}
\multiput(130.701,37.524)(-.0333916,-.0335791){16}{\line(0,-1){.0335791}}
\multiput(130.167,36.987)(-.0316835,-.0351953){16}{\line(0,-1){.0351953}}
\multiput(129.66,36.424)(-.0318904,-.0391731){15}{\line(0,-1){.0391731}}
\multiput(129.182,35.836)(-.0320426,-.0436156){14}{\line(0,-1){.0436156}}
\multiput(128.733,35.226)(-.0321331,-.0486258){13}{\line(0,-1){.0486258}}
\multiput(128.315,34.594)(-.032153,-.054341){12}{\line(0,-1){.054341}}
\multiput(127.929,33.941)(-.03209,-.060949){11}{\line(0,-1){.060949}}
\multiput(127.576,33.271)(-.031927,-.068714){10}{\line(0,-1){.068714}}
\multiput(127.257,32.584)(-.031641,-.078016){9}{\line(0,-1){.078016}}
\multiput(126.972,31.882)(-.031195,-.089427){8}{\line(0,-1){.089427}}
\multiput(126.723,31.166)(-.030534,-.103846){7}{\line(0,-1){.103846}}
\multiput(126.509,30.439)(-.029564,-.122772){6}{\line(0,-1){.122772}}
\multiput(126.332,29.703)(-.02812,-.148906){5}{\line(0,-1){.148906}}
\put(126.191,28.958){\line(0,-1){.7506}}
\put(126.088,28.208){\line(0,-1){.7548}}
\put(126.022,27.453){\line(0,-1){1.5148}}
\put(126.002,25.938){\line(0,-1){.7562}}
\put(126.049,25.182){\line(0,-1){.753}}
\multiput(126.133,24.429)(.03037,-.18697){4}{\line(0,-1){.18697}}
\multiput(126.255,23.681)(.031688,-.148188){5}{\line(0,-1){.148188}}
\multiput(126.413,22.94)(.032504,-.122027){6}{\line(0,-1){.122027}}
\multiput(126.608,22.208)(.033019,-.103082){7}{\line(0,-1){.103082}}
\multiput(126.839,21.486)(.033333,-.088652){8}{\line(0,-1){.088652}}
\multiput(127.106,20.777)(.033505,-.077233){9}{\line(0,-1){.077233}}
\multiput(127.407,20.082)(.033568,-.067927){10}{\line(0,-1){.067927}}
\multiput(127.743,19.403)(.033544,-.060161){11}{\line(0,-1){.060161}}
\multiput(128.112,18.741)(.033449,-.053553){12}{\line(0,-1){.053553}}
\multiput(128.513,18.098)(.0332916,-.0478401){13}{\line(0,-1){.0478401}}
\multiput(128.946,17.476)(.0330808,-.0428335){14}{\line(0,-1){.0428335}}
\multiput(129.409,16.877)(.032822,-.0383959){15}{\line(0,-1){.0383959}}
\multiput(129.902,16.301)(.0325196,-.0344243){16}{\line(0,-1){.0344243}}
\multiput(130.422,15.75)(.0341884,-.0327675){16}{\line(1,0){.0341884}}
\multiput(130.969,15.226)(.0381578,-.0330985){15}{\line(1,0){.0381578}}
\multiput(131.541,14.729)(.0425934,-.0333894){14}{\line(1,0){.0425934}}
\multiput(132.138,14.262)(.0475983,-.0336364){13}{\line(1,0){.0475983}}
\multiput(132.756,13.824)(.0492094,-.031232){13}{\line(1,0){.0492094}}
\multiput(133.396,13.418)(.054924,-.031147){12}{\line(1,0){.054924}}
\multiput(134.055,13.045)(.06153,-.030962){11}{\line(1,0){.06153}}
\multiput(134.732,12.704)(.06929,-.030656){10}{\line(1,0){.06929}}
\multiput(135.425,12.397)(.078585,-.030198){9}{\line(1,0){.078585}}
\multiput(136.132,12.126)(.089986,-.029542){8}{\line(1,0){.089986}}
\multiput(136.852,11.889)(.121789,-.033385){6}{\line(1,0){.121789}}
\multiput(137.583,11.689)(.147955,-.032757){5}{\line(1,0){.147955}}
\multiput(138.323,11.525)(.18675,-.03172){4}{\line(1,0){.18675}}
\put(139.07,11.398){\line(1,0){.7524}}
\put(139.822,11.309){\line(1,0){.7559}}
\put(140.578,11.257){\line(1,0){1.5149}}
\put(142.093,11.265){\line(1,0){.7553}}
\put(142.848,11.326){\line(1,0){.7513}}
\multiput(143.599,11.424)(.149105,.027043){5}{\line(1,0){.149105}}
\multiput(144.345,11.559)(.122983,.028676){6}{\line(1,0){.122983}}
\multiput(145.083,11.731)(.104064,.029782){7}{\line(1,0){.104064}}
\multiput(145.811,11.94)(.08965,.030548){8}{\line(1,0){.08965}}
\multiput(146.528,12.184)(.078242,.031076){9}{\line(1,0){.078242}}
\multiput(147.233,12.464)(.068943,.03143){10}{\line(1,0){.068943}}
\multiput(147.922,12.778)(.06118,.031649){11}{\line(1,0){.06118}}
\multiput(148.595,13.126)(.054572,.031759){12}{\line(1,0){.054572}}
\multiput(149.25,13.507)(.0488567,.031781){13}{\line(1,0){.0488567}}
\multiput(149.885,13.92)(.043846,.0317266){14}{\line(1,0){.043846}}
\multiput(150.499,14.364)(.0394025,.0316065){15}{\line(1,0){.0394025}}
\multiput(151.09,14.839)(.0377849,.0335236){15}{\line(1,0){.0377849}}
\multiput(151.657,15.341)(.0338195,.0331482){16}{\line(1,0){.0338195}}
\multiput(152.198,15.872)(.0321322,.0347862){16}{\line(0,1){.0347862}}
\multiput(152.712,16.428)(.0323901,.038761){15}{\line(0,1){.038761}}
\multiput(153.198,17.01)(.0325992,.0432012){14}{\line(0,1){.0432012}}
\multiput(153.654,17.615)(.032754,.0482098){13}{\line(0,1){.0482098}}
\multiput(154.08,18.241)(.032847,.053924){12}{\line(0,1){.053924}}
\multiput(154.474,18.888)(.032869,.060533){11}{\line(0,1){.060533}}
\multiput(154.836,19.554)(.032805,.068299){10}{\line(0,1){.068299}}
\multiput(155.164,20.237)(.032638,.077604){9}{\line(0,1){.077604}}
\multiput(155.457,20.936)(.032339,.089019){8}{\line(0,1){.089019}}
\multiput(155.716,21.648)(.031863,.103446){7}{\line(0,1){.103446}}
\multiput(155.939,22.372)(.031136,.122383){6}{\line(0,1){.122383}}
\multiput(156.126,23.106)(.030027,.148533){5}{\line(0,1){.148533}}
\put(156.276,23.849){\line(0,1){.7492}}
\put(156.389,24.598){\line(0,1){.7539}}
\put(156.465,25.352){\line(0,1){1.148}}
\qbezier(75.75,36.75)(81.625,46.875)(111,50.5)
\qbezier(75.5,38)(81.375,48.125)(110.75,51.75)
\multiput(75.75,36.75)(-.03125,.125){8}{\line(0,1){.125}}
\put(77.5,47){\makebox(0,0)[cc]{$L$}}
\put(75.75,37){\circle*{1.5}}
\put(72.75,35.75){\makebox(0,0)[cc]{$z_0$}}
\put(87.5,46){\circle*{1.5}}
\put(89,42.){\makebox(0,0)[cc]{$w(t)$}}
\put(29,26.5){\circle*{1.5}} \put(86.5,28){\circle*{1.5}}
\qbezier(147,40.75)(144.125,37)(136.75,33.25)
\qbezier(148.25,39.75)(145.375,36)(138,32.25)
\multiput(137.25,33.25)(.0333333,-.0833333){15}{\line(0,-1){.0833333}}
\put(148,45.25){\makebox(0,0)[cc]{$\lambda(t)$}}
\put(123.25,33){\vector(3,-1){.07}}\qbezier(100.25,34.25)(111.75,36.875)(123.25,33)
\put(141.25,26.25){\circle*{1.5}} \put(137.75,33.){\circle*{1.5}}
\put(147.5,40.25){\circle*{1.5}}
\put(144.25,26.75){\makebox(0,0)[cc]{$0$}}
\put(29.25,36.75){\makebox(0,0)[cc]{$E$}}
\put(32,28.25){\makebox(0,0)[cc]{$0$}}
\put(83.5,27.75){\makebox(0,0)[cc]{$0$}}
\put(63.75,32.25){\vector(3,-1){.07}}\qbezier(46.5,33)(57.375,34.875)(63.75,32.25)
\put(54.5,37.75){\makebox(0,0)[cc]{$f$}}
\put(112.75,39.25){\makebox(0,0)[cc]{$f_{\,t}^{-1}$}}
\put(80.5,4.25){\makebox(0,0)[cc]{Рис.~3}}
\end{picture}
\vskip 0.5cm

Пусть $\lambda(t)$ --- точка которая переходит в конец разреза при
отображении $f_t$. Имеют место

\begin{theorem}\label{appr} Справедливо равенство
$$
\lim_{t\to+\infty}(e^{-t}\varphi_t(z))=f(z),
$$
причем сходимость равномерна на компактах в единичном круге.
\end{theorem}

\begin{theorem}\label{loewner}
Функции $\varphi_t$ удовлетворяют дифференциальному уравнению
\begin{equation}\label{loew1}
\dot{\varphi}_t=-\varphi_t\,
\frac{\lambda(t)+\varphi_t}{\lambda(t)-\varphi_t}\,,
\end{equation}
а функции $f_t$ --- уравнению
\begin{equation}\label{loew2}
\dot{f}_t=zf'_t\, \frac{\lambda(t)+f_t}{\lambda(t)-f_t}\,.
\end{equation}
\end{theorem}

Здесь и далее точка сверху означает дифференцирование по параметру
$t$, а штрих --- по комплексной переменной~$z$.

Уравнения (\ref{loew1}) и  (\ref{loew2}) играют выдающуюся роль в
геометрической теории функций и носят название уравнений Левнера.
Уравнение (\ref{loew2}) описывает динамику конформных отображений
круга на плоскость с укорачивающимся разрезом, а
(\ref{loew1})
--- динамику обратных отображений.

Заметим, что в правых частях (\ref{loew1}) и (\ref{loew2})
содержится выражение $$\frac{\lambda(t)+w}{\lambda(t)-w}\,,$$
подсчитываемое при $w=\varphi_t$ и $w=f_t$. Это --- дробно-линейное
отображение круга $E$, удовлетворяющее неравенству
$\Re\frac{\lambda(t)+w}{\lambda(t)-w}>0$, $w\in E$. Более общие
уравнения получаются, если вместо дробно-линейного отображения
использовать достаточно произвольное отображение с положительной
вещественной частью. Такие уравнения называются уравнениями
Левнера-Куфарева.

Таким образом, теоремы~\ref{appr} и \ref{loewner} показывают, что
любая функция $f$ из класса $\widetilde{S}$ может быть
аппроксимирована однопараметрическим семейством $\varphi_t$,
которое является решением уравнения Левнера~(\ref{loew1}). Суть
параметрического метода заключается в исследовании функционалов на
семействах решений уравнения~(\ref{loew1}); при этом  функция
$\lambda$ является своего рода управлением, определяющим
траекторию дуги $L$, задающей границу образа $f(E)$.

Заметим, что в последние годы уравнения левнеровского типа являются
предметом пристального изучения математиков. Здесь можно упомянуть
исследования однопараметрических семейств конформных отображений
полуплоскости с гидродинамической нормировкой (так называемое
хордальное уравнение Левнера), уравнения Левнера-Шрамма и приложения
однопараметрических семейств в теории ветвящихся процессов и теории
вероятностей (см, напр., \cite{sle,sle2,sle1,gor1,gor2}).

Далее мы опишем однопараметрические семейства неоднолистных
функций и получим  для них дифференциальные уравнения, которые по
форме напоминают уравнения Левнера и Левнера-Куфарева.

\subsection{Однопараметрические семейства римановых поверхностей}

Пусть задана $n$-листная компактная риманова поверхность $S$ над
сферой Римана $\overline{\mathbb{C}}$. По известной теореме
униформизации в случае односвязной поверхности $S$ существует
рациональная функция $f$, которая отображает сферу Римана
$\overline{\mathbb{C}}$ на $S$. В случае, если род поверхности $S$
равен $1$, она униформизируется эллиптической функцией, т.~е.
существует функция  $f$ с периодами $\omega_1$ и $\omega_2$ такая,
что фактор-отображение
$\widetilde{f}:\mathbb{C}/$\mbox{\boldmath$\omega$}$\to S$
является конформным изоморфизмом. Здесь $\mbox{\boldmath$\omega$}$
--- решетка, порожденная периодами $\omega_1$ и $\omega_2$.
\vskip 1 cm

\hskip 1 cm
\unitlength 1mm 
\linethickness{0.4pt}
\ifx\plotpoint\undefined\newsavebox{\plotpoint}\fi 
\begin{picture}(135.642,73.977)(0,0)
\put(5,36.122){\vector(1,0){48}}
\put(25.5,21.122){\vector(0,1){37}}
\put(36.75,45.872){\circle*{1.5}} \put(39.5,28.872){\circle*{1.5}}
\put(15.5,40.372){\circle*{1.5}}
\put(41,46.372){\makebox(0,0)[cc]{$a_j$}}
\put(42.75,29.372){\makebox(0,0)[cc]{$a_k$}}
\put(18.75,41.122){\makebox(0,0)[cc]{$a_l$}}
\multiput(119.946,43.727)(.0336538462,.0519230769){260}{\line(0,1){.0519230769}}
\multiput(78.217,43.622)(.0336538462,.0519230769){260}{\line(0,1){.0519230769}}
\multiput(126.5,62.372)(.03372093,.053488372){215}{\line(0,1){.053488372}}
\multiput(127.446,62.477)(.03372093,.053488372){215}{\line(0,1){.053488372}}
\multiput(128.392,62.372)(.03372093,.053488372){215}{\line(0,1){.053488372}}
\multiput(85.717,62.372)(.03372093,.053488372){215}{\line(0,1){.053488372}}
\qbezier(126.975,59.706)(126.187,61.125)(126.45,62.334)
\qbezier(127.921,59.811)(127.133,61.23)(127.396,62.439)
\qbezier(86.192,59.706)(85.404,61.125)(85.666,62.334)
\qbezier(127.606,56.973)(127.869,58.235)(127.08,59.496)
\qbezier(128.552,57.078)(128.815,58.34)(128.026,59.601)
\qbezier(86.823,56.973)(87.085,58.235)(86.297,59.496)
\put(92.814,73.686){\line(1,0){40.783}}
\multiput(118.987,43.519)(.0336687041,.052556026){256}{\line(0,1){.052556026}}
\put(78.203,43.309){\line(1,0){39.627}}
\put(78.115,42.69){\line(1,0){39.627}}
\put(78.027,41.983){\line(1,0){39.627}}
\multiput(117.831,43.309)(.033620399,.0510102605){272}{\line(0,1){.0510102605}}
\qbezier(127.396,59.706)(127.501,61.073)(128.447,62.439)
\qbezier(127.501,59.601)(127.816,58.182)(126.87,56.973)
\put(127.606,58.655){\line(-1,0){10.091}}
\put(117.62,58.55){\circle*{1.5}}
\multiput(78.224,42.563)(.0372161,.0325641){19}{\line(1,0){.0372161}}
\multiput(78.135,41.944)(.0336718,.0336718){21}{\line(1,0){.0336718}}
\multiput(118.995,43.48)(-.0318538,-.0530896){14}{\line(0,-1){.0530896}}
\put(117.657,42.588){\line(1,0){.8176}}
\put(117.583,41.994){\line(1,0){.8919}}
\multiput(120.031,43.8)(-.0331456,-.0524806){32}{\line(0,-1){.0524806}}
\put(118.087,42.033){\line(1,0){.884}}
\put(86.515,59.163){\line(1,0){8.027}}
\put(94.542,59.014){\circle*{1.5}}
\multiput(104.799,67.785)(.03356633,.046753102){124}{\line(0,1){.046753102}}
\put(104.948,67.339){\circle*{1.5}}
\multiput(118.772,13.081)(.0336885868,.0631661003){353}{\line(0,1){.0631661003}}
\multiput(77.15,13.081)(.0336885868,.0631661003){353}{\line(0,1){.0631661003}}
\put(118.921,12.933){\line(-1,0){41.771}}
\put(89.042,35.082){\line(1,0){41.622}}
\put(102.718,30.771){\circle*{1.5}}
\put(116.245,23.338){\circle*{1.5}}
\put(92.164,23.784){\circle*{1.5}}
\put(99.299,40.433){\vector(0,-1){11.297}}
\put(95.353,22.298){\makebox(0,0)[cc]{$A_l$}}
\put(106.948,30.919){\makebox(0,0)[cc]{$A_j$}}
\put(120.029,25.){\makebox(0,0)[cc]{$A_k$}}
\put(125.461,31.703){\makebox(0,0)[cc]{$\overline{\mathbb{C}}$}}
\put(63.325,3.568){\makebox(0,0)[cc]{Рис.~4}}
\put(65,43){\makebox(0,0)[cc]{$f$}}
\put(54,45.){\vector(2,-1){20}}
\end{picture}

\vskip 0.4 cm

Так как поверхность $S$ задана, известны точки $A_j$, которые
являются проекциями на плоскость ее точек ветвления, т.~е.
критические значения $f$ (рис.~4). Важной задачей является
проблема нахождения критических точек $a_j$ отображения $f$ по
заданным значениям $A_j$.

В случае рациональных функций система $f(a_j)=A_j$ является
системой алгебраических уравнений, ее решение представляет
известные трудности. В случае эллиптических функций ситуация
усложняется тем, что кроме критических точек требуется определить
также решетку периодов~$\mbox{\boldmath$\omega$}$.

Отметим, что по заданным $A_j$ точки $a_j$ определяются, вообще
говоря, не единственным образом. А.~Гурвиц \cite{hur}, \cite{hur2}
впервые поставил проблему об определении числа неэквивалентных
накрытий сферы с заданным типом ветвления и предложил некоторые
подходы к ее решению. Проблема Гурвица и ее обобщения
исследовались, в частности, в работах
\cite{hur,hur2,lloyd,weyl,med,med2,med3,med4,nas2}; описание
недавних исследований, основанных на изучении геометрии
расслоений, можно найти в обзоре \cite{lando} и монографии
\cite{lando2}.

В \cite{nas} мы предложили приближенный метод нахождения параметров
$a_j$  в случае полиномов. Суть его состоит в соединении
(посредством гладкой гомотопии) поверхности $S$ с поверхностью
$S_0$, для которой соответствующие параметры $a^{0}_j$ известны, в
пространстве односвязных $n$-листных разветвленных накрытий сферы с
топологией, индуцируемой сходимостью к ядру по Каратеодори (по
поводу сходимости к ядру римановых поверхностей см., напр.
\cite{nas}).

Прообраз этой гомотопии в пространстве <<акцессорных>> параметров
$a_j$ дает гладкий путь, который является интегральной кривой для
системы обыкновенных дифференциальных уравнений. Решая эту систему
с начальным условием $(a^{0}_j)$, получаем уравнение этой кривой,
концевая точка которой даст значение параметров $a_j$, $1\le j\le
n-2$, для поверхности~$S$.

В следующем пункте мы, следуя работе \cite{nas3}, рассмотрим
однопараметрические семейства рациональных функций.

\subsection{Поверхности рода нуль. Униформизация рациональными функциями}

Пусть $R(z,t)$ -- семейство рациональных функций, имеющих в
плоскости $\mathbb{C}$ критические точки $a_l=a_l(t)$  порядков
$m_l-1$, $1\le l \le M$, и полюсы $b_j=b_j(t)$ порядков $n_j$,
$1\le j\le N$. Обозначим
$$m=\sum_{k=1}^M(m_k-1),\quad n=\sum_{j=1}^N(n_j+1).$$

В случае $m>n$ будем рассматривать семейство
\begin{equation}\label{ra}
R(z,t)=\int_{a_1}^z\frac{\prod_{l=1}^M(\zeta-a_l)^{m_l-1}d\zeta}{\prod_{j=1}^N(\zeta-b_j)^{n_j+1}},
\end{equation}
в случае $m<n-1$ --- семейство
\begin{equation}\label{ra1}
R(z,t)=\int_{\infty}^z\frac{\prod_{l=1}^M(\zeta-a_l)^{m_l-1}d\zeta}{\prod_{j=1}^N(\zeta-b_j)^{n_j+1}},
\end{equation}
наконец, при $m=n$
\begin{equation}\label{ra2}
R(z,t)=z+\int_{\infty}^z\left[\frac{\prod_{l=1}^M(\zeta-a_l)^{m_l-1}}{\prod_{j=1}^N(\zeta-b_j)^{n_j+1}}-1\right]d\zeta.
\end{equation}

Отметим, что $m\neq n-1$, поскольку функции $R(z,t)$ однозначны и
вычет производной $R'(z,t)$ на бесконечности равен нулю. При этом,
$R'(z,t)\sim cz^{m-n}$, где $c=1$; последнего можно добиться
применением линейного преобразования в $z$-плоскости.

Запишем разложение функции $R(z,t)$ в окрестности $a_l$:
\begin{equation}\label{7}
R(z,t)=A_l+B_l(z-a_l)^{m_l}+\ldots,\quad
\dot{R}(z,t)=\dot{A}_l-m_l\dot{a}_lB_l(z-a_l)^{m_l-1}+\ldots
\end{equation}
Имеем

\begin{equation}\label{8}
(R'(z,t))^{-1}=\frac{H_l(z,t)}{(z-a_l)^{m_l-1}},\quad \mbox{\rm
где} \quad H_l(z,t)=\frac{\prod_{j=1}^N(z-b_j)^{n_j+1}}
{\prod_{k=1, k\neq l}^M(z-a_k)^{m_k-1}}.
\end{equation}

Обозначим через $P_{l,j}$ многочлен Тейлора функции $H_l(z,t)$
степени $j$, записанный в точке $a_l$: $$
P_{\,l,j}(z,t)=\sum_{s=0}^j\frac{1}{s!}\,H_l^{(s)}(a_l,t)(z-a_l)^s.
$$

Тогда при $z\to a_l$
\begin{equation}\label{10}
h(z,t):=\frac{\dot{R}(z,t)}{R'(z,t)}={\dot{R}(z,t)}\frac{H_l(z,t)}{(z-a_l)^{m_l-1}}=\frac{P_{\,l,m_l-2}(z,t)}{(z-a_l)^{m_l-1}}\,\dot{A}_l
+O(1).
\end{equation}

В окрестности полюса $b_j$
$$R(z,t)=C_j(z-b_j)^{-n_j}+\ldots,\ \dot{R}(z,t)=n_j\dot{b}_jC_j(z-b_j)^{-n_j-1}+\ldots,\ R'(z,t)=-n_jC_j(z-b_j)^{-n_j-1}+\ldots,
$$
поэтому
\begin{equation}\label{11}
h(z,t)=-\dot{b}_j+o(1).
\end{equation}
Наконец, в окрестности бесконечности
$$
\frac{\prod_{k=1}^M(\zeta-a_k)^{m_k-1}}{\prod_{j=1}^N(\zeta-b_j)^{n_j+1}}=z^{m-n}\frac{\prod_{k=1}^M(1-a_k/\zeta)^{m_k-1}}{\prod_{j=1}^N(1-b_j/\zeta)^{n_j+1}}=
$$$$=
\zeta^{m-n}\Bigl(1-\Bigl({\sum_{k=1}^M(m_k-1)a_k-\sum_{j=1}^N{(n_j+1)b_j}}\Bigr){\zeta^{-1}}+O\left({\zeta^{-2}}\right)\Bigr).
$$

Если $m\neq n$, то при сдвиге в плоскости $\zeta$ на $\zeta_0$
величина $\sum_{k=1}^M(m_k-1)a_k-\sum_{j=1}^N{(n_j+1)b_j}$
изменяется на $(m-n)\zeta_0$, поэтому за счет свободы в выборе
отображающей функции мы можем считать, что
\begin{equation}\label{12}
\sum_{k=1}^M(m_k-1)a_k-\sum_{j=1}^N{(n_j+1)b_j}=0.
\end{equation}
Если же $m=n$, то
$\sum_{k=1}^M(m_k-1)a_k-\sum_{j=1}^N{(n_j+1)b_j}=0$ как вычет
производной отображающей функции. Таким образом, мы можем считать,
что условие (\ref{12}) выполняется для всех функций семейства.
Тогда, учитывая (\ref{10}), (\ref{11}) и (\ref{13}), последовательно
находим
$$R'(z,t)=z^{m-n}(1+O(z^{-2})), \ R(z,t)=z^{m-n+1}\left(\frac{1}{m-n+1}+O(z^{-2})\right),\ \dot{R}(z,t)=z^{m-n+1}O(z^{-2})),$$
поэтому \begin{equation}\label{13} h(z,t)=o(1),\quad z\to\infty.
\end{equation}

Рациональная функция $h(z,t)$  имеет особенности только в точках
$a_l$, $b_j$ и на бесконечности, причем, в силу (\ref{11}) и
(\ref{13}), в точках $b_j$ и на бесконечности они устранимые. Из
(\ref{10}) следует, что
$$
h(z,t)=\sum_{l=1}^M\frac{P_{\,l,m_l-2}(z,t)}{(z-a_l)^{m_l-1}}\,\dot{A}_l.
$$

\begin{theorem}\label{rat}
Семейство функций ${R}(z,t)$ удовлетворяет дифференциальному
уравнению
\begin{equation}\label{14}
\frac{\dot{R}(z,t)}{R'(z,t)}=\sum_{l=1}^M\frac{P_{\,l,m_l-2}(z,t)}{(z-a_l)^{m_l-1}}\,\dot{A}_l.
\end{equation}
\end{theorem}

Отметим, что при $m>n$ в силу (\ref{ra}) $A_1\equiv 0$, поэтому
фактически сумма в (\ref{14}) берется по $2\le l \le M$.

Запишем правую части полученного уравнения в более компактном виде.
Для этого используем следующее легко проверяемое утверждение.

\begin{lemma}\label{l} Пусть $f$ голоморфна в точке \ $a$. Тогда
многочлен Тейлора степени $n$ для функции  $f$ в точке \ $a$ равен
\begin{equation*}
P_n(x)=\frac{(x-a)^{n+1}}{n!}\frac{\partial^{\,n}}{d\xi^n}\left.\left(\frac{f(\xi)}{x-\xi}\right)\right|_{\xi=a}\,.
\end{equation*}
\end{lemma}

С использованием этой леммы запишем

\begin{equation}\label{9}
P_{\,l,j}(z,t)=-\frac{(z-a_l)^{j+1}}{j!}\frac{\partial^{j}
G_l{}(a_l,z,t)}{\partial \xi^j}, \quad \mbox{\rm где}\quad
G_l(\xi,z,t)=\frac{H_l(\xi,t)}{\xi-z}\,.
\end{equation}

Теперь очевидно, что выражения
$$\frac{P_{\,l,m_l-2}(z,t)}{(z-a_l)^{m_l-1}},$$ которые являются
рациональными функциями от переменной $z$ и параметров $a_l$, $b_j$,
с учетом  (\ref{9}) могут быть записаны в виде:
$$\frac{P_{\,l,m_l-2}(z,t)}{(z-a_l)^{m_l-1}}=-\frac{1}{(m_l-2)!}\,\frac{\partial^{m_l-2}
G_l{}(a_l,z,t)}{\partial \xi^{m_l-2}}\,.$$

Таким образом, (\ref{14}) эквивалентно равенству
\begin{equation}\label{15}
\frac{\dot{R}(z,t)}{R'(z,t)}=-\sum_{l=1}^M\frac{1}{(m_l-2)!}\,\frac{\partial^{m_l-2}
G_l{}(a_l,z,t)}{\partial \xi^{m_l-2}}\,\dot{A}_l.
\end{equation}

Теперь найдем уравнения для определения $a_l$. В силу (\ref{14})
\begin{equation}\label{16}
{\dot{R}(z,t)}={R'(z,t)}\sum_{l=1}^M\frac{P_{\,l,m_l-2}(z,t)}{(z-a_l)^{m_l-1}}\,\dot{A}_l.
\end{equation}
с другой стороны, в силу (\ref{7})
$$\dot{R}(z,t)=\dot{A}_l-m_l\dot{a}_lB_l(z-a_l)^{m_l-1}+\ldots,$$ где
$$
m_lB_l=1/H_l(A_l,t),
$$
и для нахождения $\dot{a}_l$ достаточно определить коэффициент
$(-m_l\dot{a}_lB_l)$, соответствующий степени $(z-a_l)^{m_l-1}$ в
разложении в ряд Тейлора функции $\dot{R}(z,t)$. В силу (\ref{16})
$${\dot{R}(z,t)}=\frac{(z-a_l)^{m_l-1}}{H_l(z,t)}\left[\frac{P_{\,l,m_l-2}(z,t)}{(z-a_l)^{m_l-1}}\,\dot{A}_l+\sum_{k\neq l}\frac{P_{\,k,m_k-2}(z,t)}{(z-a_k)^{m_k-1}}\,\dot{A}_k\right]=$$
\begin{equation}\label{17}
=\frac{P_{\,l,m_l-2}(z,t)}{H_l(z,t)}\,\dot{A}_l+\frac{(z-a_l)^{m_l-1}}{H_l(z,t)}\sum_{k\neq
l}\frac{P_{\,k,m_k-2}(z,t)}{(z-a_k)^{m_k-1}}\,\dot{A}_k.
\end{equation}

\begin{lemma}\label{tayl}
{Если функция $f$ аналитична в окрестности точки $a$, $f(a)\neq
0$, и  $P_{n-1}$
--- многочлен Тейлора функции $f$ в точке $a$, то}
$$
\frac{P_{n-1}(x)}{f(x)}
=1-\frac{f^{(n)}(a)}{n!f(a)}(x-a)^n+o((x-a)^n), \quad x\to a.
$$
\end{lemma}

Доказательство. Имеем
$$
\frac{f(x)}{P_{n-1}(x)}=\frac{P_{n-1}(x)+\frac{f^{(n)}(a)}{n!}(x-a)^n+o((x-a)^n))}{P_{n-1}(x)}=
$$
$$
={1+\frac{f^{(n)}(a)}{n!P_{n-1}(x)}(x-a)^n+o((x-a)^n))}={1+\frac{f^{(n)}(a)}{n!f(a)}(x-a)^n+o((x-a)^n))},
\quad x\to a,
$$ откуда следует утверждение леммы.\hfill $\square$ \medskip

Применяя лемму к первому слагаемому в последней строке (\ref{17}),
получаем
$$
-m_l\dot{a}_lB_l=\frac{1}{H_l(a_l,t)}\left[-\frac{H_l^{(m_l-1)}(a_l,t)}{(m_l-1)!}+\sum_{k\neq
l}\frac{P_{\,k,m_k-2}(a_l,t)}{(a_l-a_k)^{m_k-1}}\,\dot{A}_k\right],
$$
откуда
$$
\dot{a}_l=\frac{H_l^{(m_l-1)}(a_l,t)}{(m_l-1)!}-\sum_{k\neq
l}\frac{P_{\,k,m_k-2}(a_l,t)}{(a_l-a_k)^{m_k-1}}\,\dot{A}_k=$$$$=\frac{H_l^{(m_l-1)}(a_l,t)}{(m_l-1)!}+\sum_{k\neq
l}\frac{1}{(m_k-2)!}\,\frac{\partial^{m_k-2}
G_l{}(a_k,a_l,t)}{\partial \xi^{m_k-2}}\,\dot{A}_k=
\frac{H_l^{(m_l-1)}(a_l)}{(m_l-1)!}\,\dot{A}_l+\sum_{k=1,k\neq l}
\frac{G_{kl}^{(m_k-2)}(a_k)}{(m_k-2)!}\,\dot{A}_k,
$$
где $$ G_{kl}(x)=\frac{H_k(x)}{x-a_l}.
$$

Из (\ref{11}) находим
$$
\dot{b}_j=-h(b_l,t)=\sum_{k=1}^M\frac{I_{kj}^{(m_k-2)}(a_k)}{(m_k-2)!}\,\dot{A}_k.
$$

Итак, установлена

\begin{theorem}\label{ratparam}Критические точки $a_l$ и полюсы $b_j$ семейства
рациональных функций удовлетворяют системе дифференциальных
уравнений
$$
\dot{a}_l=\frac{H_l^{(m_l-1)}(a_l)}{(m_l-1)!}\,\dot{A}_l+\sum_{k=1,k\neq
l}^M \frac{G_{kl}^{(m_k-2)}(a_k)}{(m_k-2)!}\,\dot{A}_k, \quad
\dot{b}_j=\sum_{k=1}^M\frac{I_{kj}^{(m_k-2)}(a_k)}{(m_k-2)!}\,\dot{A}_k.$$
где
$$
H_l(x)=\frac{\prod_{j=1}^N(x-b_j)^{n_j+1}}{\prod_{k=1, k\neq
l}^M(x-a_k)^{m_k-1}}\,,\quad G_{kl}(x)=\frac{H_k(x)}{x-a_l}\,,
\quad I_{kj}(x)=\frac{H_k(x)}{x-b_j}\,.$$

При этом, предполагается, что искомые параметры удовлетворяют
соотношению
$$\sum_{k=1}^M(m_k-1)a_k-\sum_{j=1}^N{(n_j+1)b_j}=0.$$
\end{theorem}

\subsection{Пример. Униформизация семейства односвязных компактных поверхностей}

Рассмотрим функцию
\begin{equation}\label{examplerat}
f(z)=z+\frac{2}{z}-\frac{1}{3z^3}\,.
\end{equation}
Ее производная
$$
f'(z)=\frac{(z^2-1)^2}{z^4}.
$$
Критические точки  $a_1=1$, $a_2=-1$ переходят в точки $A_1=8/3$,
$A_2=-8/{3}$.

Найдем рациональную функцию, униформизирующую поверхность, которая
получается из данной движением критических значений точек по
прямолинейным отрезкам из точек $A_1$, $A_2$ в точки $2$, $-1+i$.
Тогда
$\dot{A}_1=-2/3$,  $\dot{A}_2=5/{3}+i$.

Система дифференциальных уравнений для определения параметров имеет
вид
$$\dot{a}_1=\frac{(a_1-b)^4}{(a_1-a_2)^2}\left[\frac{6}{(a_1-b)^2}+\frac{3}{(a_1-a_2)^2}-\frac{8}{(a_1-b)(a_1-a_2)}\right]\,\dot{A}_1+$$
$$+\frac{(a_2-b)^4}{(a_2-a_1)^3}\left[\frac{4}{a_2-b}-\frac{3}{a_2-a_1}\right]\,\dot{A}_2,
$$
$$\dot{a}_2=\frac{(a_2-b)^4}{(a_2-a_1)^2}\left[\frac{6}{(a_2-b)^2}+\frac{3}{(a_2-a_1)^2}-\frac{8}{(a_2-b)(a_2-a_1)}\right]\,\dot{A}_2+$$
$$
+\frac{(a_1-b)^4}{(a_1-a_2)^3}\left[\frac{4}{a_1-b}-\frac{3}{a_1-a_2}\right]\,\dot{A}_1,
$$
$$\dot{b}=\frac{(a_1-b)^3}{(a_1-a_2)^2}\left[\frac{3}{a_1-b}-\frac{2}{a_1-a_2}\right]\dot{A}_1+
\frac{(a_2-b)^3}{(a_2-a_1)^2}\left[\frac{3}{a_2-b}-\frac{2}{a_2-a_1}\right]\,\dot{A}_2.$$

Решая задачу Коши для этой системы с начальными данными
$a_1(0)=1$, $a_2(0)=-1$,  $b(0)=0$, соответствующими функции
(\ref{examplerat}), находим
$a_1(1)=0.547419\ldots-i\,0.166211\ldots,$ $
a_2(1)=-0.778733\ldots+$ $+ i\, 0.611500\ldots,$
$b(1)=-0.5790162\ldots-i\,0.222208\ldots.$

Подсчеты показывают, что для найденных параметров координаты
критических значений отличаются от заданных на величины порядка
$10^{-7}$.

\subsection{Семейства комплексных торов и эллиптических функций}

Рассмотрим теперь однопараметрическое семейство эллиптических
функций $f(z,t)$ порядка $n\ge 2$ с периодами $\omega_1(t)$,
$\omega_2(t)$, гладко зависящее от вещественного параметра $t$.
Ограничимся случаем, когда у функции $f(z,t)$ при фиксированном
значении $t$ имеется \textit{единственный полюс в начале координат},
а все \textit{критические  точки --- простые.} Будем следовать
изложению, данному в \cite{nas4}.

Очевидно, что производная $f'(z,t)$ является эллиптической
функцией с теми же периодами $\omega_1(t)$, $\omega_2(t)$. Запишем
условия периодичности:
$$
f(z+\omega_k(t),t)=f(z,t),\quad k=1,2,
$$
и продифференцируем их по $t$. Получаем
$$
f'(z+\omega_k(t),t)\dot{\omega}_k(t)+\dot{f}(z+\omega_k(t),t)=\dot{f}(z,t),
$$
откуда с учетом периодичности $f'(z,t)$ приходим к равенствам
$$
\frac{\dot{f}(z+\omega_k(t),t)}{f'(z+\omega_k(t),t)}+\dot{\omega}_k(t)=\frac{\dot{f}(z,t)}{f'(z,t)}.
$$
Следовательно, функция $h(z,t):={\dot{f}(z,t)}/{f'(z,t)}$
удовлетворяет условиям
\begin{equation}\label{periodh}
h(z+\omega_k(t),t)-h(z,t)=-\dot{\omega}_k(t),\quad k=1,2.
\end{equation}
Теперь с учетом теоремы~\ref{sigma} представим $f'(z,t)$ в виде
$$
f'(z,t)=c(t)\,\frac{\prod\limits_{k=0}^n\sigma(z-a_k(t))}{\sigma^{n+1}(z)},
$$
где $a_k(t)$ --- полный набор критических точек функции $f'(z,t)$,
причем $\sum_{k=0}^n a_k(t)=0.$ Пусть $A_k(t)=f(a_k(t),t)$ ---
образы этих критических точек. Без ограничения общности можно
считать, что $A_0(t)\equiv 0$. Как и в случае полиномов, нашей
задачей будет нахождение системы дифференциальных уравнений,
которым удовлетворяют $a_k(t)$, если заданы зависимости $A_k(t)$,
$1\le k\le n$.

Запишем в окрестности точки $a_k(t)$ тейлоровское разложение
\begin{equation}\label{taylor}
 f(z,t)=A_k(t)+\frac{D_k(t)}{2}\,(z-a_k(t))^2+\ldots,\quad \textrm{где} \quad D_k(t)=f''(a_k(t),t).
\end{equation}
$$
$$
Имеем
$$
f''(z,t)=c(t)\,\frac{\prod\limits_{j=0}^n\sigma(z-a_j(t))}{\sigma^{n+1}(z)}\,\Bigl[\,\sum_{j=0}^n\zeta(z-a_j(t))-(n+1)\zeta(z)\Bigr],
$$
поэтому при $z\to a_k(t)$ находим
$$f''(a_k(t),t)=c(t)\,\frac{\prod\limits_{j\neq k}\sigma(a_k(t)-a_j(t))}{\sigma^{n+1}(a_k(t))}.$$
Итак,
\begin{equation}\label{dk}
D_k(t)={c(t)}\,\frac{\prod\limits_{j\neq
k}\sigma(a_k(t)-a_j(t))}{\sigma^{n+1}(a_k(t))}.
\end{equation}
 Из
(\ref{taylor}) следует, что
\begin{equation}\label{fprime}
 f'(z,t)=D_k(t)(z-a_k(t))+\ldots,
\end{equation}
\begin{equation}
 \dot{f}(z,t)=\dot{A}_k(t)-\dot{a}_k(t)D_k(t)(z-a_k(t))+\ldots,\label{(3.2)}
\end{equation}
откуда
\begin{equation}
 h(z,t)=\frac{\dot{f}(z,t)}{f'(z,t)}=\frac{\gamma_k(t)}{z-a_k(t)}+O(1),
\quad z\to a_k(t), \label{(3.3)}
\end{equation}
где
\begin{equation}\label{gammak}
\gamma_k(t):=\frac{\dot{A}_k(t)}{D_k(t)}\,.\end{equation}
 В окрестности точки $z=0$ функция $\dot{f}(z,t)$ имеет полюс
порядка не выше $n$, а $f'(z,t)$ --- $(n+1)$-го порядка, поэтому
$h(z,t)$ имеет в этой точке нуль.

Функция $$ g(z,t):=h(z,t)-\sum_{j=1}^n\gamma_j(t)\zeta(z-a_j(t))$$
имеет только устранимые особенности в точках $\omega$, а в
остальных точках плоскости она голоморфна. Следовательно, она
продолжается на $\mathbb{C}$ до целой функции.

Из (\ref{periodzeta}) и (\ref{periodh}) получаем
\begin{equation}\label{periodg}
g(z+\omega_k(t),t)-g(z,t)=-\dot{\omega}_k(t)-\eta_k(t)\sum_{j=1}^n\gamma_j(t).
\end{equation}
В силу (\ref{periodg}) функция $g$ растет на бесконечности не
быстрее, чем линейная функция, поэтому
$g(z,t)=\alpha(t)z+\beta(t)$. Итак,
\begin{equation}\label{h0}
h(z,t)=\sum_{j=1}^n\gamma_j(t)\zeta(z-a_j(t))+\alpha(t)z+\beta(t).
\end{equation}
Из условия $h(0,t)=0$ находим
\begin{equation}\label{beta0}
\beta(t)=\sum_{j=1}^n\gamma_j(t)\zeta(a_j(t)).
\end{equation}
При этом, из (\ref{periodg}) следует, что
\begin{equation}\label{periods}
\alpha(t)\omega_k(t)=-\dot{\omega}_k(t)-\eta_k(t)\sum_{j=1}^n\gamma_j(t),
\quad k=1,2.
\end{equation} Если мы предположим, что $\omega_1(t)\equiv
1$ (это условие не ограничивает общности, так как при замене
$\omega_k$ на пропорциональные величины получаем конформно
эквивалентные торы), то тогда из (\ref{periods}) при $k=1$
следует, что
\begin{equation}\label{alpha0}
\alpha(t)=-\eta_1(t)\sum_{j=1}^n\gamma_j(t).
\end{equation}
Окончательно из (\ref{h0}),  (\ref{beta0}) и (\ref{alpha0}) находим
\begin{equation}\label{h}
h(z,t)=\sum_{j=1}^n\gamma_j(t)[\zeta(z-a_j(t))+\zeta(a_j(t))-\eta_1(t)z].
\end{equation}
При этом, (\ref{periods}) при $k=2$ дает
$$
\dot{\omega}_2(t)=-\alpha(t)\omega_2(t)-\eta_2(t)\sum_{j=1}^n\gamma_j(t)=(\omega_2(t)\eta_1(t)-\eta_2(t))\sum_{j=1}^n\gamma_j(t),
$$
и с учетом равенства (\ref{etaomega}) получаем
\begin{equation}\label{period2}
\dot{\omega}_2(t)=2\pi i\sum_{j=1}^n\gamma_j(t).
\end{equation}

Таким образом, доказана

\begin{theorem}\label{loew-ell}
Функции $f(z,t)$ удовлетворяют уравнению в частных производных $$
\frac{\dot{f}(z,t)}{f'(z,t)}=h(z,t),$$ где  $h(z,t)$ имеет вид
$(\ref{h})$, в котором $\gamma_k(t)$ определены равенством
$(\ref{gammak})$, а  $D_k(t)$ описываются соотношениями
$(\ref{dk})$. При этом, период $\omega_1(t)$ функций $f(z,t)$
равен единице, а период $\omega_2(t)$ удовлетворяет
дифференциальному уравнению~$(\ref{period2})$.
\end{theorem}

  Теперь выведем систему дифференциальных уравнений для определения критических точек $a_l(t)$, $1\le
l\le n$. Для этого подсчитаем величину $\dot{f}'(a_l(t),t)$ двумя
различными способами. С одной стороны, из (\ref{fprime}) следует,
что
\begin{equation}\label{aldotprime1}
\dot{f}'(a_l(t),t)=-\dot{a}_l(t)D_l(t).
\end{equation}
С другой стороны, по теореме~\ref{loew-ell} имеем
$\dot{f}(z,t)=h(z,t)f'(z,t)$, поэтому
\begin{equation}\label{dotf}
\dot{f}(z,t)=c(t)\sum_{k=1}^n\gamma_k(t)\,[\zeta(z-a_k(t))+\zeta(a_k(t))-\eta_1(t)z]\,\frac{\prod\limits_{j=0}^n\sigma(z-a_j(t))}{\sigma^{n+1}(z)}
\end{equation}
и
\begin{multline}\label{dotfpr}
\dot{f}'(z,t)=c(t)\sum_{k=1}^n\gamma_k(t)\,\Bigl[\Bigl(\zeta(z-a_k(t))+\zeta(a_k(t))-\eta_1(t)z\Bigr)
\Bigl(\sum_{s=0}^n\zeta(z-a_s(t))-(n+1)\zeta(z)\Bigr)-\\-\mathfrak{P}(z-a_k(t))-\eta_1(t)\Bigr]\,\frac{\prod\limits_{j=0}^n\sigma(z-a_j(t))}{\sigma^{n+1}(z)}.
\end{multline}
При $z\to a_l(t)$ из (\ref{dotfpr}) получаем
\begin{multline}\label{aldotprime2}
\dot{f}'(a_l(t),t)=c(t)\Biggl[\sum_{k\neq
l}\gamma_k(t)\,\bigl[\zeta(a_l(t)-a_k(t))+\zeta(a_k(t))-\eta_1(t)a_l(t)\bigr]+\\+
\gamma_l(t)\,\Bigl(\sum_{s\neq l}\zeta(a_l(t)-a_s(t))-\eta_1(t)
a_l(t)-n\zeta(a_l(t))\Bigr)\Biggr] \,\frac{\prod\limits_{j\neq
l}\sigma(a_l(t)-a_j(t))}{\sigma^{n+1}(a_l(t))}\,.
\end{multline}

Приравнивая (\ref{aldotprime1}) и  (\ref{aldotprime2}), с учетом
(\ref{dk}) находим
\begin{multline}\label{al}
\dot{a}_l(t)=\sum_{k\neq
l}\gamma_k(t)\,\bigl[\zeta(a_k(t)-a_l(t))-\zeta(a_k(t))+\eta_1(t)a_l(t)\bigr]+\\
+\gamma_l(t)\,\Bigl(\sum_{s\neq l}\zeta(a_s(t)-a_l(t))+\eta_1(t)
a_l(t)+n\zeta(a_l(t))\Bigr),\quad 1\le l \le n.
\end{multline}

Напомним, что входящие в эти дифференциальные уравнения функции
соответствуют периодам $\omega_1(t)\equiv 1$ и $\omega_2(t)$, где
$\omega_2(t)$ определяется из дифференциального уравнения
(\ref{period2}).

Наконец, найдем  дифференциальное уравнение для определения
$c(t)$. Из (\ref{dotf}) следует, что в окрестности точки $z=0$
имеем разложение
$$
\dot{f}(z,t)=(-1)^{n}c(t)\sum_{k=1}^n\gamma_k(t)\,[\mathfrak{P}(a_k(t))+\eta_1(t)]
\prod\limits_{j=0}^n\sigma(a_j(t))\,
\frac{1}{z^n}+O\left(\frac{1}{z^{n-1}}\right),$$ поэтому
\begin{equation}\label{mixed1}
\dot{f}'(z,t)\sim(-1)^{n+1}n\,c(t)\sum_{k=1}^n\gamma_k(t)\,[\mathfrak{P}(a_k(t))+\eta_1(t)]
\prod\limits_{j=0}^n\sigma(a_j(t))\, \frac{1}{z^{n+1}},\quad z\to
0.
\end{equation}

С другой стороны,
$$
f'(z,t)=(-1)^{n+1}c(t)\prod_{k=0}^n\sigma(a_k(t))\,\frac{1}{z^{n+1}}+O\left(\frac{1}{z^{n}}\right),
$$
и, таким образом,
\begin{equation}\label{mixed2}
\dot{f}'(z,t)\sim (-1)^{n+1}\prod_{k=0}^n\sigma(a_k(t))
\,\Biggl(\dot{c}(t)+c(t)\sum_{j=0}^n\Bigl[\zeta(a_j(t))\dot{a}_j(t)
+\dot{\omega}_2(t)\frac{\partial\ln\sigma(a_j(t))}{\partial\omega_2}\Bigr]\Biggr)
\, \frac{1}{z^{n+1}},
\end{equation}
$z\to 0$.

Сравнивая (\ref{mixed1}) и (\ref{mixed2}), получаем
\begin{equation}\label{c}
\dot{c}(t)/c(t)=-\sum_{j=0}^n\Bigl[\zeta(a_j(t))\dot{a}_j(t)
+\dot{\omega}_2(t)\frac{\partial\ln\sigma(a_j(t))}{\partial\omega_2}\Bigr]+
n\sum_{k=1}^n\gamma_k(t)\,\bigl(\mathfrak{P}(a_k(t))+\eta_1(t)\bigr).
\end{equation}

Итак, справедлива

\begin{theorem}\label{odu} Параметры $a_l(t)$, величина $c(t)$ и период $\omega_2(t)$
удовлетворяют системе дифференциальных уравнений $(\ref{al})$,
$(\ref{c})$, $(\ref{period2})$, где
${\partial\ln\sigma(z)}/{\partial\omega_2}$ вычисляется по формуле
$(\ref{lnsom2})$, $\gamma_k(t)$ определены равенством
$(\ref{gammak})$, а $D_k(t)$ --- формулой~$(\ref{dk})$.
\end{theorem}

Как и в случае полиномов, мы получаем в процессе решения не одно
отображение, а целое семейство эллиптических функций,
униформизирующих семейство комплексных торов. Кроме того,
описанную в теореме~\ref{odu} систему можно рассматривать как
вариационные формулы для эллиптических функций, выражающие
вариации критических точек через вариации их образов.

\subsection{Пример. Униформизация семейства комплексных торов}

В качестве примера рассмотрим однопараметрическое семейство
функций, в котором начальный элемент
$f(z,0)=\mathfrak{P}^2(z)-4\mathfrak{P}(z)$, где
$\mathfrak{P}(z)=\mathfrak{P}(z;1,i)$ --- $\mathfrak{P}$-функция
Вейерштрасса с периодами $\omega_1=1$ и $\omega_2(0)=i$. Эта
функция представима в виде
$$f(z,0)=c(0)\int_{a_0(0)}^{z}\frac{\prod_{k=0}^4\sigma(t-a_k(0))}{\sigma^5(t)}\,dt,$$
где $\sigma(z)=\sigma(z;1,i)$, $a_0(0)=-(1+i)/2$, $a_1(0)=1/2$,
$a_2(0)=i/2$, $a_3(0)=-a_4(0)=0.5+$ $+i\,0.292496\ldots,
c(0)=-56.796445+i\,7.628085\ldots.$ Отметим, что
$\mathfrak{P}(a_3(0))=\mathfrak{P}(a_4(0))=2$. Образы точек $a_k(0)$
равны $A_0(0)=0,\ A_1(0)=19.767437\ldots,\ A_2(0)=74.768923\ldots,\
A_3(0)=A_4(0)=-4.$

Рассмотрим равномерное движение точек $A_k$ вдоль прямолинейных
отрезков при изменении параметра $t$ от $0$ до $1$ так, что в
результате движения точка $A_1$ сдвигается на $i$, точка $A_2$
--- на $(-i)$, точка $A_3$ --- на $(-1)$, точка $A_4$ --- на $1$.
После численного решения соответствующей задачи Коши получаем, что
функция $f(z,1)$ имеет вид
\begin{equation}\label{ex}
f(z,1)=c(1)\int_{a_0(1)}^{z}\frac{\prod_{k=0}^4\sigma(t-a_k(1))}{\sigma^5(t)}\,dt,
\end{equation}
 где $\sigma(z)$ имеет периоды $\omega_1=1$ и
$\omega_2(1)=i\,0.995555\ldots,$ а параметры в формуле (\ref{ex})
имеют вид $$a_0(1)=-(1+i)/2, a_1(1)=0.494370\ldots,
a_2(1)=0.003009\ldots+i\, 0.497778\ldots,$$
$$a_3(1)=0.494655+ i\,
0.289238\ldots, a_4(1)=-0.505345\ldots -i\, 0.289238\ldots,$$
$$c(1)=-58.544026+i\, 7.879114\ldots$$

Подсчет расположения проекций точек ветвления $A_j(1)$ по формуле
(\ref{ex}) с параметрами, найденными приближенно, дает совпадение
с нужными с точностью до 6-го знака после запятой.

\section{Комплексные торы в задачах, связанных с аппроксимациями Эрмита-Паде}

В этом разделе с помощью эллиптических функций мы исследуем задачу
о разложении Наттолла трехлистной римановой поверхности рода $1$,
связанном с некоторым абелевым интегралом на этой поверхности. Это
разложение имеет важное применение при исследовании диагональных
аппроксимаций Эрмита-Паде. Мы дает описание этого разложения в
случае, когда трехточечная конфигурация проекций точек ветвления
поверхности обладает зеркальной симмерией.

\subsection{Аппроксимации Эрмита-Паде и римановы поверхности.
Разложение Наттолла. }

Сначала напомним определение диагональных аппроксимаций
Эрмита-Паде~II. Пусть нам даны голоморфные в окрестности
бесконечности функции $f_j$, $1\le j\le m$. Их аппроксимаций
Эрмита-Паде типа II суть рациональные функции
\begin{equation}\label{pade}
\frac{Q_{nj}(z)}{P_n(z)}\,,\quad 1\le j\le m,
\end{equation}
удовлетворяющих условию:  $$
P_n(z)f_j(z)-Q_{nj}(z)=O(z^{-(n+1)}),\  z\to\infty, \quad 1\le
j\le m,
$$
причем $\mbox{\rm deg}\, P_n\le mn$.

Предположим, что функции $f_j$ аналитически продолжимы из
бесконечности по любому пути на плоскости, лежащему вне
фиксированного компакта $E$. Важной задачей является нахождение
максимальных областей сходимости аппроксимаций к заданной функции.

В случае одной функции (случай $m=1$) задача была решена
Шталем~\cite{stahl}, \cite{stahl1}. Он показал, что такая область
представляет собой внешность некоторого компакта $K$. Этот компакт
$K$ описывается с помощью ортогональных критических траекторий
квадратичного дифференциала, связанного с функцией Грина
внешности~$K$.

В случае нескольких функций вопрос остается открытым. Представляет
интерес исследование ситуации, когда множество $E$ является
конечным; в дальнейшем будем считать $E$ таковым.

Наттолл~\cite{nutall} предположил, что асимптотика аппроксимаций
Паде-Эрмита связана с $(m+1)$-листной компактной римановой
поверхностью $S$, накрывающей сферу Римана, точки ветвления
которой лежат над множеством $E\subset \mathbb{C}$. Обозначим
через $p$ проектирующее отображение $p:S\to
\overline{\mathbb{C}}$. Рассмотрим абелев интеграл $G$ на $S$,
который регулярен в любой точке $S$, кроме точек $P_0$,
$P_1,\ldots,P_m$, лежащих над бесконечно удаленной точкой. В
точках $P_j$ этот интеграл $G$ имеет асимптотику
\begin{equation}\label{dif_nat}
G(z)\sim \left\{
\begin{array}{ccc}
m\ln z, &z\to P_0,  &  \\[2mm]
- \ln z,  &z\to P_j,& 1\le j \le m.\\
\end{array}
\right.
\end{equation}
 Кроме того, все периоды $G$ чисто мнимые.
Отметим, что такой дифференциал определяется с точностью до
константы и $g=\Re G$ является однозначной гармонической функцией.
Для каждого $z\in \mathbb{C}\setminus E$ существует ровно $(m+1)$
точка поверхности $S$, лежащая над ней. Обозначим эти точки через
$z^{(0)}$, $z^{(1)},\ldots,z^{(m)}$. Введем также обозначения
$g_j(z):=g(z^{(j)})$. Понятно, что нумерацию точек $z^{(0)}$,
$z^{(1)},\ldots,z^{(m)}$ можно выбрать так, чтобы $g_0(z)\ge
g_1(z)\ge\ldots\ge g_m(z)$, $z\in \mathbb{C}\setminus E$. Заметим,
что в точках, где некоторые из $g_j$ совпадают, такая нумерация
точек определена неоднозначно. Рассмотрим множество точек $z$ на
плоскости, в которых значения $g_j(z)$ попарно различны. Тогда для
таких~$z$
$$
g_0(z)> g_1(z)>\ldots>g_m(z).
$$

Рассмотрим множества $$S_j:=\{P\in S: P=z^{(j)}\ \mbox{\rm для
точки}\ z =p(P)\}.$$ Назовем $S_j$ $j$-м листом поверхности $S$.
Отметим, что $S$ получается из листов $S_j$, $0\le j\le m$,
склеиванием вдоль некоторых кусочно-гладких кривых. Действительно,
пусть точка $P$ такова, что $z=p(P)\in \mathbb{C}\setminus E$,
$p^{-1}(z)=\{z_0,z_1\ldots,z_m\}$, причем для некоторых $j\neq k$
имеем $g(z_j)=g(z_k)$. Проекция $p:S\to \mathbb{C}$ локально
гомеоморфна в окрестности точек $z_j$ и $z_k$. Найдем такие малые
окрестности $U_j$ и $U_k$ точек $z_j$ и $z_k$ и окрестность $V$
точки $z$, что $p|_{U_j}$ отображает гомеоморфно $U_j$ на $V$ и
$p|_{U_k}$ отображает гомеоморфно  $U_k$ на $V$. Тогда множество
точек в окрестности $U_j$, которые лежат на границе $j$-го и $k$-го
листов,  это множество точек $Q$, для которых $g(Q)= g(p_j^{-1}\circ
p_k(Q))$. Таким образом, локально --- это множество нулей
гармонической функции $g-g\circ p_j^{-1}\circ p_k$. В силу связности
и компактности $S$ это множество либо совпадает со всей $S$, что
невозможно ввиду заданной асимптотики (\ref{dif_nat}) функции $G$ на
бесконечности, либо состоит из конечного числа аналитических дуг.
Отметим, что проекция граничных кривых $\partial S_j$ листов $S_j$
на сферу Римана также состоит из конечного числа аналитических дуг.

Даже когда $m=2$, и число точек множества $E$ невелико, ситуация
остается не исследованной. Рассмотрим, например задачу об
аппроксимациях Паде для пары функций $f_1=f$ $f_2=f^2$, где
$$
f(z)=\prod_{j=1}^3(z-a_j)^{\alpha_j},
$$
все точки $a_j$ различны, числа $2\alpha_j$ --- не целые и
$\sum_{j=1^3}{\alpha_j}=0$. Такая задача была поставлена
С.~П.~Суетиным (см. \cite{apt_tul}).

Рассмотрим риманову поверхность функции
\begin{equation}\label{surf3}
w=\sqrt[3]{(z-a_1)(z-a_2)(z-a_3)}
\end{equation}
и построим для нее функции $g_j$, $j=0$, $1$, $2$. В работе
А.~И.~Аптекарева и Д.~Н.~Тулякова~\cite{apt_tul} исследована
геометрическая структура множества
$$
\Gamma:=\{z\in \mathbb{C} \mid \exists j,k\in \{0,1,2\}: j\neq k\
\mbox{\rm и} \ g_j(z)=g_k(z)\}
$$
в случае, когда треугольник с вершинами $a_1$, $a_2$, $a_3$
достаточно близок к правильному. По гипотезе Натолла, множество
$$\gamma=\{z\in \mathbb{C} \mid
g_0(z)=g_1(z)\}$$ притягивает полюсы рациональных аппроксимаций
Паде--Эрмита. Также интересно множество $$\delta =\{z\in
\mathbb{C} \mid g_0(z)=g_1(z)\}.$$

Мы рассмотрим произвольное расположение точек $a_1$, $a_2$, $a_3$
и исследуем вопрос в случае равнобедренных треугольников
$\Delta(a_1,a_2,a_3)$ c вершинами в точках $a_1$, $a_2$, $a_3$
(рис.~5).

\subsection{Униформизация трехлистной поверхности}
Рассмотрим трехлистную риманову поверхность $R(a_1,a_2,a_3)$
функции (\ref{surf3}) над расширенной плоскостью (сферой Римана)
$\overline{\mathbb{C}}$ комплексного переменного $z$. Здесь $a_1$,
$a_2$ и $a_3$
--- три попарно различные точки комплексной плоскости, не лежащие на одной прямой.  Эта
поверхность имеет три листа и род $g=1$. Проведем через точки
$a_1$, $a_2$ и $a_3$ окружность и будем считать, что точки
занумерованы таким образом, что при обходе окружности против
часовой стрелки они встречаются в порядке $a_1$, $a_3$, $a_2$.
Рассмотрим дробно-линейное преобразование $T$ в $w$-плоскости,
переводящее точки $a_1$, $a_2$ и $a_3$ в точки
$\widetilde{a}_1=1$, $\widetilde{a}_2=e^{2\pi i/3}$,
$\widetilde{a}_3=e^{-2\pi i/3}$ --- корни кубические из единицы;
соответствующие точки располагаются в вершинах равностороннего
треугольника, вписанного в единичную окружность. Это
преобразование переводит внешность окружности $\sigma$ во
внутренность единичной окружности. При преобразовании $T$
поверхность $R(a_1,a_2,a_3)$ переходит в поверхность
$R(\widetilde{a}_1,\widetilde{a}_2,\widetilde{a}_3)$. Следует
также отметить, что при таком преобразовании бесконечно удаленная
точка переходит в точку
$$
z_0=\frac{\gamma-e^{-\pi i/3}}{\gamma-e^{\pi i/3}}\,,
$$
лежащую в единичном круге; здесь
$$
\gamma=\frac{a_3-a_2}{a_3-a_1}
$$
ангармоническое отношение точек $a_1$, $a_2$, $a_3$  и $\infty$.

\vskip 1 cm \hskip 0.7 cm
\unitlength 1mm 
\linethickness{0.4pt}
\ifx\plotpoint\undefined\newsavebox{\plotpoint}\fi 
\begin{picture}(138,64.25)(0,0)
\put(52.822,33.75){\line(0,1){.8846}}
\put(52.801,34.635){\line(0,1){.8826}}
\put(52.738,35.517){\line(0,1){.8787}}
\multiput(52.635,36.396)(-.028988,.174571){5}{\line(0,1){.174571}}
\multiput(52.49,37.269)(-.030967,.14418){6}{\line(0,1){.14418}}
\multiput(52.304,38.134)(-.032322,.122198){7}{\line(0,1){.122198}}
\multiput(52.078,38.989)(-.033275,.105476){8}{\line(0,1){.105476}}
\multiput(51.811,39.833)(-.030557,.083037){10}{\line(0,1){.083037}}
\multiput(51.506,40.663)(-.031296,.074099){11}{\line(0,1){.074099}}
\multiput(51.162,41.478)(-.031848,.066501){12}{\line(0,1){.066501}}
\multiput(50.779,42.276)(-.032251,.0599359){13}{\line(0,1){.0599359}}
\multiput(50.36,43.056)(-.0325299,.0541858){14}{\line(0,1){.0541858}}
\multiput(49.905,43.814)(-.0327045,.0490906){15}{\line(0,1){.0490906}}
\multiput(49.414,44.551)(-.0327895,.0445306){16}{\line(0,1){.0445306}}
\multiput(48.89,45.263)(-.0327963,.0404145){17}{\line(0,1){.0404145}}
\multiput(48.332,45.95)(-.0327339,.0366714){18}{\line(0,1){.0366714}}
\multiput(47.743,46.61)(-.0326096,.0332455){19}{\line(0,1){.0332455}}
\multiput(47.123,47.242)(-.0360324,.033436){18}{\line(-1,0){.0360324}}
\multiput(46.475,47.844)(-.0397737,.0335706){17}{\line(-1,0){.0397737}}
\multiput(45.798,48.414)(-.0438892,.0336433){16}{\line(-1,0){.0438892}}
\multiput(45.096,48.953)(-.0484499,.0336463){15}{\line(-1,0){.0484499}}
\multiput(44.369,49.457)(-.0535475,.0335701){14}{\line(-1,0){.0535475}}
\multiput(43.62,49.927)(-.0593019,.0334023){13}{\line(-1,0){.0593019}}
\multiput(42.849,50.362)(-.065873,.033127){12}{\line(-1,0){.065873}}
\multiput(42.058,50.759)(-.073481,.032721){11}{\line(-1,0){.073481}}
\multiput(41.25,51.119)(-.082431,.032154){10}{\line(-1,0){.082431}}
\multiput(40.426,51.441)(-.093168,.031383){9}{\line(-1,0){.093168}}
\multiput(39.587,51.723)(-.106358,.030341){8}{\line(-1,0){.106358}}
\multiput(38.736,51.966)(-.123047,.028925){7}{\line(-1,0){.123047}}
\multiput(37.875,52.168)(-.173978,.032354){5}{\line(-1,0){.173978}}
\multiput(37.005,52.33)(-.21913,.03018){4}{\line(-1,0){.21913}}
\put(36.129,52.451){\line(-1,0){.8812}}
\put(35.247,52.53){\line(-1,0){.884}}
\put(34.363,52.568){\line(-1,0){.8848}}
\put(33.479,52.564){\line(-1,0){.8836}}
\put(32.595,52.519){\line(-1,0){.8805}}
\multiput(31.714,52.432)(-.21887,-.03201){4}{\line(-1,0){.21887}}
\multiput(30.839,52.304)(-.144751,-.028177){6}{\line(-1,0){.144751}}
\multiput(29.97,52.135)(-.1228,-.029956){7}{\line(-1,0){.1228}}
\multiput(29.111,51.925)(-.106099,-.031232){8}{\line(-1,0){.106099}}
\multiput(28.262,51.676)(-.092901,-.032164){9}{\line(-1,0){.092901}}
\multiput(27.426,51.386)(-.082158,-.032845){10}{\line(-1,0){.082158}}
\multiput(26.604,51.058)(-.073204,-.033336){11}{\line(-1,0){.073204}}
\multiput(25.799,50.691)(-.065593,-.033678){12}{\line(-1,0){.065593}}
\multiput(25.012,50.287)(-.0548038,-.0314775){14}{\line(-1,0){.0548038}}
\multiput(24.245,49.846)(-.049713,-.0317504){15}{\line(-1,0){.049713}}
\multiput(23.499,49.37)(-.0451555,-.0319235){16}{\line(-1,0){.0451555}}
\multiput(22.777,48.859)(-.0410403,-.0320098){17}{\line(-1,0){.0410403}}
\multiput(22.079,48.315)(-.0372966,-.0320197){18}{\line(-1,0){.0372966}}
\multiput(21.408,47.739)(-.0357506,-.0337372){18}{\line(-1,0){.0357506}}
\multiput(20.764,47.131)(-.0323294,-.033518){19}{\line(0,-1){.033518}}
\multiput(20.15,46.494)(-.032425,-.0369448){18}{\line(0,-1){.0369448}}
\multiput(19.566,45.829)(-.032456,-.0406883){17}{\line(0,-1){.0406883}}
\multiput(19.014,45.138)(-.0324146,-.0448042){16}{\line(0,-1){.0448042}}
\multiput(18.496,44.421)(-.0322914,-.0493633){15}{\line(0,-1){.0493633}}
\multiput(18.011,43.68)(-.032074,-.0544569){14}{\line(0,-1){.0544569}}
\multiput(17.562,42.918)(-.0317468,-.0602044){13}{\line(0,-1){.0602044}}
\multiput(17.15,42.135)(-.031289,-.066766){12}{\line(0,-1){.066766}}
\multiput(16.774,41.334)(-.030673,-.074359){11}{\line(0,-1){.074359}}
\multiput(16.437,40.516)(-.033176,-.092545){9}{\line(0,-1){.092545}}
\multiput(16.138,39.683)(-.032389,-.105752){8}{\line(0,-1){.105752}}
\multiput(15.879,38.837)(-.031295,-.122465){7}{\line(0,-1){.122465}}
\multiput(15.66,37.98)(-.029756,-.144434){6}{\line(0,-1){.144434}}
\multiput(15.481,37.113)(-.027522,-.174808){5}{\line(0,-1){.174808}}
\put(15.344,36.239){\line(0,-1){.8795}}
\put(15.247,35.36){\line(0,-1){.8831}}
\put(15.193,34.477){\line(0,-1){.8847}}
\put(15.179,33.592){\line(0,-1){.8844}}
\put(15.207,32.708){\line(0,-1){.8821}}
\put(15.277,31.826){\line(0,-1){.8778}}
\multiput(15.388,30.948)(.030452,-.174321){5}{\line(0,-1){.174321}}
\multiput(15.541,30.076)(.032176,-.143915){6}{\line(0,-1){.143915}}
\multiput(15.734,29.213)(.033346,-.121923){7}{\line(0,-1){.121923}}
\multiput(15.967,28.359)(.030364,-.093505){9}{\line(0,-1){.093505}}
\multiput(16.24,27.518)(.031252,-.082777){10}{\line(0,-1){.082777}}
\multiput(16.553,26.69)(.031917,-.073834){11}{\line(0,-1){.073834}}
\multiput(16.904,25.878)(.032405,-.066231){12}{\line(0,-1){.066231}}
\multiput(17.293,25.083)(.0327528,-.0596631){13}{\line(0,-1){.0596631}}
\multiput(17.719,24.307)(.0329834,-.0539109){14}{\line(0,-1){.0539109}}
\multiput(18.18,23.553)(.0331153,-.0488144){15}{\line(0,-1){.0488144}}
\multiput(18.677,22.82)(.033162,-.0442539){16}{\line(0,-1){.0442539}}
\multiput(19.208,22.112)(.0331343,-.0401379){17}{\line(0,-1){.0401379}}
\multiput(19.771,21.43)(.0330405,-.0363954){18}{\line(0,-1){.0363954}}
\multiput(20.366,20.775)(.0328874,-.0329706){19}{\line(0,-1){.0329706}}
\multiput(20.99,20.148)(.0363118,-.0331324){18}{\line(1,0){.0363118}}
\multiput(21.644,19.552)(.040054,-.0332356){17}{\line(1,0){.040054}}
\multiput(22.325,18.987)(.04417,-.0332738){16}{\line(1,0){.04417}}
\multiput(23.032,18.455)(.0487306,-.0332385){15}{\line(1,0){.0487306}}
\multiput(23.763,17.956)(.0538274,-.0331196){14}{\line(1,0){.0538274}}
\multiput(24.516,17.492)(.0595801,-.0329035){13}{\line(1,0){.0595801}}
\multiput(25.291,17.065)(.066149,-.032573){12}{\line(1,0){.066149}}
\multiput(26.085,16.674)(.073753,-.032103){11}{\line(1,0){.073753}}
\multiput(26.896,16.321)(.082698,-.031461){10}{\line(1,0){.082698}}
\multiput(27.723,16.006)(.093428,-.0306){9}{\line(1,0){.093428}}
\multiput(28.564,15.731)(.121838,-.033654){7}{\line(1,0){.121838}}
\multiput(29.417,15.495)(.143833,-.032539){6}{\line(1,0){.143833}}
\multiput(30.28,15.3)(.174244,-.030892){5}{\line(1,0){.174244}}
\put(31.151,15.145){\line(1,0){.8775}}
\put(32.028,15.032){\line(1,0){.8819}}
\put(32.91,14.96){\line(1,0){.8843}}
\put(33.794,14.93){\line(1,0){.8847}}
\put(34.679,14.941){\line(1,0){.8832}}
\put(35.562,14.993){\line(1,0){.8798}}
\multiput(36.442,15.088)(.174877,.02708){5}{\line(1,0){.174877}}
\multiput(37.317,15.223)(.144509,.029391){6}{\line(1,0){.144509}}
\multiput(38.184,15.399)(.122544,.030986){7}{\line(1,0){.122544}}
\multiput(39.041,15.616)(.105833,.032122){8}{\line(1,0){.105833}}
\multiput(39.888,15.873)(.092628,.032942){9}{\line(1,0){.092628}}
\multiput(40.722,16.17)(.08188,.033533){10}{\line(1,0){.08188}}
\multiput(41.541,16.505)(.066845,.03112){12}{\line(1,0){.066845}}
\multiput(42.343,16.879)(.0602844,.0315946){13}{\line(1,0){.0602844}}
\multiput(43.126,17.289)(.0545377,.0319363){14}{\line(1,0){.0545377}}
\multiput(43.89,17.736)(.0494448,.0321665){15}{\line(1,0){.0494448}}
\multiput(44.632,18.219)(.044886,.0323013){16}{\line(1,0){.044886}}
\multiput(45.35,18.736)(.0407702,.0323531){17}{\line(1,0){.0407702}}
\multiput(46.043,19.286)(.0370266,.0323316){18}{\line(1,0){.0370266}}
\multiput(46.709,19.868)(.0335995,.0322446){19}{\line(1,0){.0335995}}
\multiput(47.348,20.48)(.032047,.0337881){19}{\line(0,1){.0337881}}
\multiput(47.957,21.122)(.0321138,.0372156){18}{\line(0,1){.0372156}}
\multiput(48.535,21.792)(.0321134,.0409592){17}{\line(0,1){.0409592}}
\multiput(49.081,22.488)(.0320375,.0450747){16}{\line(0,1){.0450747}}
\multiput(49.593,23.21)(.031876,.0496326){15}{\line(0,1){.0496326}}
\multiput(50.071,23.954)(.0316159,.0547241){14}{\line(0,1){.0547241}}
\multiput(50.514,24.72)(.0312405,.0604687){13}{\line(0,1){.0604687}}
\multiput(50.92,25.506)(.033521,.073119){11}{\line(0,1){.073119}}
\multiput(51.289,26.311)(.033053,.082075){10}{\line(0,1){.082075}}
\multiput(51.619,27.131)(.032399,.09282){9}{\line(0,1){.09282}}
\multiput(51.911,27.967)(.031501,.10602){8}{\line(0,1){.10602}}
\multiput(52.163,28.815)(.030266,.122724){7}{\line(0,1){.122724}}
\multiput(52.375,29.674)(.028543,.144679){6}{\line(0,1){.144679}}
\multiput(52.546,30.542)(.03257,.21879){4}{\line(0,1){.21879}}
\put(52.676,31.417){\line(0,1){.8803}}
\put(52.765,32.298){\line(0,1){1.4524}}
\put(47.75,46.5){\circle*{1.5}} \put(25,50.5){\circle*{1.5}}
\put(18,23.75){\circle*{1.5}}
\put(75.5,37){\vector(4,-1){.07}}\qbezier(58.25,36.75)(66.125,38.625)(75.5,37)
\put(134.447,36.75){\line(0,1){.8014}}
\put(134.427,37.551){\line(0,1){.7995}}
\put(134.369,38.351){\line(0,1){.7957}}
\multiput(134.271,39.147)(-.027222,.158004){5}{\line(0,1){.158004}}
\multiput(134.135,39.937)(-.029074,.130408){6}{\line(0,1){.130408}}
\multiput(133.961,40.719)(-.030338,.110432){7}{\line(0,1){.110432}}
\multiput(133.749,41.492)(-.031223,.095219){8}{\line(0,1){.095219}}
\multiput(133.499,42.254)(-.031845,.083186){9}{\line(0,1){.083186}}
\multiput(133.212,43.003)(-.032274,.073382){10}{\line(0,1){.073382}}
\multiput(132.889,43.736)(-.032556,.065202){11}{\line(0,1){.065202}}
\multiput(132.531,44.454)(-.03272,.058244){12}{\line(0,1){.058244}}
\multiput(132.139,45.153)(-.0327871,.0522277){13}{\line(0,1){.0522277}}
\multiput(131.712,45.832)(-.0327721,.046956){14}{\line(0,1){.046956}}
\multiput(131.254,46.489)(-.0326865,.0422831){15}{\line(0,1){.0422831}}
\multiput(130.763,47.123)(-.0325388,.0381002){16}{\line(0,1){.0381002}}
\multiput(130.243,47.733)(-.0323357,.0343241){17}{\line(0,1){.0343241}}
\multiput(129.693,48.316)(-.0339698,.0327077){17}{\line(-1,0){.0339698}}
\multiput(129.115,48.872)(-.0377434,.0329519){16}{\line(-1,0){.0377434}}
\multiput(128.512,49.4)(-.0419245,.0331452){15}{\line(-1,0){.0419245}}
\multiput(127.883,49.897)(-.0465962,.0332817){14}{\line(-1,0){.0465962}}
\multiput(127.23,50.363)(-.0518674,.0333541){13}{\line(-1,0){.0518674}}
\multiput(126.556,50.796)(-.057884,.033353){12}{\line(-1,0){.057884}}
\multiput(125.861,51.196)(-.064844,.033265){11}{\line(-1,0){.064844}}
\multiput(125.148,51.562)(-.073026,.033072){10}{\line(-1,0){.073026}}
\multiput(124.418,51.893)(-.082835,.032749){9}{\line(-1,0){.082835}}
\multiput(123.672,52.188)(-.094874,.032258){8}{\line(-1,0){.094874}}
\multiput(122.913,52.446)(-.110094,.031539){7}{\line(-1,0){.110094}}
\multiput(122.143,52.667)(-.130084,.030493){6}{\line(-1,0){.130084}}
\multiput(121.362,52.85)(-.157698,.028942){5}{\line(-1,0){.157698}}
\put(120.574,52.994){\line(-1,0){.7946}}
\put(119.779,53.1){\line(-1,0){.7988}}
\put(118.98,53.168){\line(-1,0){.8012}}
\put(118.179,53.196){\line(-1,0){.8016}}
\put(117.378,53.185){\line(-1,0){.8001}}
\put(116.577,53.135){\line(-1,0){.7967}}
\multiput(115.781,53.047)(-.19786,-.03187){4}{\line(-1,0){.19786}}
\multiput(114.989,52.919)(-.156861,-.033182){5}{\line(-1,0){.156861}}
\multiput(114.205,52.753)(-.110756,-.029133){7}{\line(-1,0){.110756}}
\multiput(113.43,52.549)(-.095554,-.030183){8}{\line(-1,0){.095554}}
\multiput(112.665,52.308)(-.083528,-.030937){9}{\line(-1,0){.083528}}
\multiput(111.914,52.029)(-.07373,-.031473){10}{\line(-1,0){.07373}}
\multiput(111.176,51.715)(-.065553,-.031844){11}{\line(-1,0){.065553}}
\multiput(110.455,51.364)(-.058597,-.032084){12}{\line(-1,0){.058597}}
\multiput(109.752,50.979)(-.0525818,-.0322161){13}{\line(-1,0){.0525818}}
\multiput(109.068,50.561)(-.0473103,-.0322586){14}{\line(-1,0){.0473103}}
\multiput(108.406,50.109)(-.0426367,-.0322239){15}{\line(-1,0){.0426367}}
\multiput(107.767,49.626)(-.0384524,-.0321217){16}{\line(-1,0){.0384524}}
\multiput(107.151,49.112)(-.0346744,-.0319598){17}{\line(-1,0){.0346744}}
\multiput(106.562,48.568)(-.0330758,-.0336114){17}{\line(0,-1){.0336114}}
\multiput(106,47.997)(-.0333612,-.0373822){16}{\line(0,-1){.0373822}}
\multiput(105.466,47.399)(-.0336,-.0415609){15}{\line(0,-1){.0415609}}
\multiput(104.962,46.775)(-.0315349,-.0431488){15}{\line(0,-1){.0431488}}
\multiput(104.489,46.128)(-.0314946,-.0478223){14}{\line(0,-1){.0478223}}
\multiput(104.048,45.459)(-.0313674,-.0530925){13}{\line(0,-1){.0530925}}
\multiput(103.64,44.768)(-.031138,-.059104){12}{\line(0,-1){.059104}}
\multiput(103.266,44.059)(-.030787,-.066056){11}{\line(0,-1){.066056}}
\multiput(102.928,43.333)(-.03365,-.082473){9}{\line(0,-1){.082473}}
\multiput(102.625,42.59)(-.03329,-.094516){8}{\line(0,-1){.094516}}
\multiput(102.359,41.834)(-.032737,-.109744){7}{\line(0,-1){.109744}}
\multiput(102.129,41.066)(-.031909,-.129744){6}{\line(0,-1){.129744}}
\multiput(101.938,40.287)(-.030659,-.157374){5}{\line(0,-1){.157374}}
\put(101.785,39.501){\line(0,-1){.7934}}
\put(101.67,38.707){\line(0,-1){.7981}}
\put(101.594,37.909){\line(0,-1){3.2007}}
\put(101.68,34.708){\line(0,-1){.7928}}
\multiput(101.799,33.916)(.031471,-.157213){5}{\line(0,-1){.157213}}
\multiput(101.956,33.13)(.032579,-.129577){6}{\line(0,-1){.129577}}
\multiput(102.152,32.352)(.033304,-.109574){7}{\line(0,-1){.109574}}
\multiput(102.385,31.585)(.030025,-.083861){9}{\line(0,-1){.083861}}
\multiput(102.655,30.83)(.030668,-.074068){10}{\line(0,-1){.074068}}
\multiput(102.962,30.09)(.031128,-.065896){11}{\line(0,-1){.065896}}
\multiput(103.304,29.365)(.031443,-.058943){12}{\line(0,-1){.058943}}
\multiput(103.682,28.657)(.0316413,-.0529297){13}{\line(0,-1){.0529297}}
\multiput(104.093,27.969)(.0317412,-.0476589){14}{\line(0,-1){.0476589}}
\multiput(104.537,27.302)(.0317575,-.0429853){15}{\line(0,-1){.0429853}}
\multiput(105.014,26.657)(.0317009,-.0388001){16}{\line(0,-1){.0388001}}
\multiput(105.521,26.037)(.0335539,-.0372093){16}{\line(0,-1){.0372093}}
\multiput(106.058,25.441)(.0332491,-.0334401){17}{\line(0,-1){.0334401}}
\multiput(106.623,24.873)(.034839,-.0317802){17}{\line(1,0){.034839}}
\multiput(107.215,24.333)(.0386179,-.0319226){16}{\line(1,0){.0386179}}
\multiput(107.833,23.822)(.0428027,-.0320031){15}{\line(1,0){.0428027}}
\multiput(108.475,23.342)(.0474764,-.0320137){14}{\line(1,0){.0474764}}
\multiput(109.14,22.894)(.0527476,-.031944){13}{\line(1,0){.0527476}}
\multiput(109.826,22.478)(.058762,-.03178){12}{\line(1,0){.058762}}
\multiput(110.531,22.097)(.065717,-.031505){11}{\line(1,0){.065717}}
\multiput(111.254,21.75)(.073891,-.031092){10}{\line(1,0){.073891}}
\multiput(111.993,21.439)(.083687,-.030505){9}{\line(1,0){.083687}}
\multiput(112.746,21.165)(.095709,-.029689){8}{\line(1,0){.095709}}
\multiput(113.511,20.927)(.129389,-.03332){6}{\line(1,0){.129389}}
\multiput(114.288,20.727)(.15703,-.032371){5}{\line(1,0){.15703}}
\multiput(115.073,20.566)(.19803,-.03085){4}{\line(1,0){.19803}}
\put(115.865,20.442){\line(1,0){.7972}}
\put(116.662,20.358){\line(1,0){.8004}}
\put(117.463,20.312){\line(1,0){.8016}}
\put(118.264,20.305){\line(1,0){.801}}
\put(119.065,20.338){\line(1,0){.7985}}
\put(119.864,20.409){\line(1,0){.794}}
\multiput(120.658,20.519)(.157547,.029757){5}{\line(1,0){.157547}}
\multiput(121.445,20.668)(.129925,.031165){6}{\line(1,0){.129925}}
\multiput(122.225,20.855)(.10993,.032108){7}{\line(1,0){.10993}}
\multiput(122.994,21.08)(.094706,.032748){8}{\line(1,0){.094706}}
\multiput(123.752,21.342)(.082664,.033177){9}{\line(1,0){.082664}}
\multiput(124.496,21.64)(.072855,.033449){10}{\line(1,0){.072855}}
\multiput(125.225,21.975)(.064671,.033599){11}{\line(1,0){.064671}}
\multiput(125.936,22.344)(.057711,.033651){12}{\line(1,0){.057711}}
\multiput(126.629,22.748)(.0516944,.0336217){13}{\line(1,0){.0516944}}
\multiput(127.301,23.185)(.0464236,.0335221){14}{\line(1,0){.0464236}}
\multiput(127.951,23.655)(.0417527,.0333614){15}{\line(1,0){.0417527}}
\multiput(128.577,24.155)(.0375726,.0331465){16}{\line(1,0){.0375726}}
\multiput(129.178,24.685)(.0338003,.0328828){17}{\line(1,0){.0338003}}
\multiput(129.753,25.244)(.0321579,.0344907){17}{\line(0,1){.0344907}}
\multiput(130.299,25.831)(.0323415,.0382678){16}{\line(0,1){.0382678}}
\multiput(130.817,26.443)(.0324676,.0424515){15}{\line(0,1){.0424515}}
\multiput(131.304,27.08)(.032529,.0471248){14}{\line(0,1){.0471248}}
\multiput(131.759,27.739)(.0325167,.0523964){13}{\line(0,1){.0523964}}
\multiput(132.182,28.421)(.032419,.058412){12}{\line(0,1){.058412}}
\multiput(132.571,29.122)(.032219,.06537){11}{\line(0,1){.06537}}
\multiput(132.925,29.841)(.031895,.073548){10}{\line(0,1){.073548}}
\multiput(133.244,30.576)(.031415,.08335){9}{\line(0,1){.08335}}
\multiput(133.527,31.326)(.03073,.095379){8}{\line(0,1){.095379}}
\multiput(133.773,32.089)(.029767,.110587){7}{\line(0,1){.110587}}
\multiput(133.981,32.863)(.0284,.130557){6}{\line(0,1){.130557}}
\multiput(134.152,33.647)(.03301,.19768){4}{\line(0,1){.19768}}
\put(134.284,34.437){\line(0,1){.7962}}
\put(134.377,35.234){\line(0,1){1.5163}}
\put(50.75,49){\makebox(0,0)[cc]{$a_1$}}
\put(23,53.25){\makebox(0,0)[cc]{$a_3$}}
\put(12.5,24.25){\makebox(0,0)[cc]{$a_2$}}
\put(138,38.25){\makebox(0,0)[cc]{$\tilde{a}_1$}}
\put(108.25,55){\makebox(0,0)[cc]{$\tilde{a}_2$}}
\put(100,22){\makebox(0,0)[cc]{$\tilde{a}_3$}}
\put(33.25,64.25){\makebox(0,0)[cc]{$\infty$}}
\put(134,37.5){\circle*{1.5}} \put(109.75,50.75){\circle*{1.5}}
\put(26,63.25){\circle*{1.5}} \put(115,44.5){\circle*{1.5}}
\put(109.5,23){\circle*{1.5}}
\put(118.5,44){\makebox(0,0)[cc]{$z_0$}}
\put(71.75,7){\makebox(0,0)[cc]{Рис.~5}}
\multiput(17.75,23.75)(.04362170088,.03372434018){682}{\line(1,0){.04362170088}}
\multiput(47.5,46.75)(-.203125,.033482143){112}{\line(-1,0){.203125}}
\multiput(24.75,50.5)(-.033653846,-.127403846){208}{\line(0,-1){.127403846}}
\put(34.,45.3){\makebox(0,0)[cc]{$\Delta(a_1,a_2,a_3)$}}
\end{picture}

Таким образом, задача униформизации поверхности $R(a_1,a_2,a_3)$
сводится к аналогичной задачи для поверхности
$R_0:=R(\widetilde{a}_1,\widetilde{a}_2,\widetilde{a}_3)$, которой
соответствует функция $$w=\sqrt[3]{z^3-1}\,.$$ Обозначим через
$A_1$, $A_2$ и $A_3$ точки поверхности $R_0$, лежащие над
$\widetilde{a}_1$, $\widetilde{a}_2$ и $\widetilde{a}_3$.

Поскольку род поверхности $R_0$ равен $1$ (параболический случай),
в качестве универсального накрытия  $R_0$ можно взять комплексную
плоскость $\mathbb{C}$. Универсальное накрытие
$\pi:{\mathbb{C}}\to R_0$ осуществляется некоторой
двоякопериодической функцией $\pi$  с периодами $\omega_1$,
$\omega_2$. Мы можем считать, что $\pi(0)=\widetilde{a}_1$; этого
можно добиться сдвигом плоскости.

У римановой поверхности $R_0$ имеется нетривиальная группа
(голоморфных) автоморфизмов $f$, сохраняющих проекции точек, т.~е.
таких $f:R_0\to R_0$, что $p\circ f=p$. В частности, автоморфизмами
являются непрерывные функции $f$, переставляющие листы. Обозначим
через $Aut(R_0)$ группу таких отображений. (В число элементов
$Aut(R_0)$ естественно включить и тождественный автоморфизм, хотя он
листы и не переставляет!) Заметим, что точки $A_1$, $A_2$ и $A_3$
являются неподвижными точками для $f\in Aut(R_0)$.

Любой автоморфизм $f\in Aut(R_0)$ обладает поднятиями
$\widetilde{f}$ на универсальное накрытие. Таким образом, существует
$\widetilde{f}:\mathbb{C}\to \mathbb{C}$ такое, что
$f\circ\pi=\pi\circ\widetilde{f}$, т.~е. коммутативна диаграмма
\[
\begin{CD}
\mathbb{C} @>\widetilde{f}>> \mathbb{C} \\
{\pi}@VVV @VVV\pi \\
R_0 @>f>> R_0
\end{CD}
\]
Поскольку $\pi(0)=\widetilde{a}_1$ и
$f(\widetilde{a}_1)=\widetilde{a}_1$,  можно считать, что
$\widetilde{f}(0)=0$. Последнее равенство определяет $\widetilde{f}$
единственным образом.

Голоморфный автоморфизм $\widetilde{f}$ плоскости $\mathbb{C}$
является линейным отображением, следовательно,
$\widetilde{f}(z)=bz$, $z\in \mathbb{C}$, для некоторого $b\neq
0$. Анализ локального поведения автоморфизмов в окрестности
неподвижных точек $A_1$, $A_2$ и $A_3$, показывает, что подгруппа
автоморфизмов $\widetilde{f}$ плоскости, соответствующая группе
автоморфизмов, переставляющих листы $R_0$ и удовлетворяющих
условию $\widetilde{f}(0)=0$, является циклической группой,
порожденной поворотом на $2\pi/3$.

Теперь рассмотрим множество точек $\pi^{-1}(A_1)$. Нетрудно
показать, что это множество является решеткой
$\mbox{\boldmath$\omega$}$ на плоскости, порожденной двумя линейно
независимыми векторами $\omega_1$ и $\omega_2$. Из определения этой
решетки  следует, что
 $\mbox{\boldmath$\omega$}$
инвариантна относительно поворотов на углы, кратные $2\pi/3$, вокруг
начала координат. Аналогично показывается, что эта решетка
инвариантна относительно поворотов вокруг любой точки плоскости,
соответствующей $A_1$, $A_2$ и $A_3$. Отсюда следует

\begin{lemma}\label{hexagon}
Решетка  $\mbox{\boldmath$\omega$}$ располагается в вершинах
некоторой триангуляции плоскости правильными треугольниками. При
этом, если решетка инвариантна относительно поворота на угол
$2\pi/3$ относительно некоторой точки плоскости, то либо эта точка
является вершиной триангуляции, либо она является центром некоторого
треугольника триангуляции.
\end{lemma}

Поскольку отображение $\pi$ определяется с точностью до линейного
автоморфизма плоскости, можно считать, что образующими решетки
являются числа
\begin{equation}\label{per_sqrt2}
\omega_1=\sqrt{3},\quad \omega_2=\sqrt{3}e^{\pi i/3}.
\end{equation}

Точкам решетки $\mbox{\boldmath$\omega$}$ при отображении $\pi$
соответствует точка $A_1$. Из леммы~\ref{hexagon} следует, что
центру треугольника $T_1$ с вершинами $0$, $\sqrt{3}$,
$\sqrt{3}e^{i\pi/3}$ и эквивалентным ему точкам по модулю решетки
при отображении $\pi$ соответствует одна из точек $A_2$, $A_3$, а
центру треугольника $T_2$ с вершинами $\sqrt{3}$,
$\sqrt{3}e^{i\pi/3}$, $\sqrt{3}(1+e^{i\pi/3})$ и эквивалентным ему
точкам
--- другая. Пусть, для определенности центру $T_1$ соответствует
$\widetilde{a}_2$, а центру $T_2$ --- $\widetilde{a}_3$.
Треугольники $T_1$ и $T_2$ образуют фундаментальный параллелограмм
функции $\pi$, униформизирующей поверхность $R_0$. Обозначим через
$B$ центр треугольника $T_1$, а через $C$ --- центр треугольника
$T_2$. В силу двояко-периодичности униформизирующей функции $\pi$
 будем обозначать через $B$ и $C$ также точки, эквивалентные
им относительно решетки $\mbox{\boldmath$\omega$}$ (рис.~6).

\hskip -2.9 cm
\unitlength 1mm 
\linethickness{0.4pt}
\ifx\plotpoint\undefined\newsavebox{\plotpoint}\fi 
\begin{picture}(190.75,99.25)(0,0)
\put(47.75,20.5){\vector(0,1){78.75}}
\multiput(47.75,47.5)(.0336938436,.0574043261){601}{\line(0,1){.0574043261}}
\multiput(87.75,47.5)(.0336938436,.0574043261){601}{\line(0,1){.0574043261}}
\multiput(68,13.25)(.0336938436,.0574043261){601}{\line(0,1){.0574043261}}
\multiput(47.75,47.5)(-.0336938436,.0574043261){601}{\line(0,1){.0574043261}}
\multiput(67.5,13.5)(-.0336938436,.0574043261){601}{\line(0,1){.0574043261}}
\multiput(87.75,47.5)(-.0337268128,.0581787521){593}{\line(0,1){.0581787521}}
\put(67.75,81.75){\line(-1,0){40.5}}
\put(107.75,81.75){\line(-1,0){40.5}}
\multiput(47.68,47.68)(.0571429,.0321429){15}{\line(1,0){.0571429}}
\multiput(49.394,48.644)(.0571429,.0321429){15}{\line(1,0){.0571429}}
\multiput(51.108,49.608)(.0571429,.0321429){15}{\line(1,0){.0571429}}
\multiput(52.823,50.573)(.0571429,.0321429){15}{\line(1,0){.0571429}}
\multiput(54.537,51.537)(.0571429,.0321429){15}{\line(1,0){.0571429}}
\multiput(56.251,52.501)(.0571429,.0321429){15}{\line(1,0){.0571429}}
\multiput(57.965,53.465)(.0571429,.0321429){15}{\line(1,0){.0571429}}
\multiput(59.68,54.43)(.0571429,.0321429){15}{\line(1,0){.0571429}}
\multiput(61.394,55.394)(.0571429,.0321429){15}{\line(1,0){.0571429}}
\multiput(63.108,56.358)(.0571429,.0321429){15}{\line(1,0){.0571429}}
\multiput(64.823,57.323)(.0571429,.0321429){15}{\line(1,0){.0571429}}
\multiput(66.537,58.287)(.0571429,.0321429){15}{\line(1,0){.0571429}}
\multiput(68.251,59.251)(.0571429,.0321429){15}{\line(1,0){.0571429}}
\multiput(69.965,60.215)(.0571429,.0321429){15}{\line(1,0){.0571429}}
\multiput(71.68,61.18)(.0571429,.0321429){15}{\line(1,0){.0571429}}
\multiput(73.394,62.144)(.0571429,.0321429){15}{\line(1,0){.0571429}}
\multiput(75.108,63.108)(.0571429,.0321429){15}{\line(1,0){.0571429}}
\multiput(76.823,64.073)(.0571429,.0321429){15}{\line(1,0){.0571429}}
\multiput(78.537,65.037)(.0571429,.0321429){15}{\line(1,0){.0571429}}
\multiput(80.251,66.001)(.0571429,.0321429){15}{\line(1,0){.0571429}}
\multiput(81.965,66.965)(.0571429,.0321429){15}{\line(1,0){.0571429}}
\multiput(83.68,67.93)(.0571429,.0321429){15}{\line(1,0){.0571429}}
\multiput(85.394,68.894)(.0571429,.0321429){15}{\line(1,0){.0571429}}
\multiput(87.108,69.858)(.0571429,.0321429){15}{\line(1,0){.0571429}}
\multiput(88.823,70.823)(.0571429,.0321429){15}{\line(1,0){.0571429}}
\multiput(90.537,71.787)(.0571429,.0321429){15}{\line(1,0){.0571429}}
\multiput(92.251,72.751)(.0571429,.0321429){15}{\line(1,0){.0571429}}
\multiput(93.965,73.715)(.0571429,.0321429){15}{\line(1,0){.0571429}}
\multiput(95.68,74.68)(.0571429,.0321429){15}{\line(1,0){.0571429}}
\multiput(97.394,75.644)(.0571429,.0321429){15}{\line(1,0){.0571429}}
\multiput(99.108,76.608)(.0571429,.0321429){15}{\line(1,0){.0571429}}
\multiput(100.823,77.573)(.0571429,.0321429){15}{\line(1,0){.0571429}}
\multiput(102.537,78.537)(.0571429,.0321429){15}{\line(1,0){.0571429}}
\multiput(104.251,79.501)(.0571429,.0321429){15}{\line(1,0){.0571429}}
\multiput(105.965,80.465)(.0571429,.0321429){15}{\line(1,0){.0571429}}
\multiput(27.68,81.68)(.0571429,-.032381){15}{\line(1,0){.0571429}}
\multiput(29.394,80.708)(.0571429,-.032381){15}{\line(1,0){.0571429}}
\multiput(31.108,79.737)(.0571429,-.032381){15}{\line(1,0){.0571429}}
\multiput(32.823,78.765)(.0571429,-.032381){15}{\line(1,0){.0571429}}
\multiput(34.537,77.794)(.0571429,-.032381){15}{\line(1,0){.0571429}}
\multiput(36.251,76.823)(.0571429,-.032381){15}{\line(1,0){.0571429}}
\multiput(37.965,75.851)(.0571429,-.032381){15}{\line(1,0){.0571429}}
\multiput(39.68,74.88)(.0571429,-.032381){15}{\line(1,0){.0571429}}
\multiput(41.394,73.908)(.0571429,-.032381){15}{\line(1,0){.0571429}}
\multiput(43.108,72.937)(.0571429,-.032381){15}{\line(1,0){.0571429}}
\multiput(44.823,71.965)(.0571429,-.032381){15}{\line(1,0){.0571429}}
\multiput(46.537,70.994)(.0571429,-.032381){15}{\line(1,0){.0571429}}
\multiput(48.251,70.023)(.0571429,-.032381){15}{\line(1,0){.0571429}}
\multiput(49.965,69.051)(.0571429,-.032381){15}{\line(1,0){.0571429}}
\multiput(51.68,68.08)(.0571429,-.032381){15}{\line(1,0){.0571429}}
\multiput(53.394,67.108)(.0571429,-.032381){15}{\line(1,0){.0571429}}
\multiput(55.108,66.137)(.0571429,-.032381){15}{\line(1,0){.0571429}}
\multiput(56.823,65.165)(.0571429,-.032381){15}{\line(1,0){.0571429}}
\multiput(58.537,64.194)(.0571429,-.032381){15}{\line(1,0){.0571429}}
\multiput(60.251,63.223)(.0571429,-.032381){15}{\line(1,0){.0571429}}
\multiput(61.965,62.251)(.0571429,-.032381){15}{\line(1,0){.0571429}}
\multiput(63.68,61.28)(.0571429,-.032381){15}{\line(1,0){.0571429}}
\multiput(65.394,60.308)(.0571429,-.032381){15}{\line(1,0){.0571429}}
\multiput(67.108,59.337)(.0571429,-.032381){15}{\line(1,0){.0571429}}
\multiput(68.823,58.365)(.0571429,-.032381){15}{\line(1,0){.0571429}}
\multiput(70.537,57.394)(.0571429,-.032381){15}{\line(1,0){.0571429}}
\multiput(72.251,56.423)(.0571429,-.032381){15}{\line(1,0){.0571429}}
\multiput(73.965,55.451)(.0571429,-.032381){15}{\line(1,0){.0571429}}
\multiput(75.68,54.48)(.0571429,-.032381){15}{\line(1,0){.0571429}}
\multiput(77.394,53.508)(.0571429,-.032381){15}{\line(1,0){.0571429}}
\multiput(79.108,52.537)(.0571429,-.032381){15}{\line(1,0){.0571429}}
\multiput(80.823,51.565)(.0571429,-.032381){15}{\line(1,0){.0571429}}
\multiput(82.537,50.594)(.0571429,-.032381){15}{\line(1,0){.0571429}}
\multiput(84.251,49.623)(.0571429,-.032381){15}{\line(1,0){.0571429}}
\multiput(85.965,48.651)(.0571429,-.032381){15}{\line(1,0){.0571429}}
\put(67.5,58.75){\circle*{1.5}} \put(67.5,36){\circle*{1.5}}
\put(67.43,59.18){\line(0,-1){.9583}}
\put(67.43,57.263){\line(0,-1){.9583}}
\put(67.43,55.346){\line(0,-1){.9583}}
\put(67.43,53.43){\line(0,-1){.9583}}
\put(67.43,51.513){\line(0,-1){.9583}}
\put(67.43,49.596){\line(0,-1){.9583}}
\put(67.43,47.68){\line(0,-1){.9583}}
\put(67.43,45.763){\line(0,-1){.9583}}
\put(67.43,43.846){\line(0,-1){.9583}}
\put(67.43,41.93){\line(0,-1){.9583}}
\put(67.43,40.013){\line(0,-1){.9583}}
\put(67.43,38.096){\line(0,-1){.9583}}
\multiput(27.68,59.18)(.0595238,.0334821){14}{\line(1,0){.0595238}}
\multiput(29.346,60.117)(.0595238,.0334821){14}{\line(1,0){.0595238}}
\multiput(31.013,61.055)(.0595238,.0334821){14}{\line(1,0){.0595238}}
\multiput(32.68,61.992)(.0595238,.0334821){14}{\line(1,0){.0595238}}
\multiput(34.346,62.93)(.0595238,.0334821){14}{\line(1,0){.0595238}}
\multiput(36.013,63.867)(.0595238,.0334821){14}{\line(1,0){.0595238}}
\multiput(37.68,64.805)(.0595238,.0334821){14}{\line(1,0){.0595238}}
\multiput(39.346,65.742)(.0595238,.0334821){14}{\line(1,0){.0595238}}
\multiput(41.013,66.68)(.0595238,.0334821){14}{\line(1,0){.0595238}}
\multiput(42.68,67.617)(.0595238,.0334821){14}{\line(1,0){.0595238}}
\multiput(44.346,68.555)(.0595238,.0334821){14}{\line(1,0){.0595238}}
\multiput(46.013,69.492)(.0595238,.0334821){14}{\line(1,0){.0595238}}
\put(27.93,59.18){\line(0,-1){.9891}}
\put(27.93,57.201){\line(0,-1){.9891}}
\put(27.93,55.223){\line(0,-1){.9891}}
\put(27.93,53.245){\line(0,-1){.9891}}
\put(27.93,51.267){\line(0,-1){.9891}}
\put(27.93,49.288){\line(0,-1){.9891}}
\put(27.93,47.31){\line(0,-1){.9891}}
\put(27.93,45.332){\line(0,-1){.9891}}
\put(27.93,43.354){\line(0,-1){.9891}}
\put(27.93,41.375){\line(0,-1){.9891}}
\put(27.93,39.397){\line(0,-1){.9891}}
\put(27.93,37.419){\line(0,-1){.9891}}
\multiput(47.68,24.93)(.0602679,.0327381){14}{\line(1,0){.0602679}}
\multiput(49.367,25.846)(.0602679,.0327381){14}{\line(1,0){.0602679}}
\multiput(51.055,26.763)(.0602679,.0327381){14}{\line(1,0){.0602679}}
\multiput(52.742,27.68)(.0602679,.0327381){14}{\line(1,0){.0602679}}
\multiput(54.43,28.596)(.0602679,.0327381){14}{\line(1,0){.0602679}}
\multiput(56.117,29.513)(.0602679,.0327381){14}{\line(1,0){.0602679}}
\multiput(57.805,30.43)(.0602679,.0327381){14}{\line(1,0){.0602679}}
\multiput(59.492,31.346)(.0602679,.0327381){14}{\line(1,0){.0602679}}
\multiput(61.18,32.263)(.0602679,.0327381){14}{\line(1,0){.0602679}}
\multiput(62.867,33.18)(.0602679,.0327381){14}{\line(1,0){.0602679}}
\multiput(64.555,34.096)(.0602679,.0327381){14}{\line(1,0){.0602679}}
\multiput(66.242,35.013)(.0602679,.0327381){14}{\line(1,0){.0602679}}
\multiput(27.93,36.18)(.0587798,-.0327381){14}{\line(1,0){.0587798}}
\multiput(29.576,35.263)(.0587798,-.0327381){14}{\line(1,0){.0587798}}
\multiput(31.221,34.346)(.0587798,-.0327381){14}{\line(1,0){.0587798}}
\multiput(32.867,33.43)(.0587798,-.0327381){14}{\line(1,0){.0587798}}
\multiput(34.513,32.513)(.0587798,-.0327381){14}{\line(1,0){.0587798}}
\multiput(36.159,31.596)(.0587798,-.0327381){14}{\line(1,0){.0587798}}
\multiput(37.805,30.68)(.0587798,-.0327381){14}{\line(1,0){.0587798}}
\multiput(39.451,29.763)(.0587798,-.0327381){14}{\line(1,0){.0587798}}
\multiput(41.096,28.846)(.0587798,-.0327381){14}{\line(1,0){.0587798}}
\multiput(42.742,27.93)(.0587798,-.0327381){14}{\line(1,0){.0587798}}
\multiput(44.388,27.013)(.0587798,-.0327381){14}{\line(1,0){.0587798}}
\multiput(46.034,26.096)(.0587798,-.0327381){14}{\line(1,0){.0587798}}
\put(87.5,70.25){\circle*{1.5}} \put(47.5,70.5){\circle*{1.5}}
\put(28,59.5){\circle*{1.5}} \put(28.25,36){\circle*{1.5}}
\put(47.5,25.25){\circle*{1.5}} \put(87.75,47.5){\circle*{1.5}}
\put(67.75,81.4){\circle*{1.5}} \put(107.75,81.5){\circle*{1.5}}
\put(28,81.5){\circle*{1.5}} \put(47.75,47.5){\circle*{1.5}}
\put(25,47.5){\vector(1,0){86}}
\multiput(140,92)(-.0337259101,-.0524625268){467}{\line(0,-1){.0524625268}}
\multiput(140,48.5)(-.0337259101,-.0524625268){467}{\line(0,-1){.0524625268}}
\put(188.25,48){\line(-1,0){49}}
\put(188.25,91.5){\line(-1,0){49}}
\put(189.25,88.25){\line(-1,0){2.75}}
\put(190.25,85.25){\line(-1,0){3.25}}
\put(124.25,67){\line(1,0){47.75}}
\put(124.25,24.5){\line(1,0){47.75}}
\put(125.25,64.5){\line(1,0){47.75}}
\put(126.75,61.25){\line(1,0){47.75}}
\multiput(172,67)(.0337136929,.0497925311){482}{\line(0,1){.0497925311}}
\multiput(172,24.5)(.0337136929,.0497925311){482}{\line(0,1){.0497925311}}
\multiput(173,64.5)(.0337136929,.0497925311){482}{\line(0,1){.0497925311}}
\multiput(174.5,61.25)(.0337136929,.0497925311){482}{\line(0,1){.0497925311}}
\multiput(125.25,64.5)(.03333333,.05555556){45}{\line(0,1){.05555556}}
\multiput(126.75,61.25)(.03333333,.07777778){45}{\line(0,1){.07777778}}
\qbezier(155.5,78.5)(158.75,79.625)(166,78.25)
\qbezier(149,84.68)(148.75,80.5)(155.5,78.75)
\qbezier(155.5,78.75)(150.125,78.375)(148.25,74.5)
\put(166,78.25){\circle*{1.5}} \put(148.5,74.5){\circle*{1.5}}
\put(149,84.25){\circle*{1.5}}
\put(168.25,71.5){\makebox(0,0)[cc]{$A_1$}}
\put(147,87.5){\makebox(0,0)[cc]{$A_2$}}
\put(144.25,72.75){\makebox(0,0)[cc]{$A_3$}}
\put(177.25,83.5){\makebox(0,0)[cc]{$R_0$}}
\put(168.25,29.5){\makebox(0,0)[cc]{$\widetilde{a}_1$}}
\put(144,31.5){\makebox(0,0)[cc]{$\widetilde{a}_3$}}
\put(147.5,42.75){\makebox(0,0)[cc]{$\widetilde{a}_2$}}
\put(166,34.25){\circle*{1.5}} \put(148.5,28.5){\circle*{1.5}}
\put(149,38.25){\circle*{1.5}}
\put(167,39.75){\makebox(0,0)[cc]{$\overline{\mathbb{C}}$}}
\put(155.25,58.5){\vector(0,-1){8.5}}
\put(153.25,53.75){\makebox(0,0)[cc]{$p$}}
\put(67.25,87.5){\makebox(0,0)[cc]{$A$}}
\put(26,85.5){\makebox(0,0)[cc]{$A$}}
\put(108.5,85.75){\makebox(0,0)[cc]{$A$}}
\put(90.5,43){\makebox(0,0)[cc]{$A$}}
\put(44.5,44.75){\makebox(0,0)[cc]{$A$}}
\put(68,13){\circle*{1.5}} \put(73,14.25){\makebox(0,0)[cc]{$A$}}
\put(67.25,64){\makebox(0,0)[cc]{$B$}}
\put(86.25,74.25){\makebox(0,0)[cc]{$C$}}
\put(72.75,37){\makebox(0,0)[cc]{$C$}}
\put(50.75,72.75){\makebox(0,0)[cc]{$C$}}
\put(25.75,63.75){\makebox(0,0)[cc]{$B$}}
\put(27.75,31.5){\makebox(0,0)[cc]{$C$}}
\put(45,23.75){\makebox(0,0)[cc]{$B$}}
\put(70,71.75){\makebox(0,0)[cc]{$T_1$}}
\put(88,62.75){\makebox(0,0)[cc]{$T_2$}}
\multiput(47.43,47.43)(.0613354,-.0333851){14}{\line(1,0){.0613354}}
\multiput(49.147,46.495)(.0613354,-.0333851){14}{\line(1,0){.0613354}}
\multiput(50.864,45.56)(.0613354,-.0333851){14}{\line(1,0){.0613354}}
\multiput(52.582,44.625)(.0613354,-.0333851){14}{\line(1,0){.0613354}}
\multiput(54.299,43.691)(.0613354,-.0333851){14}{\line(1,0){.0613354}}
\multiput(56.017,42.756)(.0613354,-.0333851){14}{\line(1,0){.0613354}}
\multiput(57.734,41.821)(.0613354,-.0333851){14}{\line(1,0){.0613354}}
\multiput(59.451,40.886)(.0613354,-.0333851){14}{\line(1,0){.0613354}}
\multiput(61.169,39.951)(.0613354,-.0333851){14}{\line(1,0){.0613354}}
\multiput(62.886,39.017)(.0613354,-.0333851){14}{\line(1,0){.0613354}}
\multiput(64.604,38.082)(.0613354,-.0333851){14}{\line(1,0){.0613354}}
\multiput(66.321,37.147)(.0613354,-.0333851){14}{\line(1,0){.0613354}}
\put(100.25,7.5){\makebox(0,0)[cc]{Рис.~6}}
\end{picture}

Опишем  универсальное накрытие в явном (аналитическом) виде.
Рассмотрим функцию, конформно отображающую треугольник
$\Delta=ABC$ с вершинами в точках $0$, $e^{-i\pi/6}$, $e^{i\pi/6}$
на единичный круг с соответствием точек:
$$
0\mapsto 1, \quad e^{-i\pi/6}\mapsto  e^{i2\pi/3}, \quad
e^{i\pi/6}\mapsto  e^{-i2\pi/3}.
$$
Продолжая эту функцию по принципу симметрии, получим отображение,
накрывающее риманову поверхность $R_0$ комплексной плоскостью.
Заметим, что центр $z_1=\sqrt{3}/3$ треугольника $\Delta$
переходит в начало координат. В начало координат переходят также и
все точки, эквивалентные $z_1$ по модулю решетки. Точка $(-z_1)$
получается из $z_1$ путем нечетного числа отражений относительно
сторон треугольников триангуляции, поэтому $(-z_1)$ и все точки,
эквивалентные ей, переходят в бесконечно удаленную точку. По
теореме~\ref{zeroes} можно записать униформизирующую функцию через
нули и полюсы: $$\pi(z)=- \frac{\sigma(z-z_1)\sigma(z-e^{2\pi
i/3}z_1)\sigma(z-e^{4\pi i/3}z_1)}{\sigma(z+z_1)\sigma(z+e^{2\pi
i/3}z_1)\sigma(z+e^{4\pi i/3}z_1)}.
$$
Здесь $\sigma(z)$ --- $\sigma$-функция Вейерштрасса  периодами
$\omega_1=\sqrt{3}$ и $\omega_2=\sqrt{3}e^{i\pi/3}$; при этом,  мы
использовали тот факт, что $\pi(0)=1$.\hskip 1. cm

Обратная функция может быть выражена через интеграл
Кристоффеля-Шварца
$$
z=C_1\int_0^w\frac{dw}{(w^3-1)^{2/3}}+z_1,\quad
C_1=-\,\frac{2\pi}{(\Gamma(1/3))^3}\,=-0.326807...
$$
где $\Gamma$ --- гамма-функция Эйлера.

\subsection{Описание абелева интеграла}

Теперь обсудим, как построить абелев интеграл на римановой
поверхности $R_0$ с нужной асимптотикой (\ref{dif_nat}) при $m=2$.
Заметим, что с помощью дробно-линейного отображения  мы заменили
бесконечно удаленные точки на конечные, В окрестности точки $z_0$
$$
T^{-1}(z)=a+\,\frac{b}{z-z_0}\,\sim\, \frac{b}{z-z_0}\,,
$$
где  $a$ и  $b\neq 0$ -- некоторые константы. Поэтому
$\widetilde{G}=G\circ T^{-1}$ в окрестности точки, лежащей над
$z_0$ на  нулевом листе, имеет асимптотику
\begin{equation}\label{dif_nat1}
\widetilde{G}(z)\sim 2\ln \frac{b}{z-z_0}\,\sim -2\ln (z-z_0).
\end{equation}
В окрестности остальных двух точек, лежащих над $z_0$,
\begin{equation}\label{dif_nat2}
\widetilde{G}(z)\sim -\ln \frac{b}{z-z_0}\,\sim \ln (z-z_0).
\end{equation}

Теперь построим функцию $\widetilde{G}_1$ на комплексной
плоскости, которая дает абелев интеграл на торе
$\mathbb{C}/\mbox{\boldmath$\omega$}$ такой, что
$\widetilde{G}_1=\widetilde{G}\circ \pi$.

В треугольнике $\Delta$ имеется ровно одна точка $\alpha$, которая
при отображении $\pi$ переходит в точку $z_0$. (Напомним, что $z_0$
--- образ бесконечно удаленной точки при дробно-линейном преобразовании $T$, переводящем  исходный треугольник в
правильный.) Заметим, что исходная поверхность  трехлистна, и при
преобразовании $T$ три точки, лежащие над бесконечно  удаленной,
переходят в три точки, лежащие над точкой $z_0$. Тогда точки
$e^{2\pi i/3}\alpha$ и   $e^{-2\pi i/3}\alpha$, которые получаются
поворотом $\alpha$ вокруг нуля, соответствуют двум остальным
точкам на поверхности $R_0$, которые лежат над $z_0$.

Рассмотрим функцию
\begin{equation}\label{harm}
\widetilde{G}_1(z)=-2\ln\sigma(z-\alpha)+\ln\sigma(z-e^{2\pi
i/3}\alpha)+\ln\sigma(z-e^{-2\pi i/3}\alpha).
\end{equation}
В окрестности точек $\alpha$, $e^{2\pi i/3}\alpha$ и $e^{-2\pi
i/3}\alpha$ она имеет нужную асимптотику:
\begin{equation*}\label{dif_nat3}
\widetilde{G}_1(z)\sim \left\{%
\begin{array}{cc}
-2\ln (z-\alpha), &z\to \alpha,   \\[2mm]
\ln (z-e^{\pm2\pi i/3}\alpha),  &z\to e^{\pm2\pi i/3}\alpha.\\
\end{array}%
\right.
\end{equation*}

В силу (\ref{persi1}) и (\ref{persi2})  имеем
\begin{equation}\label{ln_sig}
\ln\sigma(z+\omega_k)-\ln\sigma(z)=\eta_k(z+\omega_k/2)\pm \pi
i,\quad  k=1,2,
\end{equation}
где $\eta_k=2\zeta(\omega_k/2)$ и знак последнего выражения
определяется выбором ветви логарифма.

Покажем, что
\begin{equation}\label{eta12}
\eta_1=2\pi/3,\quad \eta_2=(2\pi/3)e^{-\pi i/3}.
\end{equation}
 Действительно, $\omega_1=\sqrt{3}$, $\omega_2=\sqrt{3}e^{\pi
i/3}$. Кроме того, в силу инвариантности решетки относительно
поворотов на углы $\pi/3$ имеем $\eta_2=\eta_1e^{-\pi i/3}$.
Используя (\ref{etaomega}), получаем
$\eta_1\omega_2-\eta_2\omega_1=\eta_1\sqrt{3}e^{\pi
i/3}-\eta_1e^{-\pi i/3}\sqrt{3}=3\eta_1 i= 2\pi i$, откуда следует
(\ref{eta12}).

Из (\ref{harm}) и (\ref{ln_sig})  находим
$$
\widetilde{G}_1(z+\omega_k)-\widetilde{G}_1(z)=3\alpha\eta_k \quad
(\mbox{\rm mod}\, 2\pi i), \quad k=1,2,
$$
так как $$
\widetilde{G}_1(z+\omega_k)-\widetilde{G}_1(z)=-2(\ln\sigma(z-\alpha+\omega_k)-\ln\sigma(z-\alpha))+(\ln\sigma(z-e^{2\pi
i/3}\alpha+\omega_k)-\ln\sigma(z-e^{2\pi
i/3}\alpha))+$$$$+(\ln\sigma(z-e^{-2\pi
i/3}\alpha+\omega_k)-\ln\sigma(z-e^{-2\pi
i/3}\alpha))=-2\eta_k(z-\alpha+\omega_k/2)+ \eta_k(z-\alpha e^{2\pi
i/3}+\omega_k/2)+$$$$+ \eta_k(z-\alpha e^{-2\pi
i/3}+\omega_k/2)=\eta_k\alpha(2-e^{2\pi i/3}-e^{-2\pi
i/3})=3\alpha\eta_k \quad (\mbox{\rm mod}\, 2\pi i).
$$

Тогда в силу (\ref{per_sqrt2}) и (\ref{eta12}) действительная
часть $u(z)$ функции
$$
\Phi(z)=\widetilde{G}_1(z)-\sqrt{3}\eta_1\overline{\alpha}z
$$
является двоякопериодической гармонической функцией с периодами
$\omega_1$ и $\omega_2$.

Напомним, что нулевой лист поверхности $R_0$ состоит из точек,
которые соответствуют множеству $S_0$ точек $z$ на универсальном
накрытии, для которых
$$u(z)>u(ze^{2\pi i/3}) \quad \mbox{\rm и} \quad u(z)>u(ze^{-2\pi i/3}). $$
Для множества точек $S_2$, соответствующих  второму листу,
$$u(z)<u(ze^{2\pi i/3}) \quad \mbox{\rm и} \quad u(z)<u(ze^{-2\pi i/3}). $$
Наконец, для точек  $S_1$, соответствующих  первому листу,
$$u(ze^{-2\pi i/3})<u(z)<u(ze^{2\pi i/3}) \quad \mbox{\rm либо} \quad u(ze^{2\pi i/3})<u(z)<u(ze^{-2\pi i/3}). $$

Пусть $0\le j$, $k\le2$, $j\neq k$. Обозначим
$$L_{jk}:=\{z:u(ze^{2\pi j i/3})>u(ze^{2\pi k i/3})\},$$
$$\Gamma_{jk}=\{z:u(ze^{2\pi j i/3})=u(ze^{2\pi k i/3})\}.$$
Поскольку группа поворотов на углы, кратные $2\pi k /3$, изоморфна
группа $\mathbb{Z}/3\mathbb{Z}$, мы можем определить $L_{jk}$ и
$\Gamma_{jk}$ для любых целых $j$ и $k$ не сравнимых по модулю 3;
для этого достаточно взять вместо $j$ и $k$ их остатки при делении
на $3$. В этих обозначениях
\begin{equation}\label{sheets}
S_0=L_{01}\cap L_{02},\quad  S_2=L_{10}\cap L_{20},\quad
S_1=(L_{01}\cap L_{20})\cup (L_{02}\cap L_{10}).
\end{equation}
Отметим, что для любых $j\neq k$ три множества $S_{jk}$, $S_{kj}$ и
$\Gamma_{jk}=\Gamma_{kj}$ попарно не пересекаются и дают в
объединении всю комплексную плоскость.

При повороте на угол $2\pi k /3$ множества $L_{jk}$ и
$\Gamma_{jk}$ переходят в множества $L_{j-1,k-1}$ и
$\Gamma_{j-1,k-1}$. Действительно,
$$
e^{2\pi i/3}L_{jk}=e^{2\pi i/3}\{z:u(ze^{2\pi j i/3})>u(ze^{2\pi k
i/3})\}=\{e^{2\pi i/3}z:u(ze^{2\pi j i/3})>u(ze^{2\pi k
i/3})\}=$$$$=\{w:u(we^{2\pi (j-1) i/3})>u(we^{2\pi (k-1)
i/3})\}=L_{j-1,k-1}.
$$
Тогда и для любого $n\in \mathbb{Z}$
\begin{equation}\label{rot}
  e^{2\pi n i/3}L_{jk}=L_{j-n,k-n}.
\end{equation}
С использованием (\ref{rot}) мы можем записать (\ref{sheets}) в
виде
\begin{equation}\label{sh}
S_0=L_{01}\cap e^{-2\pi i/3} L_{10},\quad  S_2=L_{10}\cap e^{-2\pi
i/3}L_{01},\quad S_1=(L_{01}e^{-2\pi i/3}\cap L_{01})\cup
(e^{-2\pi i/3}L_{10}\cap L_{10}).
\end{equation}

Исследуем множества $L_{jk}$. В силу инвариантности решетки
относительно поворотов на углы, кратные $2\pi/3$,  и однородности
функции $\sigma(z)=\sigma(z;\omega_1,\omega_2 )$ как функции $z$ и
периодов $\omega_1$ и $\omega_2$ имеем $$\ln|\sigma(ze^{\pm 2\pi
i/3})|=\ln|\sigma(z)|,$$  следовательно,
$$u(ze^{2\pi i/3})=-2\ln|\sigma(ze^{2\pi
i/3}-\alpha)|+\ln|\sigma(ze^{2\pi i/3}-e^{2\pi
i/3}\alpha)|+\ln|\sigma(ze^{2\pi i/3}-e^{4\pi
i/3}\alpha)|-$$$$-\sqrt{3}\eta_1\Re(\overline{\alpha}e^{2\pi
i/3}z)=-2\ln|\sigma(z-e^{-2\pi
i/3}\alpha)|+\ln|\sigma(z-\alpha)|+\ln|\sigma(z-e^{2\pi
i/3}\alpha)|-\sqrt{3}\eta_1\Re(\overline{\alpha}e^{2\pi i/3}z).$$
 Поэтому неравенство
 $u(z)>u(ze^{2\pi i/3})$ эквивалентно неравенству
$$
\ln|\sigma(z-\alpha)|-\ln|\sigma(z-e^{-2\pi
i/3}\alpha)|<(\sqrt{3}/3)\eta_1\Re(\overline{\alpha}(e^{2\pi
i/3}-1)z),
$$
аналогично неравенство
 $u(z)>u(ze^{4\pi i/3})$ эквивалентно неравенству
$$
\ln|\sigma(z-\alpha)|-\ln|\sigma(z-e^{2\pi
i/3}\alpha)|<(\sqrt{3}/3)\eta_1\Re(\overline{\alpha}(e^{-2\pi
i/3}-1)z),
$$
а неравенство
 $u(ze^{2\pi i/3})>u(ze^{-2\pi i/3})$ --- неравенству
\begin{equation*}\label{}
\ln|\sigma(z-e^{-2\pi i/3}\alpha)|-\ln|\sigma(z-e^{2\pi
i/3}\alpha)|<(\sqrt{3}/3)\eta_1\Re(\overline{\alpha}(e^{-2\pi
i/3}-e^{2\pi i/3})z).
\end{equation*}

\subsection{Случай вещественного $\alpha$}

В дальнейшем ограничимся исследованием случая, когда $\alpha$
вещественно.

В силу (\ref{rot}) достаточно изучить геометрическую структуру
множества $L_{12}$, поскольку остальные множества $L_{jk}$
получаются из него поворотами на углы, кратные $2\pi /3$, и
операцией дополнения (с отбрасыванием границы). При вещественном
$\alpha$ множество $\Gamma_{12}$ является границей $L_{12}$, оно
симметрично относительно оси абсцисс, содержит эту ось и поэтому его
легче описывать. Поэтому начнем с описания $\Gamma_{12}$.

Рассмотрим гармоническую функцию
$$
g(z)=g(z;\alpha)=\ln|\sigma(z-e^{-2\pi
i/3}\alpha)|-\ln|\sigma(z-e^{2\pi
i/3}\alpha)|-(\sqrt{3}/3)\eta_1\Re(\overline{\alpha}(e^{-2\pi
i/3}-e^{2\pi i/3})z),
$$
нули которой дают множество решений уравнения
$$
u(ze^{2\pi i/3})=u(ze^{-2\pi i/3}),
$$
и исследуем ее нулевое множество.

Пусть $\alpha$ -- вещественное число, $\alpha\in (0,\sqrt{3}/2)$.
Тогда  уравнение $g(z)=0$ эквивалентно уравнению
\begin{equation}\label{set}
\ln|\sigma(z-e^{-2\pi i/3}\alpha)|-\ln|\sigma(z-e^{2\pi
i/3}\alpha)|-\eta_1\alpha\Im z=0,
\end{equation}
множество решений которого содержит вещественную ось. Критические
точки функции $g$, т.~е. точки, где ее градиент обращается в нуль,
удовлетворяют уравнению
\begin{equation}\label{crit}
\zeta(z-e^{-2\pi i/3}\alpha)-\zeta(z-e^{2\pi i/3}\alpha)+\eta_1
\alpha i=0,
\end{equation}
Поскольку функция $\zeta(z-e^{-2\pi i/3}\alpha)-\zeta(z-e^{2\pi
i/3}\alpha)$ является двоякопериодической и в любом
параллелограмме периодов имеет два полюса, она двулистна в
параллелограмме, поэтому принимает значение $(-\eta_1 \alpha i)$
ровно в двух точках (с учетом кратности). Следовательно, уравнение
(\ref{crit}) имеет ровно два решения в параллелограмме периодов,
т.~е. критических точек  функции $g$ там тоже две.

Исследуем, при каких $\alpha$ эти критические точки лежат на
вещественной оси. Если $z=x$ вещественно, то в силу вещественности
$\alpha$ и симметричности $\zeta(z)$:
$$
\zeta(\overline{z})=\overline{\zeta(z)},
$$
уравнение (\ref{crit})
эквивалентно равенству
$$
\varphi(x):=\Im \zeta(x-e^{2\pi i/3}\alpha)-(\eta_1/2){\alpha}=0.
$$
Функция $\varphi$ в силу (\ref{periodzeta}) является периодической
на вещественной оси с периодом $\sqrt{3}$. Поскольку
$$\varphi'(x)=-\Im \mathfrak{P}(x-e^{2\pi i/3}\alpha),$$ видим,
что $\varphi'(x)=0$ в точках, где $$\mathfrak{P}(x-e^{2\pi
i/3}\alpha)=\mathfrak{P}(x+\alpha/2-i\alpha \sqrt{3}/2)$$
принимает вещественные значения. Значит, $x+\alpha/2=0$ или
$x+\alpha/2=\sqrt{3}/2$ (\mbox{\rm mod}\,
\mbox{\boldmath$\omega$}).

Рассмотрим функцию $\varphi$ на отрезке
$\Delta:=[-\alpha/2,-\alpha/2+\sqrt{3}]$ шириной, равной периоду
$\omega_1$. Нетрудно видеть, что на концах отрезка $\Delta$ функция
$\varphi$ принимает положительные значения, так как
$$
\varphi(-\alpha/2)=\Im
\zeta(-i\alpha\sqrt{3}/2)-(\eta_1/2){\alpha}>\Im
\zeta(-i)-(\eta_1/2){\alpha}=(\eta_1/2)(2/\sqrt{3}-\alpha)>0.$$
Здесь мы воспользовались тем, что функция $\zeta(z)$ отображает
отрезок мнимой оси с концами в точках $0$ и $(-i)$ на луч, идущий
из точки $\zeta(-i)=i\eta_1/\sqrt{3}$ вверх (см.
замечание~\ref{rect3}). Функция $\varphi'(x)=-\Im
\mathfrak{P}(x-e^{2\pi i/3}\alpha)$ меняет знак с <<$+$>> на
<<$-$>> при переходе через точку $x=-\alpha/2+\sqrt{3}/2$, поэтому
на концах отрезка $\Delta $ функция $\varphi$ принимает
максимальные значения, а $x=-\alpha/2+\sqrt{3}/2$
--- это ее точка минимума. Значение функции $\varphi$
в этой точке равно
$$\psi(\alpha):=\varphi(-\alpha/2+\sqrt{3}/2)=\Im \zeta((\sqrt{3}/2)(1-i\alpha))-(\eta_1/2){\alpha}.$$

Из анализа функции $\mathfrak{P}$ на отрезке с концами
$\sqrt{3}/2$ и $\sqrt{3}/2-i3/2$ видно (см.
замечание~\ref{rect3}), что производная
$$
\psi'(\alpha)=(\sqrt{3}/2)\mathfrak{P}((\sqrt{3}/2)(1-i\alpha))-\eta_1/2
$$
строго монотонно убывает на отрезке $[0,\sqrt{3}]$ от значения
$$\psi'(0)=(\sqrt{3}/2)\mathfrak{P}(\sqrt{3}/2))-\eta_1=0.655\ldots>0$$ до
$\psi'(\sqrt{3})=-\infty$.  Следовательно, функция $\psi$ на
отрезке $[0,\sqrt{3}]$ сначала строго монотонно возрастает, а
потом строго убывает. Поскольку $\psi(0)=\Im\zeta(\sqrt{3}/2)=0$,
делаем вывод, что на первом участке она положительна, а на втором
имеет не более одного нуля. Так как
$$\psi(\sqrt{3}/3)=\Im
\zeta((\sqrt{3}/2)(1-i\sqrt{3}/3))-\eta_1\sqrt{3}/6=\Im
\zeta(\sqrt{3}/2-i/2)-\eta_1\sqrt{3}/6=0,$$ делаем вывод, что на
отрезке $[0,\sqrt{3}/2]$ функция $\psi$ имеет единственный нуль в
точке $\sqrt{3}/3$. Подведем итоги в следующем утверждении.

\begin{theorem}\label{cr_points}
$1)$ Если $0<\alpha<\sqrt{3}/3$, то критических точек на
вещественной оси у функции $u(z)$ нет. Этот случай соответствует
равнобедренному треугольнику $\Delta(a_1,a_2,a_3)$ с углом при
вершине, меньшем $\pi/3$.

$2)$ Если $\sqrt{3}/3<\alpha<\sqrt{3}/2$, то у функции $u(z)$
имеется ровно две критические точки на любом отрезке, по ширине
равном периоду $\omega_1=\sqrt{3}$. Этот случай соответствует
равнобедренному треугольнику $\Delta(a_1,a_2,a_3)$ с углом при
вершине, меньшем $\pi/3$.

$3)$ Наконец, если
$$\alpha=\sqrt{3}/3=0.57735\ldots,$$ то критическая точка ---
двойная. Это --- случай правильного треугольника
$\Delta(a_1,a_2,a_3)$.
\end{theorem}

Можно показать, что при $0<\alpha<\sqrt{3}/3$ топологическая
структура множеств  $\Gamma_{jk}$ одинакова, то же касается и случая
$\sqrt{3}/3<\alpha<\sqrt{3}/2$.

На рис.~7 приведена структура  $\Gamma_{12}$ для
$\alpha=0.5<\sqrt{3}/3$ (случай 1) теоремы~\ref{cr_points}) и
$\alpha=0.6>\sqrt{3}/3$ (случай~2) теоремы~\ref{cr_points}). На
рис.~7~а видно, что в случае 1) множество $\Gamma_{12}$ состоит из
непересекающихся линий
--- прямых, параллельных оси абсцисс, и простых периодических
кривых, не пересекающих эти прямые.  На рис.~7~б линии уровня
состоят из прямых  $\Im z=(3/2)m$, $m\in \mathbb{Z}$, и семейства
кривых пересекающих  ортогонально прямые в критических точках,
описанных выше, и эквивалентных с ними по модулю решетки. Заметим
также, что в случае $\alpha=\sqrt{3}/3$ получается замощение
плоскости правильными треугольниками.

На рис.~8 для тех же значений параметра $\alpha$ изображены листы
$S_j$, $0\le j\le 2$. Видно, что нулевой лист является всегда
связным множеством. Лист  $S_2$ является связным в случае 1) и
состоит из трех компонент связности в случае 2). Лист $S_1$
всегда несвязен, и состоит либо из четырех, либо из шести
компонент.

Более подробный анализ показывает, что в случае
$0<\alpha<\sqrt{3}/3 $ множество $\Gamma_{12}$ состоит из прямых
$\Im z=(3/2)m$, $m\in \mathbb{Z}$, и множества кривых $\gamma_m$,
$m\in \mathbb{Z}$, обладающих свойством: кривая $\gamma_0$ лежит в
полосе $0<\Im z<3/2$ и является графиком некоторой гладкой функции
$y=F(x)$ с периодом $\sqrt{3}$; функция $F$ достигает
максимального значения в точках $-\alpha/2+\sqrt{3} n$,  и
минимального значения в точках $-\alpha/2+\sqrt{3}(n+1/2)$, $n\in
\mathbb{Z}$. В остальных точках прямой $\psi'(x)\neq 0$.

В заключение выражаем благодарность члену-корр. РАН проф.
А.~А.~Аптекареву, обратившему наше внимание на задачу, описанную в
третьей части настоящей работы, и проф. С.~К.~Водопьянову за
приглашение прочесть курс лекций в Региональном математическом
центре при Новосибирском государственном университете и
прекрасные условия, созданные для написания этой статьи.

\hskip -0.7 cm\includegraphics[width=6.8 in,%
]{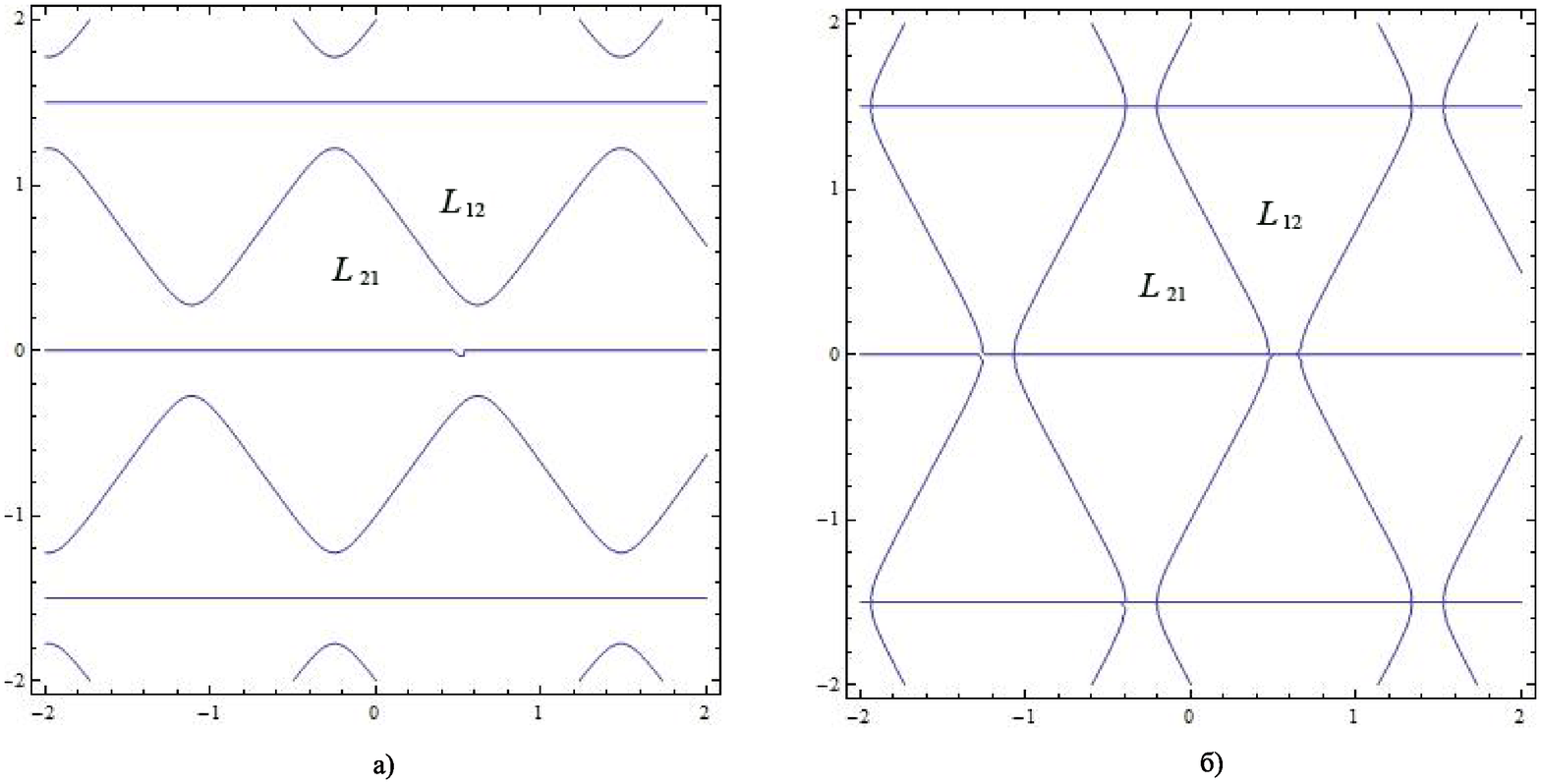} \vskip -3 cm\begin{center} Рис. 7\end{center} \vskip
1 cm

\includegraphics[width=6.4 in,%
]{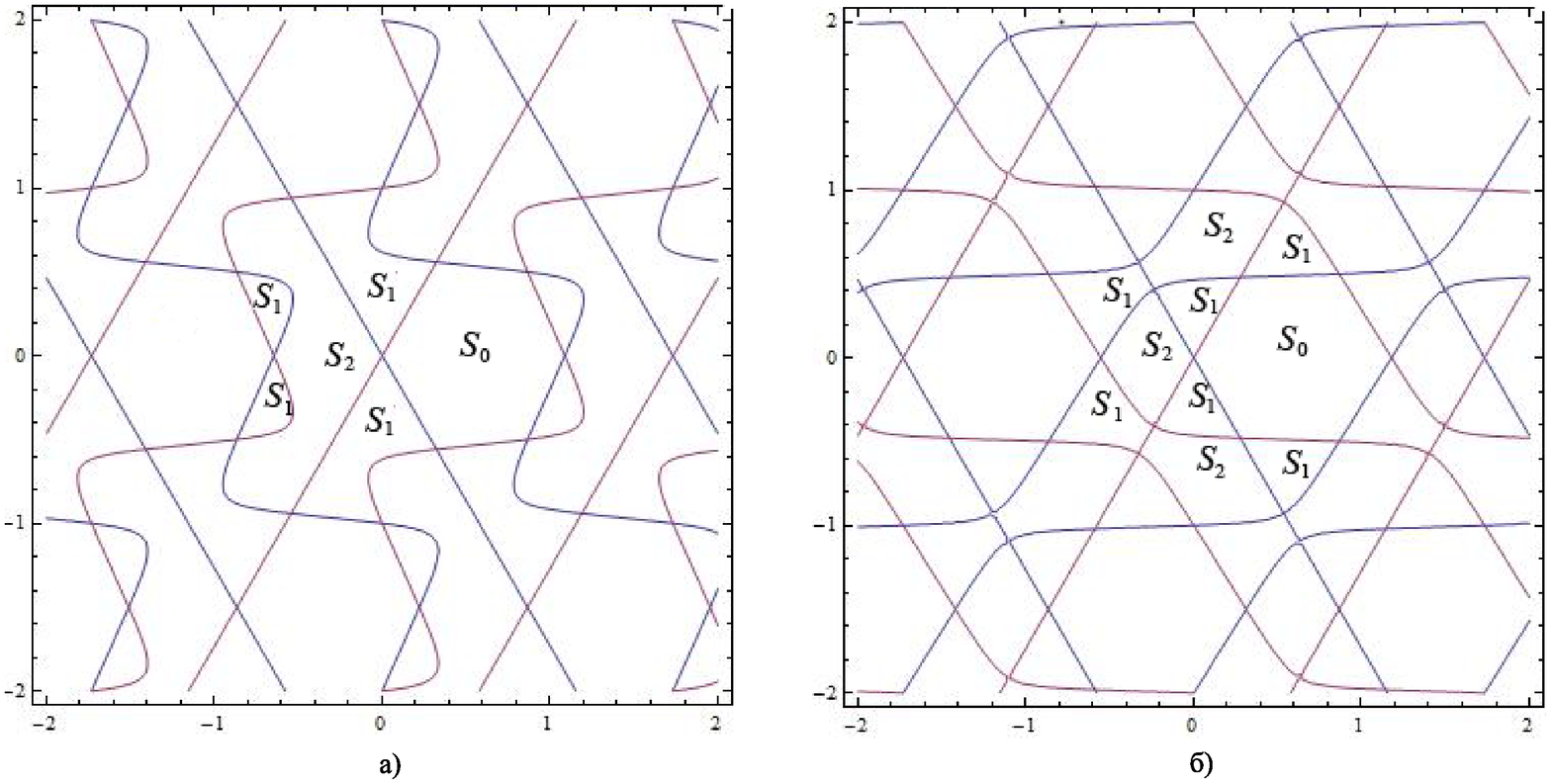}\
\begin{center} Рис. 8\end{center}

\end{document}